\documentclass[11pt]{extarticle}
\usepackage{amsmath,amsthm,amssymb,amsbsy,amsopn,mathtools,comment,mathptmx,bm,bbm,fancyhdr,lineno}

\def\esssup{\mathop{\rm esssup}}

\usepackage[ruled,vlined,algo2e]{algorithm2e}
\usepackage{textcomp,breqn,mathrsfs}
\usepackage{graphicx,xcolor,float,subcaption,caption,algorithm}
\usepackage[noend]{algpseudocode}
\usepackage[bookmarks=true,bookmarksnumbered=true,bookmarkstype=toc]{hyperref}
\usepackage[toc,page]{appendix}

\def\argmin{\mathop{\rm argmin}}

\newtheorem{Def}{Definition}[section]
\newtheorem{theorem}[Def]{Theorem}

\newtheorem{assumption}[Def]{Assumption}

\theoremstyle{definition}

\makeatletter
\@addtoreset{equation}{section}
\makeatother

\title{Numerical methods for backward stochastic differential equations: A survey\footnotetext[0]{
The authors are grateful to the anonymous referees and Bernhard Hientzsch for their valuable feedback that helped improve the quality of this manuscript.
This work was initiated while JC and RK were based in School of Mathematics and Statistics at the University of Sydney, Australia. 
The opinions expressed by the authors are solely their own and do not reflect those of Appian Corporation, Bank of Japan or their related entities.
This work was partially supported by JSPS Grants-in-Aid for Scientific Research (Grant Numbers 20K22301 and 21K03347) and by JST PRESTO (Grant Number JPMJPR2029).}}
\author{\sc Jared Chessari\footnote{
Email: jared.chessari@appian.com. Appian Corporation, Australia.},\, Reiichiro Kawai\footnote{Corresponding author. 
Email: raykawai@g.ecc.u-tokyo.ac.jp. Affiliations: Graduate School of Arts and Sciences / Mathematics and Informatics Center, the University of Tokyo, Japan.},\, Yuji Shinozaki\footnote{Email: sinozaki.yuji@gmail.com.  Bank of Japan, and Tokyo Institute of Technology, Japan.} \, and Toshihiro Yamada\footnote{Email: toshihiro.yamada@r.hit-u.ac.jp. Affiliations: Graduate School of Economics, Hitotsubashi University and Japan Science and Technology Agency (JST), Japan.}}
\date{}
\begin{document}


\maketitle


\begin{abstract}
\noindent 
Backward Stochastic Differential Equations (BSDEs) have been widely employed in various areas of social and natural sciences, such as the pricing and hedging of financial derivatives, stochastic optimal control problems, optimal stopping problems and gene expression. 
Most BSDEs cannot be solved analytically and thus numerical methods must be applied to approximate their solutions.
There have been a variety of numerical methods proposed over the past few decades as well as many more currently being developed.
For the most part, they exist in a complex and scattered manner with each requiring a variety of 
assumptions and conditions.
The aim of the present work is thus to systematically survey various numerical methods for BSDEs, and in particular, compare and categorize them, for further developments and improvements.
To achieve this goal, we focus primarily on the core features of each method based on an extensive collection of 333 references: the main assumptions, the numerical algorithm itself, key convergence properties and advantages and disadvantages, to provide an up-to-date coverage of numerical methods for BSDEs, with insightful summaries of each and a useful comparison and categorization. 
\vspace{0.3em}


\end{abstract}

\tableofcontents


\section{Introduction}\label{section introduction}

A Backward Stochastic Differential Equation (BSDE) was first introduced by Jean-Michel Bismut in 1973 \cite{Bismut1973Conjugate}.
The paper used a linear BSDE as an equation for the adjoint process in the stochastic version of the Pontryagin maximum principle. 
BSDEs were generalized by Pardoux and Peng in 1990 \cite{Adapted_PaPe} to a general non-linear BSDE of the following form:
    \begin{equation*}
    Y_t = \xi_T + \int_t^T f(s, Y_s, Z_s)\,ds - \int_t^T Z_s\,dW_s, \quad t \in [0,T],
    \end{equation*}
    where questions regarding the existence and uniqueness of the solution $(Y_t, Z_t)_{t \in [0,T]}$ were addressed. We also mention here the closely related Forward BSDEs (FBSDEs) given as follows:
	\begin{equation*} 
    \begin{cases} 
         Y_t = \Phi(X_T) + \int_t^T f(s, X_s, Y_s, Z_s)ds - \int_t^T Z_sdW_s, \\
          X_t = {\bf x} + \int_0^t \mu(s,X_s)ds + \int_0^t\sigma(s,X_s)dW_s,
         \end{cases}
	\end{equation*}
for $t\in [0,T]$.
	They are essentially a BSDE coupled with a (forward) Stochastic Differential Equation (SDE), in the sense that the driver can depend on the solution to the SDE, and the terminal condition of the BSDE is a function of the terminal value of the SDE. We also look in detail at FBSDEs, as these types of BSDEs appear widely in literature (especially in regards to numerical methods). The solution of an FBSDE can be formulated in terms of the solution of a semilinear parabolic PDE, and this equivalence is exploited in many numerical methods. For the remainder of the introduction when we mention BSDEs, this can refer to BSDEs or FBSDEs. In following sections, however, we will be strict on which type of BSDE we are referring to. 
	
BSDEs have been used widely in a number of areas of social and natural sciences.
The problem of pricing and hedging a European option can be formulated in terms of a BSDE.
In fact, any pricing problem using a replication argument can be written in terms of a BSDE \cite{BSDEFin_Kar_Pen,  https://doi.org/10.48550/arxiv.1911.12231}.
A BSDE can also take into account portfolio constraints in pricing problems \cite{Hedge_constraints1, Hedge_constraints2, Hedge_constraints3, Hedge_constraints4}.
	Thus, BSDEs have often been employed for the valuation of many financial derivatives in both complete and incomplete or constrained markets, including European and American options \cite{BSDE_American, BSDE_American2}.
	BSDEs can provide necessary and sufficient conditions for optimality \cite{BSDE_opti1, BSDE_opti2} and enhance implementability \cite{doi:10.1137/17M1123559, hawkins2021value} in optimal control problems such as utility maximization control problems with constraints and risk-sensitive control problems \cite{BSDE_risksens1, https://doi.org/10.48550/arxiv.1911.12231, BSDE_utility1, BSDE_utility2, BSDE_utility3, BSDE_risksens2}. 
BSDEs have also been tied to nonlinear expectations and thus have been used to construct the risk measures \cite{BSDE_riskmeas}. The more complex reflected BSDEs \cite{BSDE_American} and doubly reflected BSDEs \cite{Double_refl1, Double_refl2} have been developed in connection with zero-sum Dynkin games, optimal stopping problems, recallable options, mixed differential games, and mixed stochastic control \cite{Complex_Egs1, Complex_Egs2}. 
In addition, BSDEs have found applications in filtering \cite{CiCP-27-589, Bao_2019, doi:10.1142/S0219530520400102} and physics \cite{PhysRevE.96.042123, PhysRevE.95.032418}.
More examples of applications in these areas can be found in \cite{crisan_manolarakis_2010, BSDE_Applications, BSDEFin_Kar_Pen} and image processing \cite{borkowski_janczakborkowska_2021}.
	
Due to the complex nature of BSDEs, it is rarely possible to find an analytical solution.
Thus, one must often resort to numerical methods in order to solve these equations.
There have been a variety of numerical methods proposed over the past few decades and also many more currently in progress. For the most part they exist in a complex and scattered manner, with each requiring a variety of 
assumptions and conditions.
The ultimate goal of this survey is to further facilitate the research activity on numerical methods for BSDEs by providing an up-to-date overview of the different types in an organized structure.
Before providing the systematic survey of numerical methods, we will first review some basic theory regarding BSDEs (Section \ref{section preliminaries}) in order to better prepare for the goal.
Broadly speaking, we will then look at the three main classes of numerical methods, namely backward (Section \ref{section backward numerical methods}), forward (Section \ref{section forward numerical methods}), and deep learning based methods (Section \ref{section deep learning}).
By a backward method, we mean a numerical method which works backwards in time and requires the computation of conditional expectations (Section \ref{section computation of conditional expectations}), whereas for a forward method, we refer to a method that does not inherently work backward in time so as to (originally, at least) avoid the computation of conditional expectations.
In sharp contrast to backward and forward methods, the method of deep learning (Section \ref{section deep learning}) is quite distinctive in its rather mixed structure and effectiveness in high-dimensional problems.

We then summarize in Section \ref{section discussion} the major components of each category
to discuss the power and limitations of various classes of numerical methods in collective comparison.
We next provide in Section \ref{section BSDEs with nonstandard features} a further survey on numerical methods for BSDEs with a variety of nonstandard features, such as coupled FBSDEs (Section \ref{subsection coupled FBSDEs}), reflected BSDEs (Section \ref{subsection reflected BSDEs numerical methods}), BSDEs with jumps (Section \ref{subsection BSDEs with jumps}), non-Lipschtiz BSDEs (Section \ref{subsection Lipschitz BSDEs}), quadratic BSDEs (Section \ref{subsection qBSDEs}), second-order BSDEs (Section \ref{section 2bsde}), McKean-Vlasov BSDEs (Section \ref{subsection mckean-vlasov}) and backward stochastic partial differential equations (Section \ref{subsection BSPDEs}), each of which has attracted increasing interest for addressing various emerging problems in social and natural sciences and nonlinear partial differential equations.
Finally, in Section \ref{section concluding remarks}, we summarize our discourse and highlight some future directions.

In order to achieve our main goal of describing and discussing the main idea of each category, we will only focus on one or two relevant representative methods in each subsection, followed by an overview of other methods.
To avoid overloading the paper with non-essential technical details, we omit the rather lengthy technical intricacies of the methods in most instances.
In particular, the algorithms we provide in each section are described without going into detail and are presented for illustrative purposes only. 

\section{Preliminaries}\label{section preliminaries}

The aim of this section is to review the basic theory of BSDEs (Section \ref{subsection BSDEs in preliminaries}) and FBSDEs (Section \ref{subsection FBSDEs in preliminaries}) in brief. 
In particular, we give the standard sufficient conditions that ensure the existence and uniqueness of their solutions.
Let $(\Omega, \mathcal{F}, \mathbb{F}, \mathbb{P})$ be a given filtered probability space which satisfies the usual conditions of completeness and right continuity.
On this space, we define the $d$-dimensional standard Brownian motion $(W_t)_{t \in [0,T]}$, whose natural filtration, augmented by the class of $\mathbb{P}$-null sets of $\mathcal{F}_T$, is $\mathbb{F}$ = $(\mathcal{F}_t)_{t \in [0,T]}$.
We reserve $T$ for a strictly positive constant which indicates a fixed terminal time of the interval $[0,T]$, and denote by $\pi_n:=\{0=:t_0<\cdots<t_k<t_{k+1}<\cdots<t_n:=T\}$ an arbitrary yet fixed partition of the interval $[0,T]$ for all $n\in\mathbb{N}$, with $|\pi_n|:=\max_{k\in \{0,\cdots,n-1\}} (t_{k+1}-t_k)$.
For the sake of brevity, we often suppress the subscript $n$ from the notation $\pi_n$.
We define the following probabilistic spaces and sets:
    \begin{itemize}
    \setlength{\parskip}{0cm}
    \setlength{\itemsep}{0cm}
        \item $\mathcal{P}$: the $\sigma$-field of predictable sets in $\Omega \times [0,T]$.
        \item $\mathcal{L}_m^2(\mathcal{F}_t)$: the set of all $\mathbb{R}^m$-valued random vectors $X$ that are $\mathcal{F}_t$-measurable and square-integrable, that is, $\mathbb{E}[\|X\|^2] < \infty$.
        \item $\mathcal{S}_m^2(0,T)$: the space of all continuous predictable processes $X: \Omega \times [0,T] \rightarrow \mathbb{R}^m$, satisfying $\mathbb{E}[\sup\nolimits_{t \in [0,T]}\|X_t\|^2]<\infty$.
        \item $\mathcal{H}_m^2(0,T)$: the space of all predictable processes $X: \Omega \times [0,T] \rightarrow \mathbb{R}^m$, which satisfy $\|X\|_0^{2} \coloneqq \mathbb{E} [\int_0^T \|X_t\|^2dt] < \infty$.
    \end{itemize}
    
     We reserve $C$ for a positive constant, whose value changes from line to line, depending on the context at hand. Let $|\cdot|$ and $\|\cdot\|$ denote the magnitude and Euclidean norm respectively, where the latter is understood to be a suitable matrix norm in the context of matrices.
    We denote $D^{\otimes 2}:=DD^{\top}$ and by $\mathcal{B}(D)$ the Borel $\sigma$-field of a set $D$.
    We define the following function spaces which appear throughout:
 \begin{itemize}
    \setlength{\parskip}{0cm}
    \setlength{\itemsep}{0cm}
        \item $\mathcal{C}_b^{k,l}$: set of continuously differentiable and bounded functions which are $k$ times continuously differentiable in their first coordinate and $l$ times in their second, with bounded partial derivatives up to order $l$.
        \item $\mathcal{C}^{k}_p$: set of $\mathcal{C}^{k-1}$ functions with piecewise continuous $k$-th derivative.
        \item $\mathcal{C}^{k + \alpha}$: set of $\mathcal{C}^k$ functions whose $k$-th derivative is H\"{o}lder continuous of order $\alpha \in [0,1]$.
		\end{itemize}

\subsection{Backward stochastic differential equations}
\label{subsection BSDEs in preliminaries}

We aim here to review the basics of the following BSDE:
    \begin{equation}
    \begin{cases} \label{eq:BSDE_Diff}
         -dY_t = f(t, Y_t, Z_t)dt - Z_tdW_t, & t \in [0,T], \\
        Y_T = \xi, \\
    \end{cases}
	\end{equation}
where $(Y,Z)_{t \in [0,T]}$ takes values in $\mathbb{R}^m \times \mathbb{R}^{m \times d}$. 
The function $f : \Omega \times [0,T] \times \mathbb{R}^m \times \mathbb{R}^{m \times d} \rightarrow \mathbb{R}^m$ is called the generator or driver. 
The random variable $\xi = Y_T$ is the terminal condition for the process $Y$ and the pair $(f, \xi)$ are called the parameters of the BSDE.
A pair $(Y,Z)$ of stochastic processes is a (not necessarily unique) solution to the BSDE \eqref{eq:BSDE_Diff} if
	\begin{itemize}
	\setlength{\parskip}{0cm}
    \setlength{\itemsep}{0cm}
	    \item $(Y_t)_{t \in [0,T]}$ is an $\mathbb{F}$-adapted continuous process in $\mathbb{R}^m$.
	    \item $(Z_t)_{t \in [0,T]}$ is an $\mathbb{F}$-predictable process in $\mathbb{R}^{m \times d}$ satisfying $\mathbb{P}(\int_0^T \|Z_t\|^2dt < \infty) = 1$.
	    \item $Y_t = \xi+ \int_t^T f(s, Y_s, Z_s)\,ds - \int_t^T Z_s\,dW_s$ for all $t \in [0,T]$, $a.s.$
	\end{itemize}
We will only consider solutions $(Y,Z)$ in the space $\mathcal{S}^2_m(0,T) \times \mathcal{H}_{m \times d}^2(0,T)$.
The following theorem \cite{Adapted_PaPe, zhang_2017} identifies the standard sufficient conditions imposed to ensure the existence and uniqueness of a solution to the BSDE \eqref{eq:BSDE_Diff} in $\mathcal{S}^2_m(0,T) \times \mathcal{H}_{m \times d}^2(0,T)$.

\begin{theorem}[\bf Existence and uniqueness of the solution to a BSDE] \label{Pa_Pe_Unique}
Consider the BSDE given by \eqref{eq:BSDE_Diff} with parameters $(f,\xi)$, satisfying:
\begin{itemize}
    \setlength{\parskip}{0cm}
    \setlength{\itemsep}{0cm}
    \item $f : \Omega \times [0,T] \times \mathbb{R}^m \times \mathbb{R}^{m \times d} \rightarrow \mathbb{R}^m$ is $\mathcal{P} \otimes \mathcal{B}(\mathbb{R}^m) \otimes \mathcal{B}(\mathbb{R}^{m 
    \times d})$-measurable.
    \item The process $(f(t,0,0))_{t \in [0,T]}$ belongs to $\mathcal{H}_m^2(0,T)$.
    \item The driver is uniformly Lipschitz continuous with respect to $({\bf y},{\bf z})$, that is, there exists a constant $C>0$ such that for all ${\bf y}_1, {\bf y}_2 \in \mathbb{R}^m$ and ${\bf z}_1, {\bf z}_2 \in \mathbb{R}^{m \times d}$, $\|f(\omega,t,{\bf y}_1,{\bf z}_1) - f(\omega, t, {\bf y}_2, {\bf z}_2)\| \leq C(\|{\bf y}_1-{\bf y}_2\| + \|{\bf z}_1 - {\bf z}_2\|),$ $dt \otimes d\mathbb{P}$-$a.e.$
    \item $\xi \in \mathcal{L}_m^2(\mathcal{F}_T)$.
\end{itemize}
Then, there exists a unique solution $(Y,Z) \in \mathcal{S}_m^2(0,T) \times \mathcal{H}_{m \times d}^2(0,T)$ which solves the BSDE \eqref{eq:BSDE_Diff}.
\end{theorem}

	Before moving on, we briefly present some intuition regarding the solution $(Y,Z)$ of the BSDE \eqref{eq:BSDE_Diff}.
	Consider the simple case when $f \equiv 0$ in \eqref{eq:BSDE_Diff} giving us the following BSDE:
	\begin{equation}\label{eq:BSDE_Diff_simple}
    \begin{cases}    
        dY_t = Z_t dW_t , \quad t \in [0,T], \\
        Y_T = \xi.
    \end{cases}
	\end{equation}
	Clearly, the degenerate driver $f \equiv 0$ satisfies the sufficient assumptions for uniqueness and existence.
	Hence, under the assumption that $\xi \in \mathcal{L}_m^2(\mathcal{F}_T)$, it holds by Theorem \ref{Pa_Pe_Unique} that the BSDE with null generator \eqref{eq:BSDE_Diff_simple} admits a unique solution $(Y,Z) \in \mathcal{S}^2_m(0,T) \times \mathcal{H}_{m \times d}^2(0,T)$.
	If $(Y,Z)$ is the unique solution in $\mathcal{S}^2_m(0,T) \times \mathcal{H}_{m \times d}^2(0,T)$ of the BSDE with null driver \eqref{eq:BSDE_Diff_simple}, then $Y$ satisfies 
        \begin{equation*}
            Y_t = \mathbb{E}\left[Y_T | \mathcal{F}_t\right] = \mathbb{E}\left[\xi | \mathcal{F}_t\right],
        \end{equation*}
        and thus $Y$ is an $(\mathcal{F}_t)_{t\in[0,T]}$-martingale.
	Moreover, based on the formulation \eqref{eq:BSDE_Diff_simple}, the martingale $Y$ can be written as follows:
    \[
          Y_t = \mathbb{E}\left[\xi | \mathcal{F}_t\right] = Y_0 + \int_0^t Z_sdW_s 
          = \mathbb{E}\left[\xi | \mathcal{F}_0\right] + \int_0^t Z_sdW_s 
          = \mathbb{E}\left[\xi\right] + \int_0^t Z_sdW_s.
    \]
As the process $Y$ is a martingale, it is certainly a local martingale.
Thus, by the martingale representation theorem, there exists a unique process $\gamma \in \mathcal{H}_{m \times d}^2(0,T)$ satisfying $Y_t = Y_0 + \int_0^t \gamma_s dW_s$ for all $t\in [0,T]$.
Hence, we obtain
\[
         \xi = Y_T = Y_0 + \int_0^T \gamma_sdW_s = Y_0 + \int_0^t \gamma_sdW_s+ \int_t^T \gamma_sdW_s = Y_t + \int_t^T \gamma_sdW_s,
\]
which shows that $(Y,\gamma)$ solves the BSDE \eqref{eq:BSDE_Diff_simple}.
By the uniqueness of the solution to the BSDE, we get $Z = \gamma$.
For this simple case, we see that an adapted solution to the BSDE can only be given by a pair, that is, the $Z$ component is needed to ensure that the process $Y$ is adapted.
In a sense, the $Z$ component ``steers'' the system and is thus called the control process.
One cannot simply revert time as in deterministic ODEs, as the filtration can only go in one direction. 

To better illustrate the interpretations of $Y$, $Z$ and $\xi$, we provide a simple example in a tangible problem in finance.
Consider a financial market with one risky asset $S$, whose price at time t, $S_t$, follows the following SDE:
    \begin{equation}\label{basic sde finance}
        dS_t = S_t \mu_t dt + S_t \sigma_t dW_t.
    \end{equation}
     A trader can either invest in the risky asset $S$ or borrow/invest money at an instantaneous risk free interest rate, denoted by $r$.
    Here, we assume that $\mu$, $\sigma$ and $r$ are bounded and predictable processes.
    If the amount of money invested in $S$ at time $t$ is $\phi_t$ and the total wealth of the trade is $Y$, then the magnitude $|Y_t-\phi_t|$ represents the amount of money that is borrowed (if $Y_t-\phi_t$ is negative) or invested (if $Y_t-\phi_t$ is positive) at time $t$.
    The wealth process $Y$ can thus be shown to follow the dynamics:
    \begin{equation}\label{wealth process finance}
        dY_t = \frac{\phi_t}{S_t} dS_t + r_t (Y_t - \phi_t) dt = \left(\phi_t(\mu_t - r_t) +r_t Y_t\right)dt + \phi_t \sigma_t dW_t.
    \end{equation}
    Consider a European option with payoff at time $T$ given by a random variable $\xi \in \mathcal{L}_m^2(\mathcal{F}_T)$. We note that the payoff for a European option $\xi$ will usually be a function of the asset price $S$ at time $T$, that is $\xi = \Phi(S_T)$.
Although this would move us to the context of a FBSDE (Section \ref{subsection FBSDEs in preliminaries}), we ignore this fact here for illustrative purposes.
    A trader who wants to sell this option must identify the minimal initial amount of capital $Y_0$ in such a way that the payoff $\xi= \Phi(S_T)$ can be replicated.
    If a process $\phi$ is found such that $Y_T = \xi$, then this initial amount is $Y_0$. In other words, we look for a couple $(Y,\phi)$ such that
    \begin{equation*}
        Y_t = \xi - \int_t^T \left(\phi_s(\mu_s - r_s) +r_s Y_s\right)ds - \int_t^T \phi_s \sigma_s dW_s.
    \end{equation*}
    If there exists a predictable process $\pi$ such that $(\mu - r) = \sigma \pi$, then by setting $Z= \phi \sigma$, the above equation becomes
    \[
        Y_t = \xi - \int_t^T \left(Z_s \pi_s + r_s Y_s\right)ds - \int_t^T Z_s dW_s.
    \]
    Hence, the problem of pricing and hedging is reduced to finding a solution to the above (linear) BSDE, where $Y$ is the wealth process (which is equivalent to the price of the option), $\xi= \Phi(S_T)$ is the payoff of a European option and $Z$ is proportional to the volatility in the model and the amount invested in the risky asset $S$.

\color{black}
Another typical example in quantitative finance, that is more general than the simple European case and which justifies the use of nonlinear BSDEs, is the so-called the differential rates problem (for instance, \cite{https://doi.org/10.48550/arxiv.2006.07635} among many others). 
Consider again, a financial market with one risky asset $S$ given by the SDE \eqref{basic sde finance},
with $\phi_t$ representing the amount of money invested in $S$ at time $t$.
We thus have that the dynamics of the total wealth $Y$ are again given by \eqref{wealth process finance}.
Let a lending interest rate $r$ be applied to the lending case $Y_t-\phi_t\geq 0$ and let a borrowing interest rate $R$ such that $R>r$ be applied to the borrowing case $Y_t-\phi_t \leq 0$. 
In a similar manner to the European case above, if a process $\phi$ is found such that $Y_T = \xi$ with target payoff $\xi\in \mathcal{L}_m^2(\mathcal{F}_T)$ and with $Z= \phi \sigma$, we then get a backward process:  
    \[
        Y_t = \xi - \int_t^T \left[ \frac{\mu_s}{\sigma_s}Z_s +r_s \left(Y_s-\frac{Z_s}{\sigma_s}\right)_+ - R_t \left(Y_s-\frac{Z_s}{\sigma_s}\right)_- \right]ds - \int_t^T Z_s dW_s,
    \]
which lies in the framework of nonlinear BSDEs.

\color{black}
Stochastic control is undoubtedly a major field of application.
To briefly describe stochastic control in our notation, consider the stochastic optimization problem $\sup\nolimits_k J(k)$ in the weak formulation \textcolor{black}{\cite{BSDEFin_Kar_Pen}}:
\begin{equation}\label{stochastic control}
 J(k):=\mathbb{E}_k\left[\Phi(X_T)+\int_0^T f(s,X_s,k_s)ds\right],
\end{equation}
where $\Phi$ and $f$ take values in $\mathbb{R}$, and $k$ controls the SDE as in
\[
 X_t={\bf x}+\int_0^t b(s,X_s,k_s)ds+\int_0^t \sigma (s,X_s)dB^k_s,\quad t\in [0,T],
\]
and $\mathbb{E}_k$ is the expectation with respect to the probability measure $\mathbb{P}_k$ under which $B_t^k:=B_t-\int_0^t \sigma^{-1}(s,X_s)b(s,X_s,k_s)ds$ is a Brownian motion. 
Then, the following linear BSDE under $\mathbb{P}_k$ has a unique solution $(Y^k,Z^k)$:
\[
 Y_t^k=\Phi(X_T)+\int_t^Tf(s,X_s,k_s)ds-\int_t^TZ_s^k dB_s^k,
\]
due to Theorem \ref{Pa_Pe_Unique} with $\xi=\Phi(X_T)$, provided that all the conditions are met.
In view of $J(k)=Y_0^k$, the stochastic control problem \eqref{stochastic control} is well within the framework \eqref{eq:BSDE_Diff}.

\subsection{Forward backward stochastic differential equations}
\label{subsection FBSDEs in preliminaries}


We next review forward-backward stochastic differential equations (FBSDEs), which have already made a brief implicit appearance in the example on European option pricing (Section \ref{subsection BSDEs in preliminaries}). 
Consider the following BSDE:
	\begin{equation} 
    \begin{cases} \label{eq:FBSDE_Diff}
         -dY_t = f(t, X_t, Y_t, Z_t)dt - Z_tdW_t, \quad t \in [0,T], \\
        Y_T = \Phi(X_T),
    \end{cases}
	\end{equation}
	where $(Y,Z)_{t \in [0,T]}$ takes values in $\mathbb{R}^m \times \mathbb{R}^{m \times d}$, $X$ is the $\mathbb{R}^q$-valued diffusion process solving the standard stochastic differential equation:
	\begin{equation} 
    \begin{cases} \label{eq:FBSDE_SDE}
          dX_t = \mu(t,X_t)dt + \sigma(t,X_t)dW_t, \quad t \in [0,T], \\
        X_0 = {\bf x},
    \end{cases}
	\end{equation}
	and $f : [0,T] \times \mathbb{R}^q \times \mathbb{R}^m \times \mathbb{R}^{m \times d} \rightarrow \mathbb{R}^m$, $\Phi:\mathbb{R}^q \rightarrow \mathbb{R}^m$, $\mu : [0,T] \times \mathbb{R}^q \rightarrow \mathbb{R}^q$ and $\sigma : [0,T] \times \mathbb{R}^q \rightarrow \mathbb{R}^{q \times d}$ are all given functions, and ${\bf x}_0$ is a suitable point in $\mathbb{R}^q$.
We note that the dimension of the $Y$ component is set to one ($m=1$) in many instances, particularly when dealing with deterministic PDEs (Theorems \ref{Th:solution and viscosity} and \ref{Th:Equivalence_PDE}).
\textcolor{black}{The pair of equations \eqref{eq:FBSDE_Diff}-\eqref{eq:FBSDE_SDE} is known as the Markovian BSDE, or the (uncoupled) FBSDE.}
	We now introduce the ``Markovian'' counterparts of the assumptions made in Theorem \ref{Pa_Pe_Unique} which ensure the existence and uniqueness of solutions to the FBSDE \eqref{eq:FBSDE_Diff}-\eqref{eq:FBSDE_SDE}, as well as an equivalence between the solution of the FBSDE to a viscosity solution of the parabolic PDE \eqref{eq:PDE_semilinear}.
	We refer the reader to, for instance, \cite{Pardoux2014book, zhang_2017} for more detail and proofs.

\begin{assumption}\label{assumption zhang_5.0.1}{\rm
	Assume the following conditions on the FBSDE \eqref{eq:FBSDE_Diff}-\eqref{eq:FBSDE_SDE}:
    \begin{itemize}
    \setlength{\parskip}{0cm}
    \setlength{\itemsep}{0cm}
    \item $\mu$, $\sigma$, $f$ and $\Phi$ are uniformly Lipschitz continuous in $({\bf x},{\bf y},{\bf z})$.
    \item $\mu(\cdot,0)$, $\sigma(\cdot,0)$, $f(\cdot,0,0,0)$ and $\Phi(0)$ are bounded.
    \item $\mu$, $\sigma$ and $f$ are uniformly H\"older-$(1/2)$ continuous in $t$.\qed
    \end{itemize}
    }\end{assumption}

	    \begin{theorem}[\bf Existence, uniqueness and PDE equivalence]\label{Th:solution and viscosity}
	Under Assumption \ref{assumption zhang_5.0.1}, there exists a unique solution $(X,Y,Z)$ to the FBSDE \eqref{eq:FBSDE_Diff}-\eqref{eq:FBSDE_SDE}. 
	Define $v:[0,T] \times \mathbb{R}^q \to \mathbb{R}$ by 
	\begin{equation}
	v(t,{\bf x}):=Y_t^{t,{\bf x}}=\mathbb{E}\left[ \Phi(X_T^{t,{\bf x}}) + \int_t^T f(s,X_s^{t,{\bf x}},Y_s^{t,{\bf x}},Z_s^{t,{\bf x}}) ds \right], \quad (t,{\bf x}) \in [0,T]\times \mathbb{R}^q,
	\end{equation}
	 where $(X^{t,{\bf x}},  Y^{t,{\bf x}}, Z^{t,{\bf x}})$ denotes the adapted solution to the FBSDE \eqref{eq:FBSDE_Diff}-\eqref{eq:FBSDE_SDE}, restricted to $[t,T]$ with $X_t^{t,{\bf x}}= {\bf x}$, a.s.
	 Then, $v$ is a viscosity solution to the parabolic PDE:
	\begin{equation} 
    \begin{cases} \label{eq:PDE_semilinear}
         ((\partial/\partial t)+ \mathcal{L}_t) v(t,{ \bf x}) + f(t, {\bf x}, v(t,{\bf x}), (\nabla v(t,{\bf x}))^{\top}\sigma(t,{\bf x})) = 0, & (t,{\bf x})\in [0,T)\times \mathbb{R}^q,\\
        v(T,{\bf x}) = \Phi({\bf x}),& {\bf x}\in \mathbb{R}^q,
    \end{cases}
	\end{equation}
	where $\mathcal{L}$ denotes the second-order differential operator $\mathcal{L}_tv(t,{\bf x}) := \langle \mu(t,{ \bf x}), \nabla v(t,{\bf x})\rangle+\frac{1}{2}{\rm tr}[\sigma^{\otimes 2}(t,{\bf x}){\rm Hess}(v(t,{\bf x}))]$.
	\end{theorem}


\textcolor{black}{Analytic solutions are rarely available for FBSDEs, especially those of practical relevance.
However, for illustrative purposes, we present a few examples of the FBSDE in the form \eqref{eq:FBSDE_Diff}-\eqref{eq:FBSDE_SDE}, for which unique solutions are available in closed form, to indicate how simple the FBSDE must be for an analytic solution to be available.}
First, set $\mu (t, x) = \theta(a- x)$ and $\sigma (t, x) = b$.
With a suitable initial state $X_0=x_0\in\mathbb{R}$, the forward component \eqref{eq:FBSDE_SDE} is then a diffusion process in the form of $X_t=e^{-\theta t}x_0+a(1-e^{-\theta t})+b \int_0^t e^{-\theta (t-s)}dW_s$, that is, an Ornstein-Uhlenbeck process if $\theta>0$.
With the quadratic terminal $\Phi(x) = x^2$ fixed, if the driver is given by $f(t, x, y, z) = ry+e^{-r(T-t)}b^2+\theta (a-x)z/b$,
for some $r, a, b \in \mathbb{R}$, then the FBSDE \eqref{eq:FBSDE_Diff} admits a unique solution given by $(Y_t,Z_t)=(e^{-r(T-t)} X_t^2,2 b e^{-r(T-t)} X_t)$ for $t\in [0,T]$.
Or, if the driver is set, instead, to $f(t, x, y, z) = 2 \theta (a- x)x+b^2$, for some $a, b, \theta \in \mathbb{R}$, then the unique solution to the FBSDE \eqref{eq:FBSDE_Diff} is available in closed form as $(Y_t,Z_t)=(X_t^2,2b X_t)$ for $t\in [0,T]$.
Next, if the terminal, the driver and the forward component are given by $\Phi(x) = \ln(x)$, $f(t, x, y, z) = a -
b^2/2$, $\mu (t, x) = a x$, $\sigma(t, x) = b x$, and $X_0=x_0>0$ for some $a, b \in \mathbb{R}$, then the FBSDE \eqref{eq:FBSDE_Diff}-\eqref{eq:FBSDE_SDE} admits the unique solution $(Y_t,Z_t) = (\ln(X_t),b)$ for $t\in [0,T]$ and the geometric Brownian motion $X_t=x_0\exp[(a-b^2/2)t+bW_t]$ for $t\in [0,T]$.
\textcolor{black}{We once again stress that these simple examples are only presented for demonstration purposes.
FBSDEs of practical interest are generally not explicitly solvable and thus do require numerical treatment.}

We next state the nonlinear Feynman-Kac formula with the representation of the derivative of the PDE solution. 
We need the following additional assumption to give the representation under the Lipschitz continuous condition of $f$ and $\Phi$.  
For details (including the definition of the solution $\nabla X$ of a variational equation) and applications we refer the reader to, for instance, \cite{Mall_rep, Zhang_thesis}. 

\begin{assumption}\label{Assumption_3}{\rm
	Assume the following conditions:
    \begin{itemize}
    \setlength{\parskip}{0cm}
    \setlength{\itemsep}{0cm}
    \item $\sigma$ is uniformly elliptic on $[0,T] \times \mathbb{R}^q$, that is, there exists $C>1$ such that $C^{-1}\|\xi\|^2\le \langle \xi, \sigma^{\otimes 2}(t,{\bf x})\xi \rangle \le C \|\xi\|^2$ for all $t \in [0,T]$ and ${\bf x},\xi \in \mathbb{R}^q$.
    \item $\mu$ and $\sigma$ are in $\mathcal{C}^{1}_b$ in ${\bf x}$.\qed 
    \end{itemize}
    }\end{assumption}
    
    \begin{theorem}[\bf Nonlinear Feynman-Kac formula with the representation of the derivative of the PDE solution]\label{Th:Equivalence_PDE} 
	Under Assumptions \ref{assumption zhang_5.0.1} and \ref{Assumption_3}, 
	the viscosity solution $v$ to the parabolic PDE \eqref{eq:PDE_semilinear} is in $\mathcal{C}^{0,1}_b([0,T)\times \mathbb{R}^q;\mathbb{R})$ and it holds that
	\begin{equation}\label{BSDE_PDE}
	    (Y_t, Z_t) = (v(t,X_t), (\nabla v(t,X_t))^{\top}\sigma(t,{X_t})),\quad t\in [0,T),
	\end{equation}
	and 
	\begin{equation}\label{representation_BSDE}
	\nabla v(t,{\bf x})=\mathbb{E}\left[ \Phi(X_T^{t,{\bf x}}) N_T^{t,{\bf x}} + \int_t^T f(s,X_s^{t,{\bf x}},Y_s^{t,{\bf x}},Z_s^{t,{\bf x}}) N_s^{t,{\bf x}} ds \right],\quad (t,{\bf x})\in [0,T)\times \mathbb{R}^q,
	\end{equation}
	where $N^{t,{\bf x}}_s:=(s-t)^{-1}( \int_t^s (\sigma^{-1}(r,X_r^{t,{\bf x}}) \nabla X_r^{t,{\bf x}} )^\top dW_r)^\top$ for $s\in (t,T]$.
	\end{theorem}


From here on out, with some exceptions (including Sections \ref{subsection branching diffusion methods} and \ref{subsection Lipschitz BSDEs}), all the conditions of Theorems \ref{Pa_Pe_Unique} and \ref{Th:solution and viscosity} are imposed as standing assumptions.
Further assumptions (such as Assumption \ref{Assumption_3}) may be imposed in different forms for different numerical methods depending on the BSDEs and FBSDEs of interest.





    

From a practical perspective, FBSDEs have long been actively studied in the pricing and hedging problem in the incomplete market \cite{Hedge_constraints3}, which remains an ongoing topic, for instance, \cite{B_n_zet_2021}. 
With the aid of the relevant techniques developed there, it has been known that FBSDEs play an essential role for formulating X-Value Adjustment (XVA), which is a collective term for various valuation adjustments for derivative instruments, for instance, credit risk and funding costs.
Ever since the catastrophic wave of bankruptcies of big financial firms in the 2008 global financial crisis, the XVA pricing has been realized as one of most urgent problems in derivatives pricing.
To be more precise, XVAs can be formulated in the form of a FBSDE with a Lipschitz continuous driver $f$, which is non-differentiable at some points, such as $f(t, x, y, z) =   a(t) \max \{0, y\} + b(t).$
Ever since the pioneering work \cite{piterbarg2010funding}, a variety of extensions and refinements of the XVA formulation have been investigated, certainly in the framework of FBSDEs, from both theoretical and practical standpoints \cite{bichuch2016arbitragefree, https://doi.org/10.1111/mafi.12146, doi:10.1142/S0219024913500064, crepey_song, Lesniewski2016}.
In particular, various numerical methods for the XVA pricing have also been proposed, such as a branching algorithm \cite{riskhenry2012}, a Fourier-based discretization method \cite{doi:10.1137/16M1099005}, higher-order discretization methods \cite{doi:10.1080/1350486X.2019.1637268}, a dual algorithm for the stochastic control problem \cite{doi:10.1137/15M1019945} and the deep learning algorithm \cite{gnoatto2020deep}.

Despite being outside the scope of Assumption \ref{assumption zhang_5.0.1}, the American option pricing problem \cite{BSDE_American, BSDE_American2} can also be formulated in the framework of FBSDEs with a discontinuous driver of the form $f(t, x, y, z) =a(t) \mathbbm{1}(y< \Phi(x))+ b(t)$, based on which various numerical methods have been proposed. Examples of such methods include the perturbative expansion and particle method \cite{Fujii:2015aa} and a local polynomial approximation and branching processes method \cite{bouchard2019monte}. 

As such, there exist a variety of advanced features and relevant terminologies around FBSDEs and their applications in the literature, which we summarize here in brief.
First, the FBSDE in the form \eqref{eq:FBSDE_Diff}-\eqref{eq:FBSDE_SDE} is called a decoupled FBSDE, in contrast to ``coupled'' when the coefficients $\mu$ and $\sigma$ in the forward component \eqref{eq:FBSDE_SDE} depend on the backward component \eqref{eq:FBSDE_Diff} (Section \ref{subsection coupled FBSDEs}).
Moreover, the term ``fully coupled'' is often attached when those coefficients depend on the outcome $\omega$.
The backward component may be reflected at a given stochastic process (Section \ref{subsection reflected BSDEs numerical methods}), while FBSDEs may contain jumps in the backward and/or forward component (Section \ref{subsection BSDEs with jumps}).
If the drift of the backward component \eqref{eq:FBSDE_Diff} contains the second order derivative of the corresponding PDE \eqref{eq:PDE_semilinear}, then such a coupled FBSDE is referred to as a second-order FBSDE (Section \ref{section 2bsde}).
If not only the coefficients $\mu$ and $\sigma$ depend on the law of the processes $(X,Y,Z)$ but also the driver $f$, then such a coupled FBSDE is called a McKean-Vlasov FBSDE (Section \ref{subsection mckean-vlasov}).


\section{Backward numerical methods}
\label{section backward numerical methods}

With the basic background reviewed in Section \ref{section preliminaries}, we begin with various numerical methods for BSDEs which work backwards in time under two categories with respect to the degree of discretization, namely, backward Euler methods (Section \ref{subsection backward euler methods}) and higher-order methods (Section \ref{subsection higher-order methods}).
In each subsection, we present one or two representative numerical methods in some detail, while we skip lengthy technical explanations.
In particular, when describing the algorithms in each subsection, we do not go into much technical detail in order to avoid digression from the main idea. We then give a brief overview of various other methods in the category.

\textcolor{black}{In general, a discretization method is said to be of order $p$ if the discretization error has the convergence rate $\mathcal{O}(n^{-p})$ in a suitable norm where the given time interval is discretized into $n$ subintervals.
However, this convergence rate does not necessarily represent the efficiency of the numerical algorithm under consideration.
That is, in most instances, backward numerical methods require the additional computation of conditional expectations in each subinterval, 
whose efficiency depends largely on many factors, such as the dimension and the method applied.
As we are only focusing on backward methods alone here, we will discuss various methods for computing conditional expectations collectively in Section \ref{section computation of conditional expectations}, such as least-squares regression (Section \ref{subsection LS regression}), Malliavin calculus based methods (Section \ref{subsection malliavin calculus based methods}), quantization methods (Section \ref{subsection quantization methods}), tree based methods (Section \ref{subsection tree based methods}) and cubature methods (Section \ref{subsection cubature methods}).}

\subsection{Backward Euler methods}
\label{subsection backward euler methods}


The first class of backward numerical methods we review are the so-called backward Euler schemes. 
\textcolor{black}{The main references on the topic are 
\cite{backward_discrete_4, backward_discrete_3}, although the ideas of the backward Euler scheme date back to \cite{backward_discrete_1}. The convergence of the backward Euler scheme is of order $1/2$ under standard Lipschitz assumptions and with Lipschitz terminal condition \cite{backward_discrete_4, backward_discrete_3}. 
However, the scheme has a convergence rate of $1$ if the forward SDE can be simulated perfectly on the grid or if it is approximated by a higher order scheme under the conditions that the coefficients are sufficiently smooth \cite{backward_discrete_6}.}
We note that discretization methods of higher order are often collectively called higher-order discretization methods, but we do not discuss them here and instead devote Section \ref{subsection higher-order methods} to them.

\color{black}
In the present context, there are certainly two general categories of explicit and implicit discretization schemes for the FBSDE \eqref{eq:FBSDE_Diff}-\eqref{eq:FBSDE_SDE} (see, also, \cite[Section 5.3.2]{zhang_2017}), which can be summarized as follows:
\begin{algorithm}[H]
{\small
\SetAlgoLined
 \textbf{Initialization}: Approximate the terminal condition $Y^n_{t_n}=\Phi(X^n_{t_n})$ with the Euler-Maruyama scheme $X^n$. \\
 \For{$k = (n-1)$ {\rm to} $0$}{
  \begin{align}\label{backward Z} 
    Z^n_{t_k}&=\frac{1}{t_{k+1}-t_k}\mathbb{E}\left[Y^n_{t_{k+1}}(W_{t_{k+1}}-W_{t_k})^{\top}\Big|\,\mathcal{F}_{t_k}\right],\\
    Y^n_{t_k}&=\begin{dcases}
 \mathbb{E}\left[Y^n_{t_{k+1}}+f(t_k,X^n_{t_k},Y^n_{t_{k+1}},Z^n_{t_k})(t_{k+1}-t_k)\Big|\,\mathcal{F}_{t_k}\right],& \text{(explicit)}\\
 \mathbb{E}\left[Y^n_{t_{k+1}}\Big|\,\mathcal{F}_{t_k}\right]+f(t_k,X^n_{t_k},Y^n_{t_k},Z^n_{t_k})(t_{k+1}-t_k).& \text{(implicit)}
\end{dcases}\label{backward Y}
  \end{align}
 }
 \caption{}}
\end{algorithm}

As usual, the implicit scheme often provides better properties and performance relative to the explicit scheme, with these benefits coming in exchange for the additional computing effort for solving the defining equation for $Y^n_{t_k}$, which appears on both sides of the implicit scheme.

\color{black}
\textcolor{black}{
With $(Y^{n},Z^{n})_{t\in [0,T]}$ understood to be the pair of the step processes whose state at $t_k$ corresponds to \eqref{backward Z}-\eqref{backward Y},  
it was first derived 
\cite{backward_discrete_4, backward_discrete_3}
that \begin{equation}
\sup_{t\in [0,T]} \mathbb{E}\left[  \left\|Y_{t}-Y^{n}_{t}\right\|^2\right] + \mathbb{E}\left[ \int_0^T \left\|Z_{t}-Z^{n}_{t}\right\|^2 \right]  \leq \frac{C}{n}.\label{backward_Euler}
\end{equation}
The backward Euler scheme in the reflected case converges at order $1/2$, when $f$ is independent of $Z$ \cite{Quantiz2}, but only at order $1/4$ in the general case  \cite{bouchard2008discrete, ma2005representations}.
We note that convergence may depend on the regularity of the terminal condition in general through the $L^2$-regularity of $Z$ \cite{GEISS20122078}.
Stability analysis is conducted in \cite{10.2307/24512614} for the backward Euler method in a similar manner to the Euler method for ordinary differential equations.}

Backward Euler methods have been successfully generalized to address a broader class of BSDEs with nonstandard features (Section \ref{section BSDEs with nonstandard features}), such as BSDEs with a general terminal condition and driver \cite{10.1214/11-AAP762}, BSDEs driven by an infinite activity Poisson random measure in addition to the Brownian motion \cite{MADAN20161553}, reflected BSDEs \cite{bouchard2008discrete, ma2005representations, sun2020quantitative}, BSDEs with jumps \cite{Mall_Calc1, massing2021approximation} and quadratic BSDEs \cite{richou2011}.
Backward Euler methods have also been tailored to, for instance, mean-field BSDEs \cite{Zhang:2022vq}, $G$-FBSDEs (FBSDEs driven by $G$-Brownian motion) \cite{HU2021578}, FBSDEs driven by c\`adl\`ag martingales \cite{KHEDHER2016508}, and backward stochastic Volterra integral equations of type I \cite{bender2013discretization,1605.04865} and of type II \cite{hamaguchi2021approximations}.
The convergence of the backward Euler method for BSDEs is discussed in \cite{10.1214/11-AAP762} with general terminal condition and driver.

We close this section with 
\cite{BouchardElieTouzi+2009+91+124, 10.1214/10-AAP723, 10.1214/14-AAP1030}, 
where 
some studies focus to construct numerical methods for nonlinear PDEs based on backward Euler methods for FBSDEs.
One may also look to \cite{gong_zhao_2017, zhang_zheng_2002} for related studies on discretization, as well as \cite{Finite_diff1} where a number of finite difference methods are reviewed for solving BSDEs, similar to backward Euler methods for the most part.
 

\subsection{Higher-order discretization methods}
\label{subsection higher-order methods}


In the literature, a variety of discretization methods have been proposed to achieve a higher order of convergence 
than the backward Euler method of order $1$ (Section \ref{subsection backward euler methods}). 
\textcolor{black}{In principle, to achieve a higher order of convergence, the following three steps need to be planned and implemented carefully:
\begin{enumerate}
\setlength{\parskip}{0cm}
\setlength{\itemsep}{0cm}
\item We first discretize the BSDE.
A typical result on the discretization of a solution to a BSDE $(Y,Z)$ is given by  
\begin{equation}
 \sup_{k\in \{0,1,\cdots,n\}}\mathbb{E}\left[ \left\|Y_{t_k}-Y^{\pi}_{t_k}\right\|^2+\frac{1}{n}\left\|Z_{t_k}-Z^{\pi}_{t_k}\right\|^2 \right]^{1/2} \leq \frac{C}{n^2},\label{highorder}
\end{equation}
where discretization of the forward component $X$ has not been taken into account. 
Here, $Y^{\pi}$ and $Z^{\pi}$ denote suitable simple schemes involving conditional expectation representations for $Y$ and $Z$, respectively (see, for instance, \cite{Crisan2014}).
\item The corresponding forward component $X$ is to be discretized using a suitable high-order discretization scheme 
and in such a way that the overall error remains within the order of the backward component (like, the inequality \eqref{highorder}).
\item By then applying a computation method for the conditional expectation (Section \ref{section computation of conditional expectations}), an efficient numerical scheme 
is obtained which has a faster convergence rate than a scheme with the backward Euler method (Section \ref{subsection backward euler methods}). 
\end{enumerate}
A variety of relevant numerical methods have been developed in conjunction with existing techniques in each of those three steps.
We devote the present section to a review of such higher-order discretization methods.}

For illustrative purposes, we begin with a simple numerical method of higher-order type based upon a random walk approximation of the Brownian motion. 
This method is given in the context of FBSDEs whose driver is independent of $X$, although the method can be easily extended to when one does not have this independence.
The forward (SDE) component is approximated using the Euler-Maruyama method or Milstein method.
The solution to the FBSDE is approximated by a discretization scheme that is derived by using the so-called theta method.
For example, the numerical approximation for $Y$ is derived by first showing that
    \begin{equation*}
         Y_{t_k} = \mathbb{E}\left[Y_{t_{k+1}}\big|\,\mathcal{F}_{t_k}\right] + \int_{t_k}^{t_{k+1}} \mathbb{E}\left[f(s,Y_s,Z_s)\big|\,\mathcal{F}_{t_k}\right]ds,
    \end{equation*}
    and then approximating the integral using a theta approximation: a convex combination of an explicit term (term at time $t_{k+1}$) and an implicit term (term at time $t_k$):
    \begin{align}\label{theta_approximation}
        Y_{t_k} \approx \mathbb{E}\left[Y_{t_{k+1}}\big|\,\mathcal{F}_{t_k}\right] + \Delta_n\theta_1 f(t_k,Y_{t_k},Z_{t_k})
        + \Delta_n(1 - \theta_1) \mathbb{E}\left[f(t_{k+1},Y_{t_{k+1}},Z_{t_{k+1}})\big|\,\mathcal{F}_{t_k}\right], \quad \theta_1 \in [0,1],
    \end{align}
where $\Delta_n := T/n (= t_{k+1} - t_k)$ is a fixed time step. 
In a similar manner, a numerical approximation for $Z$ is derived, and then a backward numerical scheme is presented using these approximations.
Namely, the scheme for approximating $(Y,Z)$ by $(Y^{\pi}, Z^{\pi})$ is given as $Y_{t_n}^{\pi} = \Phi(X_{t_n}^{\pi})$, $Z_{t_n}^{\pi}=(\nabla\Phi(X_{t_n}^{\pi}))\sigma(t_n,{X_{t_n}^{\pi}})$, and then
    \begin{align}\label{theta_scheme_discrete}
        Z_{t_k}^{\pi} &= -\frac{1- \theta_2}{\theta_2}\mathbb{E}\left[Z_{t_{k+1}}^{\pi}\Big|\,\mathcal{F}_{t_k}\right] + \frac{1}{\theta_2\Delta_n}\mathbb{E}\left[Y_{t_{k+1}}^{\pi}\Delta W_k\Big|\,\mathcal{F}_{t_k}\right] + \frac{1- \theta_2}{\theta_2} \mathbb{E}\left[f(t_{k+1},Y_{t_{k+1}}^{\pi},Z_{t_{k+1}}^{\pi}) \Delta W_k\Big|\,\mathcal{F}_{t_k}\right], \nonumber \\
        Y_{t_k}^{\pi} &= \mathbb{E}\left[Y_{t_{k+1}}^{\pi}\Big|\,\mathcal{F}_{t_k}\right] + \theta_1\Delta_n f(t_k,Y_{t_k}^{\pi},Z_{t_k}^{\pi}) +  (1-\theta_1)\Delta_n \mathbb{E}\left[f(t_{k+1},Y_{t_{k+1}}^{\pi},Z_{t_{k+1}}^{\pi})\Big|\,\mathcal{F}_{t_k}\right],
    \end{align}
    for $k\in \{n-1,\cdots, 0\}$ backwards, with $\theta_2 \in (0,1]$.
    The FBSDE and PDE equivalence (Theorem \ref{Th:Equivalence_PDE}) is applied in computing the terminal condition for $Z$.
    Hence, Assumption \ref{Assumption_3} is essential here.
    The above scheme is proved in \cite{Error_estimates_theta} to converge with order of convergence $2$ when $\theta_1 = \theta_2 = 1/2$ and order of convergence $1$ otherwise.
    
Since its first development \cite{doi:10.1137/05063341X}, the $\theta$-scheme has been extensively generalized in various directions.
For instance, as previously mentioned, for BSDEs whose driver is independent of the control process $Z$, the $\theta$-scheme is proved \cite{Error_estimates_theta} to converge with order $2$ when $\theta_1 = \theta_2 = 1/2$ and order $1$ otherwise. 
In order to achieve the same orders ($2$ when $\theta_1 = \theta_2 = 1/2$ and order $1$ otherwise) even when the driver depends on the pair $(Y,Z)$, the $\theta$-scheme is generalized \cite{1531-3492_2012_5_1585} by introducing more parameters.  
The concept of the $\theta$-scheme has been applied to more general problem settings, for instance, mean-field BSDEs \cite{doi:10.1137/17M1161944}, BSDEs whose driver contains not only the present value $(Y_t, Z_t)$ but also the future values 
\cite{hu2019explicit}, and $G$-BSDEs \cite{hu2019numerical}.
In particular, the $\theta$-scheme with $\theta=1/2$, termed the Crank-Nicolson scheme \cite{wang_luo_zhao_2010}, is investigated along with error analysis \cite{LI20101612, IJNAM-10-876} and applications to FBSDEs \cite{li_yang_zhao_2017}, where higher-order discretization formulas are derived for the control process $Z$ via approximation of the derivative of the backward process $Y$ by the well-known Crank-Nicolson method, as well as to mean-field BSDEs \cite{Zhang:2022vq}.
Inspired by the $\theta$-scheme, a numerical method for one-dimensional FBSDEs is proposed \cite{SghirHadiri+2020+79+91} without using conditional expectations, along with a few analytical solutions of particular FBSDEs.
From an implementation point of view, the $\theta$-scheme has been further accelerated via parallel computing \cite{5662503}.

To achieve even higher orders of convergence, the multistep scheme \cite{doi:10.1137/09076979X} is developed as an extension of the $\theta$-scheme with convergence at an arbitrary order when the driver is independent of the control process.
In the multistep scheme, the pair $(Y_{t_k}^{\pi}, Z_{t_k}^{\pi})$ is approximated based on the expectation of values at multiple future time points $\{(Y_{t_l}^{\pi}, Z_{t_l}^{\pi})\}_{l\in\{k+1,\cdots, k+L\}}$, for which evidently more extensive computing is required than the backward Euler method (Section \ref{subsection backward euler methods}) and the $\theta$-scheme.
Error analyses are conducted in the cases where the driver depends on the control process \cite{doi:10.1137/120902951} and with variable step sizes \cite{hanqiang}, backed up by convincing numerical results.
Moreover, by exploiting the independence of the computation at each grid point, the multistep scheme is tailored \cite{kapllani2019multistep} for parallelized computing on GPUs.

The concept of multistep methods has also been quite prevalent in the context of FBSDEs.
A multistep scheme is developed for coupled FBSDEs \cite{doi:10.1137/130941274} along with a crude Euler-Maruyama method employed for discretizing the forward component, with convergence analysis refined in \cite{yang_zhao_2015}.
The multistep method has been generalized to address more general BSDEs with nonstandard features (Section \ref{section BSDEs with nonstandard features}), such as BSDEs with jumps \cite{EAJAM-6-253, fu_zhao_zhou_2016}, second-order FBSDEs \cite{zhao_zhou_kong_2017} and fully nonlinear parabolic PDEs \cite{kong_zhao_zhou_2015}, as well as improved in combination with relevant techniques, such as the sparse-grid method \cite{JCM-31-221}, the spectral sparse grid approximation with fast Fourier transform \cite{1531-3492_2017_9_3439}, the Lagrange interpolation polynomial \cite{liu2020fully,1937-1632_2021044}, finite difference for approximating the derivatives of the solution of FBSDEs in multi-time levels \cite{teng_zhao_2021}, 3-point Gauss-Hermite quadrature rule and non-equidistant difference scheme  \cite{pak_kim_rim_2020}, polynomial approximation \cite{zhao_zhang_ju_2016}, and predictor‐corrector and least‐squares Monte Carlo schemes \cite{Han:2022uh}.
Naturally, the multistep method for BSDEs and FBSDEs 
parallels with the relevant numerical methods for ODEs and PDEs. 
For instance, the Runge-Kutta method for ODEs is applied to BSDEs \cite{chassagneux2014AAP}, where the order barrier is exhibited to be more restrictive for BSDEs than for ODEs.
The defferred correction method for ODEs is tailored to coupled FBSDEs \cite{tang_zhao_zhou_2017} and second-order FBSDEs \cite{yang_zhao_zhou_2019}.
It is proved \cite{10.1093/imanum/drab023} that a linear multistep method for FBSDEs is stable as long as the so-called root condition is met in a similar manner to that for ODEs.
A general framework is constructed \cite{doi:10.1137/19M1260177} to investigate the stability, consistency and convergence of discretization schemes for FBSDEs in a unified manner, including the backward Euler method, the $\theta$-scheme and various types of the multistep method.
In \cite{kapllani_teng_2022, kapllani2019multistep}, the authors demonstrate an acceleration of the multistep method for FBSDEs, relying on parallel GPU computing using CUDA. 
%
It is natural to employ established higher-order discretization schemes for the forward component, such as the Milstein and stochastic Taylor methods, together with the multistep method for BSDEs \cite{fu_zhao_2014, zhao_yang_fu_2014, zhao_zhang_ju_2014}.
Along this line, various first- and second-order discretization methods have been developed by combining the multistep method for the backward component (such as the trapezoidal rule) with a suitable discretization scheme for the forward component, and then compute the resulting conditional expectations efficiently by the Gauss-Hermite quadrature rule \cite{zhao_zhang_ju_2014}, a Lagrangian interpolation \cite{fu_zhao_2014}, and an approach with Malliavin derivative \cite{zhao_yang_fu_2014}. 
Error estimates are improved in \cite{math9080848} for the schemes developed in \cite{zhao_zhang_ju_2014, zhao_zhang_ju_2016}.
Those methods are also applied to mean-field FBSDEs \cite{sun_zhao_2020}.

In \cite{Crisan2014}, a second order discretization scheme is constructed for BSDEs with a nonsmooth boundary data using Brownian weights for the approximation of the $Z$ component under weaker conditions, where the gradient estimate of \cite{Crisan_and_Delarue} is employed. 
\textcolor{black}{In order to construct numerical methods for these FBSDEs,  
the backward component is discretized by the trapezoidal or Simpson's rule, and the forward component by the cubature method on Wiener space (Section \ref{subsection cubature methods}). 
Such methods result in the orders of convergence 1 \cite{doi:10.1137/090765766}, 2 \cite{Crisan2014}, and up to 3 \cite{doi:10.1080/1350486X.2019.1637268}}.
A similar approach is employed in \cite{de2015cubature} for discretizing McKean-Vlasov FBSDEs. 
When the forward and backward components are discretized, respectively, by the cubature and backward Euler methods, one can construct a discretization method on the basis of an explicit error expansion \cite{chassagneux_trillos_2020}, resulting in an arbitrary order of convergence with the aid of the Richardson-Romberg extrapolation.
Inspired by the development of \cite{Crisan2014}, a second-order discretization method is proposed in \cite{naito_yamada_2019}, 
by introducing higher-order correction terms for the backward component in the form of a sum of polynomials of the Brownian motion through the Malliavin weights. 

As previously mentioned, the computation of conditional expectations is undoubtedly a major step in higher-order discretization methods (Section \ref{section computation of conditional expectations}).
Here, Fourier analysis may play central roles for various purposes \cite{10.1093/imanum/drx022, Ge:2022tp, HUIJSKENS2016593, Hyndman:2017aa, jrfm15090388, doi:10.1137/130913183, doi:10.1137/21M1444679}, for instance, via the characteristic function of the discretized forward component \cite{RUIJTER20161}.
A tree-based regression is suggested in \cite{teng2019review} for efficiently implementing higher-order discretization methods. 
In conjunction with the four-step scheme \cite{Four_step1}, higher-order discretization methods are investigated, for instance, on an approximation of spatial derivatives in the coefficients of the four-step scheme using finite difference \cite{10.1093/imanum/drl019}, and with the Hermite-spectral method \cite{doi:10.1137/06067393X}.
Finally, in \cite{Parallel_FD}, the authors explored the potential for parallel computing with FBSDEs when using the binomial tree type approximation.
Due to the special structure of the numerical method proposed (which is similar to that in \cite{BSDE_Discrete_Theta_Binomial_Tree}), a block allocation algorithm is developed in parallelization, where large communication overhead is avoided, backed up by numerical illustration of encouraging speedups for the parallel implementation.

\section{Computation of conditional expectations}
\label{section computation of conditional expectations}

In the present section, we survey major techniques for the computation of conditional expectations in the context of numerical methods (mostly, backward numerical methods (Section \ref{section backward numerical methods})) for BSDEs, namely, least-squares regression based methods (Section \ref{subsection LS regression}), Malliavin calculus based methods (Section \ref{subsection malliavin calculus based methods}),
quantization methods (Section \ref{subsection quantization methods}), tree based methods (Section \ref{subsection tree based methods}) and cubature methods (Section \ref{subsection cubature methods}).
\textcolor{black}{We remark that some deep learning based methods may also be interpreted as a mechanism for computing the conditional expectation in numerical methods for BSDEs, using nonlinear least-squares Monte Carlo with a deep neural network architecture. We do not review those methods here but review deep learning based methods collectively in Section \ref{section deep learning}, as this topic forms a rapidly growing field of research deserving of a separate section.}


\subsection{Least-squares regression based methods}
\label{subsection LS regression}

Methods that fit into this category are ones which use a form of least-squares regression to evaluate the conditional expectations appearing in a discretization of the BSDE. 

Here, we proceed with representative methods of \cite{backward_discrete_5, Lin_regress1} to illustrate this type of method in a clear and concise manner. 
Let $\Delta_k=T/n$, that is, $t_k=kT/n$ and $\Delta W_k=W_{t_{k+1}}-W_{t_{k}}$ for $k\in \{0,1,\cdots,n\}$.
\color{black}
For each time point $t_k$, we take $\mathbb{R}^K$-valued deterministic function bases $(e^K_{i,k})_{i\in\{0,1,\cdots,d\}}$, whose elements are given by, for instance, the sequence of Hermite polynomials or that of Laguerre polynomials of size $K$. 
 For every $k\in \{1,\cdots,n\}$, let $\{\Delta W_k^m\}_{m\in \{1,\cdots,M\}}$ be independent copies of $\Delta W_k$, and let $\{X_{k}^{\pi,m}\}_{m\in \{1,\cdots,M\}}$ be corresponding copies of $X_{k}^{\pi}$.  
An algorithm for approximating the backward Euler scheme $(Y_k^\pi,Z_k^\pi)$ 
can be described as follows. 

\begin{algorithm}[H]
{\small
\SetAlgoLined
 Set $y_n^{n,M,K}(\cdot) =\Phi$ \\
 \For{$k = (n-1)$ {\rm to} $1$}{
  \begin{equation} 
            \alpha_{i,k}^{M,K}=\argmin_\alpha \frac{1}{M} \sum_{m=1}^M \left|y_{k+1}^{n,M,K}(X_{k+1}^{\pi,m})  \frac{\Delta W_k^{i,m}}{\Delta_k} - \alpha \cdot e^{K}_{i,k}(X_{k}^{\pi,m}) \right|^2, \quad i\in \{1,\cdots,d\}. \label{lsm_z}
  \end{equation}
  Put $z_{k}^{n,M,K}=(z_{1,k}^{n,M,K},\cdots,z_{d,k}^{n,M,K})$ where $z_{i,k}^{n,M,K}(\cdot)= \alpha_{i,k}^{M,K} \cdot e^K_{i,k}(\cdot)$.
  \begin{equation} 
          \alpha_{0,k}^{M,K}=\argmin_\alpha \frac{1}{M} \sum_{m=1}^M \left|y_{k+1}^{n,M,K}(X_{k+1}^{\pi,m}) + \Delta_k f(t_k,X_{k}^{\pi,m},y_{k+1}^{n,M,K}(X_{k+1}^{\pi,m}),z_{k}^{n,M,K}(X_{k}^{\pi,m})) -  \alpha \cdot e^K_{0,k}(X_{k}^{\pi,m}) \right|^2. \label{lsm_y}
  \end{equation}
  Put $y_{k}^{n,M,K}(\cdot)= \alpha_k^{M,K} \cdot e^K_{0,k}(\cdot)$.
 }
 Return $Y_0^{\pi,M,K}=\frac{1}{M} \sum_{m=1}^M (y_{1}^{n,M,K}(X_{1}^{\pi,m}) + \Delta_1 f(t_0,{\bf x},y_{1}^{n,M,K}(X_{1}^{\pi,m}),z_{1}^{n,M,K}(X_{1}^{\pi,m})))$.
 \caption{}}
\end{algorithm}

Under the Lipschitz continuity in the state variables, the $(1/2)$-H\"older continuity in time of the coefficients of the Markovian BSDE, and the condition that for all measurable functions $\varphi$ such that $\varphi(X_k^\pi) \in L^2(\Omega)$, there is $(\beta_k^K)_k$ such that $\beta^K \cdot e_{i,k}^K(X_k^\pi) \to \varphi(X_k^\pi)$ in $L^2(\Omega,\mathscr{F}^{X_k^\pi})$ as $K \to \infty$.
It is shown in \cite{zhang_2017} that $\lim_{K\to +\infty} y_k^{n,K}(X_k^\pi) = Y_k^{\pi}$ in $L^2(\Omega)$ for all $k\in \{0,1,\cdots,n-1\}$, where $y_k^{n,K}$ and $z_k^{n,K}$ denote the theoretical means of the respective empirical ones $y_k^{n,M,K}$ and $z_k^{n,M,K}$ of \eqref{lsm_z} and \eqref{lsm_y} above, as well as
\[ 
    \lim_{M\to +\infty} Y_0^{\pi,M,K} =Y_0^{\pi,K},
\] 
almost surely, by the strong law of large numbers.
While the convergence holds by taking sufficiently large $K$ and $M$, its rate has been found to be quite complex. 
It is shown in \cite{Lin_regress1}, which extends \cite{backward_discrete_5}, that the squared error of the time discretization is of the order $\mathcal{O}(n^{-1})$, and roughly speaking, the numbers of basis functions and of the paths are to be chosen as $K \approx n^{d}$ and $M \approx n^{d+3}$ in order to control the global squared error to be $\mathcal{O}(n^{-1})$, upon a suitable choice of basis functions. 
We refer the reader to, for instance, \cite{10.1007/978-3-642-25746-9_8,backward_discrete_5, Lin_regress1,zhang_2017} for detailed analyses.

\color{black}

In a similar spirit, a numerical scheme is developed in \cite{LSMDP_noMal} for solving the multi-step forward dynamic programming equation arising from the discretization of the FBSDE.
The resulting sequence of conditional expectations is computed using empirical least-squares regressions.
Also, in \cite{backward_discrete_7}, we see another algorithm based on least-squares Monte Carlo (LSMC).
Here, the algorithm is designed to allow large scale parallelization of the computations on highly multicore processors, such as GPUs by stratifying sample paths to minimize the exposure to the memory requirements due to the storage of simulations.
Possible discontinuities at the interfaces due to piecewise polynomial bases \cite{backward_discrete_7} are avoided \cite{gobet_lopezsalas_vasquez_2020, proceedings2019021044} by instead employing smooth orthonormal basis functions.
Also, \cite{10.1007/978-3-642-25746-9_8, backward_discrete_5} propose numerical schemes based on iterative least-square regressions on function bases, where the involved coefficients are evaluated using Monte Carlo simulations. 
\textcolor{black}{We mention the so-called regression-later approach \cite{FBSDE_Discrete_Regress_Later}, which can be thought of as a variant of least-squares regression methods.}

The primal-dual methodology is generalized in \cite{German_primal_dual} to a backward dynamic programming equation associated with time discretization schemes of reflected BSDEs (Section \ref{subsection reflected BSDEs numerical methods}).
They suggest a pathwise approach to the dynamic programming equation, which avoids the evaluation of conditional expectations in backward recursion in time.
This approach leads to a minimization problem and is thus thought of as a dual minimization problem.
Under suitable assumptions, $Y$ can also be represented as the supremum over a class of classical optimal stopping problems, and this maximization problem can be seen as a primal problem.
Using the representations for $Y$ as the value of a maximization and a minimization problem, confidence intervals are constructed for $Y_0$ using LSMC (even if only to generate an initial input approximation), along with a few numerical examples to test the performance of the algorithm for multi-dimensional reflected BSDEs in nonlinear pricing problems.

In \cite{doi:10.1080/00207160.2019.1658868}, the authors apply the Stochastic Grid Bundling Method (SGBM) to develop a numerical scheme.
The SGBM algorithm involves the approximation of conditional expectations by means of bundling Monte Carlo sample paths and a local regress-later technique within each bundle. 
By employing Hermite martingales, the problem of solving a FBSDE is formulated in \cite{pelsser2019} as the problem of solving a countably infinite-dimensional system of ODEs.
On this basis, they develop a numerical scheme which involves the projection of the solution onto generalized Hermite polynomials. 
Similarly, a numerical scheme is proposed in \cite{TENG2020117} based on the projection of conditional expectations onto cubic spline polynomials.
Finally, we refer to \cite{naito_yamada_2019} for a second-order discretization method for FBSDEs based on an algorithm which utilizes polynomials of Brownian motions and is implemented by use of a least-squares Monte Carlo method.



\subsection{Malliavin calculus based methods}
\label{subsection malliavin calculus based methods}

Malliavin calculus based methods are closely related to least-squares regression based methods (Section \ref{subsection LS regression}) in the sense that a lot use least-squares regression to solve the conditional expectations resulting from a discretization of a BSDE which incorporates Malliavin weights.

\subsubsection{Representation of conditional expectations using integration-by-parts}\label{subsubsection malliavin ibp}

Malliavin calculus was first applied to BSDEs in \cite{backward_discrete_4} for the purpose of computing conditional expectations in line with \cite{BET2004, FLLL2001}, for instance, the two terms $\mathbb{E}[Y^n_{t_{k+1}}(W_{t_{k+1}}-W_{t_k})^{\top}|\,\mathcal{F}_{t_k}]$ and $\mathbb{E}[Y^n_{t_{k+1}}|\,\mathcal{F}_{t_k}]$ appearing in the implicit discretization scheme BSDEs \eqref{backward Z} and \eqref{backward Y}.
With the aid of the Markovian property, those conditional expectations boil down to the form $\mathbb{E}[\phi(X_s)|\,X_t={\bf x}]$ for $s\in (t,T]$, where $\phi:\mathbb{R}^q \to \mathbb{R}$ and $\{X_t:\,t\in [0,T]\}$ is the forward process \eqref{eq:FBSDE_SDE}.
Here, we describe in brief how this conditional expectation can be reformulated without conditioning in such a way that Monte Carlo methods can play an effective role.
In general, for two smooth Wiener functionals $F$ and $G$ in the sense of Malliavin, if $G$ is ``non-degenerate", then the conditional expectation $\mathbb{E}[F|G={\bf x}]$ can be written as fractions:
\begin{equation}\label{skorohod_integral_rep}
\mathbb{E}[F|G={\bf x}]=\frac{\mathbb{E}[F \delta_{\bf x} (G)]}{\mathbb{E}[\delta_{\bf x}(G)]}=\frac{\mathbb{E}[H_{\bf x}(G) \delta^W( F u ) ]}{\mathbb{E}[H_{\bf x}(G) \delta^W(u)]},
\end{equation}
where $\delta_{\bf x}(G)$ is understood to be the composition of the delta function at ${\bf x}$ and $G$ as a Watanabe distribution on the Wiener space, $H_{\bf x}$ denotes a Heaviside function, that is, $H_{\bf x}({\bf c}):=\prod_{k=1}^q\mathbbm{1}(c_k\ge x_k)$, and $\delta^W$ denotes the Skorohod integral operator.
The second equality is due to Malliavin integration by parts and the process $u$ satisfying $\int_0^T D_s G u(s)ds=1$.
Under the further constraint $\int_0^T D_s F u(s)ds=0$, by applying the property of the Skorohod integral operator $\delta^W(Fu)=F\delta^W(u)-\int_0^T D_s F u(s)ds$ and
setting $F=\phi(X_s)$, $G=X_t$ and $u(s)=(D_s X_t)^{-1} ( t^{-1}\mathbbm{1}(s\in (0,t])-(T-t)^{-1} \mathbbm{1}(s\in (t,T]))$, the last expression in \eqref{skorohod_integral_rep} can be rewritten in the following form, which lends itself to estimation of the conditional expectation $\mathbb{E}[\phi(X_s)|\,X_t={\bf x}]$ by Monte Carlo methods: 
\begin{equation}\label{malliavin fraction}
\mathbb{E}[\phi(X_s)|X_t={\bf x}]
=\frac{\mathbb{E}[\phi(X_s)H_{\bf x}(X_t) \delta^W(u)]}{\mathbb{E}[H_{\bf x}(X_t) \delta^W(u)]}\approx \frac{\frac{1}{M} \sum_{k=1}^M \phi(X_s^k)H_{\bf x}(X_t^k) \delta^{W^k}(u)}{\frac{1}{M}\sum_{k=1}^M H_{\bf x}(X_t^k) \delta^{W^k}(u)},
\end{equation}
by generating independent copies with sufficiently large sample size $M$.
This approach can also be applied to reflected BSDEs (Section \ref{subsection reflected BSDEs numerical methods}) in \cite{backward_discrete_4}.

As is clear, the effectiveness of this approach depends, next to the Heaviside function, on the Skorohod integral term $\delta^W(u)$.
In principle, one needs a reasonably explicit expression for this term, without which simulation would be essentially impossible.
Even if one has an explicit form, its computation can be prohibitive when the problem dimension is high.
To address this issue, a variant of the method above is proposed in \cite{Mall_improve} for a significant reduction of the numerical complexity by wisely modifying the Skorohod integral term $\delta^W(u)$ in such a way to reduce the terms involved and skip differentiation of the drift and diffusion coefficients of the forward component.


\color{black}
\subsubsection{Malliavin weights dynamic programming with regression}\label{subsection dynamic programming}

We next describe a dynamic programming scheme with the regression methods, named the Malliavin weights dynamics programming, an application of Malliavin calculus in a different spirit from Section \ref{subsubsection malliavin ibp}.  
For illustrative purposes, we proceed with the 
regression method \cite{Mall_LSregres}, among many candidates, that solves a dynamic programming equation with Malliavin weights arising from the time-discretization of FBSDEs (which may also be fit for Section \ref{section backward numerical methods}). 

\color{black}
Now, the main idea of the dynamic programming scheme \cite{Mall_LSregres} is to approximate the solution $(Y,Z)$ of a FBSDE by the discrete time stochastic processes $(Y^{\pi},Z^{\pi})$, defined on the partition $\pi_n$, which are given as follows:
		\begin{equation}\label{MWDP}
        \begin{dcases} 
         Y_k^{\pi} = \mathbb{E}\left[\Phi(X_T) + \sum_{j=k}^{n-1}f_j(Y^{\pi}_{j+1},Z^{\pi}_j)\Delta_j\Big|\,\mathcal{F}_{t_k} \right],&  k \in \{0, \cdots, n\}, \\
          Z_k^{\pi} = \mathbb{E}\left[\Phi(X_T) H_n^{(k)} + \sum_{j=k+1}^{n-1}f_j(Y^{\pi}_{j+1},Z^{\pi}_j)H_j^{(k)}\Delta_j \Big|\,\mathcal{F}_{t_k}\right],& k \in \{0, \cdots, n-1\},
        \end{dcases}
	    \end{equation}
	    where $(\omega,y,{\bf z})\to f_j(\omega,y,{\bf z})$ is $\mathcal{F}_{t_j}\otimes \mathcal{B}(\mathbb{R})\otimes \mathcal{B}(\mathbb{R}^{1\times d})$-measurable.
This system is then solved backwards in the order $Y_n^{\pi}, Z_{n-1}^{\pi},Y_{n-1}^{\pi}$ and so forth. 
	    It is called the Malliavin weights dynamic programming (MWDP) equation, as it literally takes the form of a dynamic programming equation with Malliavin weights. 
	    The MWDP equation \eqref{MWDP} is inspired by \cite{Mall_rep}, which gives the following representation of the control process $Z$:
	    \begin{equation}\label{Z_formula}
	        Z_t = \mathbb{E}\left[\Phi(X_T) H_T^{(t)} + \int_t^T f(s,X_s,Y_s,Z_s)H_s^{(t)}ds \Big|\,\mathcal{F}_t\right],
	    \end{equation}
		where the processes $(H_s^{(t)})_{0 \leq t < s \leq T}$ are the Malliavin weights defined as follows:
		\begin{equation}\label{Mall_weights}
		    H_s^{(t)} = \frac{1}{s-t}\left(\int_t^s \left(\sigma^{-1}(r,X_r)D_t X_r \right)^{\top} dW_r \right)^{\top}, \quad 0 \leq t < s \leq T,
		\end{equation}
		with $(D_t X_r)_t$ denoting the Malliavin derivative of the marginal $X_r$.
	  We note that this method adapts the least-squares multi-step forward dynamic programming algorithm to a scenario which incorporates Malliavin weights.
\textcolor{black}{Under suitable technical conditions presented in \cite{LSMDP_noMal},} it holds almost surely that $|Y_k^{\pi}| \le C(1 +(T - t_k)^{\theta_c})=:C_{y,k}$ and
		\[
		  \|Z_k^{\pi}\| \leq C\left(\frac{\esssup_{\omega}\mathbb{E}[|\Phi(X_T) - \mathbb{E} [\Phi(X_T)|\,\mathcal{F}_{t_k}]|^2|\,\mathcal{F}_{t_k}](\omega)}{(T-t_k)^{1/2}} + \frac{1}{(T-t_k)^{1/2 - \theta_c}} + (T - t_k)^{\theta_L/2}\right)=:C_{z,k},
		\]
		for all $k\in \{0,\cdots,n-1\}$.
	    Furthermore, there exists a constant $C_{y,z}$ such that
		$|Y_k^{\pi}| + \sqrt{T-t_k}|Z_k^{\pi}| \leq C_{y,z}$ almost surely for all $k\in \{0,\cdots,n\}$.
 
 The conditional expectations (and hence processes $Y^{\pi}$ and $Z^{\pi}$) appearing in the MWDP equation \eqref{MWDP} are computed using a Monte Carlo least-squares regression scheme (Section \ref{subsection LS regression}). 
 Due to the Markovian assumptions, there exist measurable and deterministic (but unknown) functions $y_k(\cdot):\mathbb{R}^q \rightarrow \mathbb{R}$ and $z_k(\cdot):\mathbb{R}^q \rightarrow \mathbb{R}^{1 \times d}$ for all $k \in \{0,\cdots,n-1\}$ such that the solution $(Y^{\pi}_k,Z^{\pi}_k)_{k \in \{0, \cdots, n-1\}}$ of the MWDP equation \eqref{MWDP} is given by $(Y^{\pi}_k,Z^{\pi}_k) := (y_k(X_k^{\pi}),z_k(X_k^{\pi}))$.
The aim of the Monte Carlo regression scheme is thus to estimate these functions. Ordinary least-squares regression (OLS) is defined in such a way that easily allows path-dependence and joint laws. The authors reformulate the MWDP equation \eqref{MWDP} in terms of the given definition of OLS.
Specifically, take $\mathcal{K}_k^{(l)}$ to be any dense subset in the $\mathbb{R}^l$-valued functions belonging to $L_2(\mathcal{B}(\mathbb{R}^q),\mathbb{P} \circ (X_k^{\pi})^{-1})$, and then for each $k \in \{0,\cdots,n-1\}$,
\begin{equation}\label{OLS_bad}
\begin{dcases}
		     y_k(\cdot) \text{ solves OLS}\left(\Phi(x_n) + \sum_{j=k}^{n-1} f_j(y_{j+1}(x_{j+1}),z_j(x_j))\Delta_j,\mathcal{K}_k^{(1)},\mathbb{P}\circ (H_{k+1}^{(k)}, \cdots, H_n^{(k)},X_k^{\pi},\cdots, X_n^{\pi})^{-1}\right),\\ 
		     z_k(\cdot) \text{ solves OLS}\left(\Phi(x_n)h_n + \sum_{j=k+1}^{n-1} f_j(y_{j+1}(x_{j+1}),z_j(x_j))h_j\Delta_j,\mathcal{K}_k^{(d)},\mathbb{P}\circ (H_{k+1}^{(k)}, \cdots, H_n^{(k)},X_k^{\pi},\cdots, X_n^{\pi})^{-1}\right),
\end{dcases}
\end{equation}
for $(h_{k+1},\cdots,h_n) \in (\mathbb{R}^{1 \times d})^{n-k}$ and $(x_k,\cdots, x_n) \in (\mathbb{R}^q)^{n-k+1}$. 
In words, for instance, $y_k(\cdot)$ is set to be the least-squares approximation of $\Phi(x_n) + \sum_{j=k}^{n-1} f_j(y_{j+1}(x_{j+1}),z_j(x_j))\Delta_j$ in the space $\mathcal{K}_k^{(1)}$ with respect to the law $\mathbb{P}\circ (H_{k+1}^{(k)}, \cdots, H_n^{(k)},X_k^{\pi},\cdots, X_n^{\pi})^{-1}$.
However, the above least-squares regressions give rise to two computational problems, that is, 
$L_2(\mathcal{B}(\mathbb{R}^q),\mathbb{P} \circ (X_k^{\pi})^{-1})$ is often infinite dimensional, and the integrals of the OLS in \eqref{OLS_bad} are presumably computed using the untraceable law of $(H_{k+1}^{(k)}, \cdots, H_n^{(k)},X_k^{\pi},\cdots, X_n^{\pi})$.
These issues are addressed by approximating $y_k(\cdot)$ and $z_k(\cdot)$ on finite-dimensional function spaces $\mathcal{K}_{Y,k}$ and $\mathcal{K}_{Z,k}$, respectively, with respect to the empirical version $\nu_{k,M}$ of the law $\mathbb{P}\circ (H_{k+1}^{(k)}, \cdots, H_n^{(k)},X_k^{\pi},\cdots, X_n^{\pi})^{-1}$ based on $M$ iid realizations. 

The global error of the algorithm is a weighted time-average of three different errors. The first is the approximation error which relates to the error involved in approximating the functions $(y_k,z_k)$ in the finite dimensional approximation spaces. This accuracy is achieved asymptotically as the number of simulations goes to infinity. The second error term is the usual statistical error, which improves with a larger number of simulations or a smaller dimension of the vector spaces. The third and final error term together is the interdependence error which is related to the inter dependencies between regressions at different times. This error is of the same magnitude as the statistical error terms, up to logarithmic factors. Thus, roughly speaking, the global error is of order of the best approximation errors plus the statistical errors.

A similar method can be found in \cite{Indep_Z_LSMC}, where the solution of the FBSDE is approximated by using a backward MWDP equation and LSMC regression (Section \ref{subsection LS regression}), along with importance sampling to minimize the conditional variance occurring in the LSMC algorithm and accelerate the convergence of Monte Carlo approximation in a similar manner to \cite{importance_samp1}.
    The Radon-Nikodym derivative is not given, as in \cite{importance_samp1}, but rather computed adaptively within the LSMC procedure.
    However, it makes sense to apply importance sampling only if the driver is independent of $Z$.
    If the driver did depend on $Z$, there would be a propagation of ``lack of variance reduction'' on the $Y$ component due to the $Z$ component through the driver.
    This means that it would not be possible to keep track of the benefit of importance sampling for $Y$.
    If the Monte Carlo estimation of $Z$ is made with appropriate variance reduction (suited to $Z$ specifically), then this problem would be avoided, and would allow the driver to depend on $Z$, but this is left to future investigation.


\color{black}

\subsection{Quantization methods}\label{subsection quantization methods}

The quantization method is a numerical scheme for computing the conditional expectation $\mathbb{E}[h(X_{k+1})|X_k={\bf x}]$ for a $\mathbb{R}^q$-valued Markov chain $\{X_k\}_{k\in \{0,1,\cdots,n\}}$ by settling the chain onto a space grid $\Gamma_k:=\{ x_k^1,\ldots,x_k^{N_k} \} (\subset \mathbb{R}^q)$ of a suitable cardinality $N_k\in \mathbb{N}$.
The quantization method can be regarded somewhere in the middle of deterministic and probabilistic methods, because it relies on space grids and weights just like deterministic methods, while the computation of weights often requires Monte Carlo methods.
Such quantization methods are developed and examined in \cite{Quantiz2, Quantiz1, BALLYPAGESPRINTEMS+2001+21+34} for probabilistically solving multidimensional optimal stopping problems, and then applied to develop discretization schemes for reflected BSDEs (Section \ref{subsection reflected BSDEs numerical methods}) based on an optimal discrete spatial quantization tree.
Further improvements on quantized BSDE schemes are provided in \cite{callegaro2021fully, Fiorin_Pages_Sagna, nmeir2021quantizationbased, PAGES2018847}.

Here, we describe quantized BSDE schemes in brief in accordance with \cite{Fiorin_Pages_Sagna}.
On the basis of the explicit backward Euler scheme \eqref{backward Z} and \eqref{backward Y},
we denote by $\widehat{X}^n_{t_k}$ the so-called quantization of $X^n_{t_k}$ taking values in a finite grid $\Gamma_k (\subset \mathbb{R}^d)$, that is, $\widehat{X}^n_{t_k}:=\mathrm{Proj}_{\Gamma_k}(X^n_{t_k})$ where $\mathrm{Proj}_{\Gamma_k}$ is a Borel projection of $\mathbb{R}^d$ onto the grid $\Gamma_k$. 
For example, for a quantization level $N_k \in \mathbb{N}$, an optimal quantization grid $\Gamma_k$ is found according to $\|X^n_{t_k} - \mathrm{Proj}_{\Gamma_k}(X^n_{t_k}) \|_2=\inf \{\|X^n_{t_k} - p(X^n_{t_k}) \|_2: \text{Borel measurable }p: \mathbb{R}^q \to \Gamma (\subset \mathbb{R}^q), \ \mathrm{card}(\Gamma) \leq N_k \}$, based on the nearest neighbor projection $\mathrm{Proj}_{\Gamma_k}$, given by $\mathrm{Proj}_{\Gamma_k}(X^n_{t_k})=\sum_{i=1}^{N_k} x_k^i \mathbbm{1}(X^n_{t_k}\in C_i^k(\Gamma_k))$,
where $N_k:=\mathrm{card}(\Gamma_k)$ and $\{C_i^k(\Gamma_k)\}_{i\in \{1,\cdots,N_k\}}$ is a sequence of Voronoi partitions of $\mathbb{R}^q$, satisfying $C_i^k(\Gamma_k) \subset \{ \xi\in \mathbb{R}^q:\, |x_k^i - \xi |=\min_{j\in\{1,\cdots,N_k\}}|x_k^j - \xi | \}$.
A quantized BSDE can then be defined as follows: 
$\widehat{y}_n(x_i)=\Phi(x_i^n)$ for $i\in \{1,\cdots,N_n\}$, and then for $k\in \{n-1,\cdots,1,0\}$ backwards, 
\[
\widehat{y}_k(x_k^i)=\widehat{\alpha}_k(x_k^i)+f(t_k,x_k^i,\widehat{\alpha}_k(x_k^i),\widehat{\beta}_k(x_k^i))(t_{k+1}-t_k), \quad x_k^i \in \Gamma_k,\quad i\in\{1,\cdots,N_k\},
\]
where 
\[
\widehat{\alpha}_k(x_k^i):=\sum_{j=1}^{N_{k+1}} \widehat{y}_{k+1}(x_{k+1}^j) p_k^{ij}, \quad \widehat{\beta}_k(x_k^i):=\frac{1}{\sqrt{t_{k+1}-t_k}}\sum_{j=1}^{N_{k+1}} \widehat{y}_{k+1}(x_{k+1}^j) \Lambda_k^{ij},
\]
with weights given in the forms of conditional probability and expectations:
\[
 p_k^{ij}:=\mathbb{P}\left(\widehat{X}_{k+1}=x_{k+1}^j|\,\widehat{X}_k=x_k^i\right),\quad \Lambda_k^{ij}:=\mathbb{E}\left[ (W_{t_{k+1}}-W_{t_k}) \mathbbm{1}(\widehat{X}_{k+1}=x_{k+1}^j)|\,\widehat{X}_k=x_k^i\right],
\]
which are to be estimated by Monte Carlo methods if analytic formulas are not available.

In principle, more accurate results are naturally expected by increasing the size of the quantizations, whereas from an implementation point of view, the size should be chosen carefully and reasonably, relative to the problem dimension, especially when the dynamics is multivariate with dependent components, since then the conditional distribution cannot be decomposed into a set of independent univariate distributions.
The complexity of the quantization methods comes largely from the computation of the optimal quantizers and their transition probabilities, if Monte Carlo simulation needs to be employed for the estimation of the transition probabilities.

\subsection{Tree based methods}
\label{subsection tree based methods}

Tree based methods have been found effective in the computation of conditional expectations, particularly in low-dimensional problems, 
with relevant analyses presented in \cite{ECP1030, Indep_Z_Discrete} (while its origin may date back to \cite{backward_discrete_1}). 
In order to illustrate the basics of tree based methods, 
we describe a numerical method proposed in \cite{Indep_Z_Discrete} for a BSDE whose driver is independent of $Z$.
The main idea is to approximate the Brownian motion $W$ by a simple random walk $B_t^n := \frac{1}{\sqrt{n}} \sum_{k=0}^{\lfloor nt \rfloor} \zeta_k^n$ for $t \in [0,1]$, where $\{\zeta^n_k\}_{k \in \{1,\cdots,n\}}$ is a sequence of iid Rademacher random variables.
Then, they use this approximation in a discretized version of the BSDE with a constant time step $\Delta:=1/n$, whose solution is denoted $(Y^n,Z^n)$. The solution to this discretized BSDE is then approximated by $(\widetilde{Y}^n,\widetilde{Z}^n)$, which satisfies the following algorithm.

\begin{algorithm}[H]
{\small
\SetAlgoLined
 \textbf{Initialization}: Approximate the terminal conditions $\widetilde{Y}_{t_n}^n = \xi^n$ and $\widetilde{Z}_{t_n}^n = 0$. \\
 \For{$k = (n-1)$ {\rm to} $0$}{
  \[ 
            \widetilde{X}_{t_k}^n = \mathbb{E}\left[\widetilde{Y}_{t_{k+1}}^n \Big|\, \mathcal{G}_k^n\right], \quad
            \widetilde{Y}_{t_k}^n = \widetilde{X}_{t_k}^n +\frac{1}{n}f(t_k,\widetilde{X}_{t_k}^n), \quad
            \widetilde{Z}_{t_k}^n = \mathbb{E}\left[\left(\widetilde{Y}_{t_{k+1}} + \frac{1}{n}f(t_k,\widetilde{Y}_{t_k}^n) - \widetilde{Y}_{t_k}^n\right) (\Delta B_{t_{k+1}}^n)^{-1} \Big|\, \mathcal{G}_k^n \right].
	   \] 
 }
 \caption{}}
\end{algorithm}
	   
The involved conditional expectations are computed using a tree structure.
	We briefly look at the error involved in approximating $(Y,Z)$ by $(Y^n,Z^n)$ (and $(Y^n,Z^n)$ by $(\widetilde{Y}^n,\widetilde{Z}^n)$).
	It is assumed that $\xi= F(W)$ and thus $\xi^n = F(B^n)$, where $F: \Omega \rightarrow \mathbb{R}$ is a bounded Lipschitz function with respect to the uniform topology on $\Omega$, that is, $|F(\omega_1) - F(\omega_2)| \leq C \sup_{t \in [0,1]}|\omega_1(t) - \omega_2(t)|$ for all $ \omega_1, \omega_2 \in \Omega$.
	   The authors prove that it can then be assumed without loss of generality that the driver is bounded and furthermore that for all $k \in \{0,\cdots,n\}$ and for large $n$, 
	   \begin{equation*}
	      \sup_{\omega \in \Omega}|Y^n_{t_k} - \widetilde{Y}_{t_k}^n| \lesssim \frac{e^{2c} -1}{n}, \quad \sup_{\omega \in \Omega}|Z^n_{t_k} - \widetilde{Z}^n_{t_k}| \lesssim \frac{(e^{2c} -1)(2 + c/n)}{\sqrt{n}},
	  \end{equation*}
	   where $c$ denotes the Lipschitz coefficient of the driver.
	   Also, the sequence $(Y^n,U^n)$ is shown to converge weakly in the Skorohod topology to $(Y, \int Z dW)$, where 
	   \begin{equation*}
	       U_{t_k}^n = \sum_{j=0}^{k-1} Z_{t_j}^n \Delta B_{t_{j+1}}^n = Y^n - F(B^n) - \frac{1}{n} \sum_{j=1}^k f(t_j,Y_{t_j}^n).
	   \end{equation*}
We also mention Donsker's theorem which has been derived in the relevant context \cite{ECP1030, 10.3150/20-BEJ1259}. 
In addition, stability and convergence of the discretized filtration have been analyzed in \cite{antonelli2000filtration, COQUET1998235, SPS_2001__35__306_0}.
With the aid of those results, 
an approximation is investigated \cite{BRIAND2002229} in the sense that
a sequence of solutions to a BSDE driven by a martingale, approximating the Brownian motion, converges to the solution to \eqref{eq:BSDE_Diff}, and is further generalized to BSDEs with random terminal time \cite{toldo2006}.


Tree based methods can also be employed in conjunction with the 
theta approximation \eqref{theta_approximation}, which is a convex combination of explicit and implicit terms:
  \begin{align*}
        Y_{t_k} \approx \mathbb{E}\left[Y_{t_{k+1}}\big|\,\mathcal{F}_{t_k}\right] + \Delta_n\theta_1 f(t_k,Y_{t_k},Z_{t_k})
        + \Delta_n(1 - \theta_1) \mathbb{E}\left[f(t_{k+1},Y_{t_{k+1}},Z_{t_{k+1}})\big|\,\mathcal{F}_{t_k}\right],
    \end{align*}
    for $\theta_1 \in [0,1]$.
 Here, the resulting conditional expectations can be approximated on the basis of a recombining tree structure of a random walk approximation for the Brownian motion with the aid of the Markov property of the approximation of the forward process $X$.
 After all conditional expectations have been approximated in this way, one obtains the following algorithm (in accordance with the $\theta$-scheme \eqref{theta_scheme_discrete}) for computing $Y_{j,k}$ and $Z_{j,k}$ at node $(j,k)$ backward in time in the recombining tree structure (for instance, \cite{BSDE_Discrete_Theta_Binomial_Tree}).

    \begin{algorithm}[H]
    {\small
\SetAlgoLined
 \textbf{Initialization}: Approximate the terminal conditions $Y_{j,n}^{\pi} = \Phi(X_{j,n}^{\pi})$ and $Z_{j,n}^{\pi} = (\nabla \Phi(X_{j,n}^{\pi}))\sigma(t_n,{X_{j,n}^{\pi}})$ for all $j\in \{0,\cdots, n\}$. \\
 \For{$k = (n-1)$ {\rm to} $0$}{
 \For{$j = 0$ {\rm to} $k$}{
  \begin{align*}
        Z_{j,k}^{\pi} &= \frac{\theta_2-1}{2\theta_2}\left[Z_{j+1,k+1}^{\pi} + Z_{j,k+1}^{\pi}\right] + \frac{1}{2 \theta_2\sqrt{\Delta_n}}\left[Y_{j+1,k+1}^{\pi} - Y_{j,k+1}^{\pi}\right]
        + \frac{1 -\theta_2}{2\theta_2}\sqrt{\Delta_n}\left[f(t_{k+1},Y_{j+1,k+1}^{\pi},Z_{j+1,k+1}^{\pi})
         - f(t_{k+1},Y_{j,k+1}^{\pi},Z_{j,k+1}^{\pi})\right],\\
        Y_{j,k}^{\pi} &= \frac{1}{2}\left[Y_{j+1,k+1}^{\pi} + Y_{j,k+1}^{\pi}\right] +  \theta_1\Delta_n f(t_k,Y_{j,k}^{\pi},Z_{j,k}^{\pi})
        + (1-\theta_1)\frac{\Delta_n}{2}\left[f(t_{k+1},Y_{j+1,k+1}^{\pi},Z_{j+1,k+1}^{\pi})
        + f(t_{k+1},Y_{j,k+1}^{\pi},Z_{j,k+1}^{\pi})\right].
    \end{align*}
    }
 }
 \caption{}}
\end{algorithm}

Since the order of convergence for the tree approximation is $1$ and this error dominates the error coming from the $\theta$-discretization, the overall scheme has an order of convergence $1$ for any values of $\theta_1 \in [0,1]$ and $\theta_2 \in (0,1]$. 
The expression for $Y_{j,k}^{\pi}$ above is given implicitly and thus needs approximation by Picard iterations.

We close this section by mentioning
random walk approximations of BSDEs  \cite{cheridito_stadje_2013,Katarzyna,doi:10.1080/07362994.2011.610162,Memin:2008aa,nakayama2002approximation,BSDE_Discrete} and 
an $L_2$-convergence of the random walk approximation of BSDEs derived in \cite{geiss_labart_luoto_2020, 10.3150/19-BEJ1120} 
which makes use of \cite{GEISS20122078}. 
Finally, we finish by noting that parallel computing on GPUs is found effective in accelerating tree-based methods \cite{6128469, Parallel_FD}.

%

\subsection{Cubature methods}
\label{subsection cubature methods}
 
Finally, we review cubature methods on the Wiener space and its application in computing conditional expectations appearing in discretization methods for BSDEs.
Consider the following Stratonovich forward SDE and backward SDE: 
\begin{align*}
dX_s^{t,{\bf x}}&=V_0(X_s^{t,{\bf x}})ds+\sum\nolimits_{i=1}^d V_i(X_s^{t,{\bf x}}) \circ dW_s^i, \quad X_t^{\bf x}={\bf x}, \\
-dY_s^{t,{\bf x}}&=f(X_s^{t,{\bf x}},Y_s^{t,{\bf x}},Z_s^{t,{\bf x}})ds- \sum\nolimits_{i=1}^d Z_s^{i,t,{\bf x}} dW_s^i, \quad Y_T^{t,{\bf x}}=\Phi(X_T^{t,{\bf x}}),  
\end{align*}
where the coefficients satisfy suitable conditions and  $\circ$ denotes the symbol for the Stratonovich integral. 
We say that the positive weights $\{\lambda_k\}_{k\in \{1,\cdots,N\}}$ and the paths of bounded variation $\omega_1,\cdots,\omega_N:\,[0,t] \to \mathbb{R}^d$ define a cubature formula of degree $m$ at time $t$ if the identity
\begin{equation}
\mathbb{E}\left[ \int_{0<t_1<\cdots<t_k<t} \circ dW_{t_1}^{\alpha_1} \cdots \circ dW_{t_k}^{\alpha_k} \right]
=
\sum_{\ell=1}^N \lambda_\ell \int_{0<t_1<\cdots<t_k<t} d\omega_{\ell,t_1}^{\alpha_1} \cdots d\omega_{\ell,t_k}^{\alpha_k}
\label{cub_formula1}
\end{equation}
holds true for any multi-index $\alpha \in \{ 0,1,\cdots, d \}^k$ such that $\#\{i:\,\alpha_i \in \{1,\cdots,d \} \}+2 \#\{i:\,\alpha_i=0 \} \leq m$, with $W^0_t:=t$.
With the cubature measure $\mathbb{Q}_t^m:=\textstyle{\sum_{\ell=1}^N} \lambda_l \delta_{\omega_{\ell}}$ where $\delta_\omega$ denotes the Dirac measure mass at $\omega \in \mathcal{C}([0,t])$, the formula \eqref{cub_formula1} can be rewritten as 
\[
\mathbb{E}\left[ \int_{0<t_1<\cdots<t_k<t} \circ dW_{t_1}^{\alpha_1} \cdots \circ dW_{t_k}^{\alpha_k} \right]
=
\mathbb{E}_{\mathbb{Q}_t^m}\left[ \int_{0<t_1<\cdots<t_k<t} \circ dW_{t_1}^{\alpha_1} \cdots \circ dW_{t_k}^{\alpha_k}\right].
\]
Let $x^{0,{\bf x}}_\cdot(\omega)$ be the solution of the ordinary differential equation obtained by replacing the $d$-dimensional Brownian motion $W$ with a path of bounded variation $\omega$ in $X_\cdot^{0,{\bf x}}$.
Then, we have $\mathbb{E}_{\mathbb{Q}_t^m}[g(X_t^{0,{\bf x}})]=\textstyle{\sum_{\ell=1}^N} \lambda_\ell g(x^{0,{\bf x}}_t(\omega_\ell))$ for a function $g:\mathbb{R}^q \to \mathbb{R}$.
We note that implementable expressions of cubature measures are available for some specific degrees (such as $m=3,5$) and dimensions $d$.
Broadly speaking, the cubature measure of degree $m$ yields the following short-time asymptotics: 
\begin{equation}
 \mathbb{E}\left[g(X_t^{0,{\bf x}})\right]-\mathbb{E}_{\mathbb{Q}_t^m}\left[g(X_t^{0,{\bf x}})\right]=\mathcal{O}(t^{(m+1)/2}), \label{cub_shorttime}
\end{equation}
for a smooth function $g$.
We note that the error term depends on the bounds of the higher order derivatives of $g$. 
For instance, the following method provides an implicit computation scheme for the iterative conditional expectations in the backward Euler scheme: 
\[
\overline{R}_{i,n-1} \Phi({\bf x})=\mathbb{E}_{\mathbb{Q}_{\Delta_{i+1}}^m}\left[\overline{R}_{i+1,n-1}\Phi(X_{t_{i+1}}^{t_i,{\bf x}})\right]+\Delta_{i+1} f\left({\bf x},\overline{R}_{i,n-1} \Phi({\bf x}),\Delta_{i+1}^{-1}\mathbb{E}_{\mathbb{Q}_{\Delta_{i+1}}^m}\left[\overline{R}_{i+1,n-1}\Phi(X_{t_{i+1}}^{t_i,{\bf x}}) (W_{t_{i+1}}-W_{t_i})^{\top}\right] \right),
\]
for $i\in \{n-1,\cdots,0\}$, with a suitable degree $m$ \cite{doi:10.1137/090765766}.
In order to reduce its computational cost, some additional techniques, such as the tree-based branching algorithm, are recommended in implementation. 
A non-uniform partition is found effective to treat possibly nonsmooth boundary data $\Phi$ based on the fact that the corresponding solution of the non-linear PDE $u(t,\cdot)$ is smooth under a suitable condition on $V_i$, even when $\Phi$ is not smooth enough, so that the cubature method \eqref{cub_shorttime} performs well.
Similarly, cubature methods are employed in \cite{Crisan2014} for constructing a second order discretization scheme, in \cite{doi:10.1080/1350486X.2019.1637268} for a third order scheme, and in \cite{de2015cubature} for discretizing McKean-Vlasov FBSDEs (Section \ref{subsection mckean-vlasov}). 
We refer the reader to \cite{chassagneux_trillos_2020} for an error expansion and the complexity control of the cubature method for solving BSDEs.

\color{black}
\section{Forward numerical methods}
\label{section forward numerical methods}

Unlike the backward numerical methods (Section \ref{section backward numerical methods}), the methods that we review in the present section do not inherently work backwards in time, and thus (originally, at least) avoid the computation of conditional expectations (Section \ref{section computation of conditional expectations}).
We hence call this category ``forward'' numerical methods in a collective way.
In a similar manner to Section \ref{section backward numerical methods}, we examine one or two representative methods in some detail in each subsection, followed by a brief overview of various other methods in the category.


\subsection{Picard iteration methods}
\label{subsection picard iteration methods}

In this section, we survey forward numerical methods based on Picard iterations. 
We start with a forward scheme proposed in \cite{Picard_it1} via Picard iterations on sample paths, which is aimed at numerically approximating sample paths of nonlinear FBSDEs, based upon the discretization of a Picard type iteration.
It is assumed that $W$ and $X$ are multidimensional, while $Y$ is univariate. Also, we impose that  	
$\mu$ and $\sigma$ are $1/2$-H\"older continuous with respect to the time variable, 
as well as  
$\Phi$ is Lipschitz. 
	    We note that it is not assumed that the matrix $\sigma$ is quadratic or that $\sigma^{\otimes 2}$ is invertible.

The limit of a Picard type iteration is used to approximate the processes $(Y,Z)$. Specifically, we set $(Y^{(0)},Z^{(0)}) \equiv (0,0)$, and $(Y^{(r)},Z^{(r)})$ as the solution of the following FBSDE: 
     \begin{equation*}
	     Y_t^{(r)} = \Phi(X_T) + \int_t^T f(s,X_s,Y_s^{(r-1)},Z_s^{(r-1)})ds - \int_t^T Z_s^{(r)}dW_s.
	 \end{equation*}
	 By taking conditional expectation, the process $Y^{(r)}$ is given as follows:
	 \begin{equation*}
	     Y_t^{(r)} = \mathbb{E}\left[\Phi(X_T) - \int_t^T f(s,X_s,Y_s^{(r-1)},Z_s^{(r-1)})ds \bigg| \mathcal{F}_t  \right],
	 \end{equation*}
	 and $Z^{(r)}$ is given via the martingale representation theorem, meaning the above Picard iteration is implicit.
	 The authors propose a time discretization of the above iteration which is explicit in time, but still requires the evaluation of conditional expectations.
	 Given a partition $\pi_n$ of $[0,T]$ and an approximation $X^{(\pi)}$ of $X$, set $(Y^{(0,\pi)},Z^{(0,\pi)}) \equiv (0,0)$.
	 Then, for $k\in\{0,1,\cdots,n\}$ forward, define the conditional expectations
	 \begin{align*}
	     Y_{t_k}^{(r,\pi)} &:= \mathbb{E}\left[\Phi(X_T^{(\pi)}) - \sum_{j=k}^{n-1} f(t_j,X^{(\pi)}_{t_j},Y_{t_j}^{(r-1,\pi)},Z_{t_j}^{(r-1,\pi)})\Delta_j \bigg| \mathcal{F}_{t_k}  \right], \\
	     Z_{l,t_k}^{(r,\pi)} &:= \mathbb{E}\left[\frac{ W_{t_{k+1}}^l-W_{t_k}^l}{\Delta_k}\left(\Phi(X_T^{(\pi)}) - \sum_{j=k+1}^{n-1} f(t_j,X^{(\pi)}_{t_j},Y_{t_j}^{(r-1,\pi)},Z_{t_j}^{(r-1,\pi)})\Delta_j\right) \bigg| \mathcal{F}_{t_k}  \right],\quad l\in \{1,\cdots,d\},
	 \end{align*}
	 where $W^l_{\cdot}$ denotes the $l$-th component of $W_{\cdot}$.
The processes $Y^{(r,\pi)}$ and $Z^{(r,\pi)}$ are extended to RCLL processes by using constant interpolation.
\textcolor{black}{By then approximating those conditional expectations, for instance, by using the LSMC regression method (Section \ref{subsection LS regression}),} 
the convergence of the above discretized Picard type iteration is given in the form:
\begin{equation*}
	     \sup_{t \in [0,T]} \mathbb{E}\left[ |Y_t - Y_t^{(r,\pi)}|^2 \right] +\mathbb{E} \left[\int_0^T |Z_s - Z_s^{(r,\pi)}|^2ds \right] \le C\left(|\pi_n| 
      + \left(\frac{1}{2} + C|\pi_n| \right)^r\right),
	 \end{equation*}
	 as $|\pi_n|$ is sufficiently small, provided that $\sup_{t \in [0,T]} \mathbb{E}[|X_t - X_t^{(\pi)}|^2] \leq C|\pi_n|$ and $\sup_{|\pi_n| \leq 1}\mathbb{E}[|\Phi(X_T^{(\pi)})|^2] \leq C$, for some positive constant $C>0$.
	 
	 The discretized Picard type iteration has no (high order) nestings of conditional expectations backwards in time, whereas it does have (lower order) nestings of conditional expectations forward in time in the number of Picard iterations.
	 This turns out to be an advantage from a numerical point of view relative to backward methods.
	 Overall, the error when approximating the involved conditional expectations (by a generic estimator) is significantly lower relative to backward methods. Specifically, it turns out the error grows moderately when the mesh of partition goes to zero and the number of Picard iterations tends to infinity.
	 This is once again an advantage over backward methods, where the error often explodes when the mesh tends to zero.
  The algorithm that puts everything all together is tested on a hedging problem to demonstrate the proven theoretical convergence of the numerical method.
	 We note that a variance reduced version of the algorithm is also mentioned and tested on a numerical example.


This work is extended in \cite{importance_samp1, FBSDE_Importance}, where the authors introduce a variance reduced version of the forward approximation scheme by means of importance sampling, or more specifically, by means of a measure transformation based on a general Radon-Nikodym derivative. 
The technique of importance sampling is generalized from simulating expectations to computing the initial value of a FBSDE.
The convergence of this modified and fully implementable numerical method is proved and the success of the involved generalized importance sampling is illustrated by numerical examples in the context of financial option pricing.

We next look at forward numerical methods which focus on solving the equivalent PDE (Theorem \ref{Th:Equivalence_PDE}) by Picard iterations.
Here, we primarily examine one numerical method, which illustrates this category well.
In \cite{BSDE_AdapCont}, the authors combine two main ingredients to an algorithm which is aimed at solving the PDE equivalent to the general FBSDE (Theorem \ref{Th:Equivalence_PDE}).
Hence, Assumption \ref{Assumption_3} is essential.
To describe the scheme, we prepare some notation. 
     Let $v:[0,T] \times \mathbb{R}^q \rightarrow \mathbb{R}$ be $\mathcal{C}^1$ in space and define $f_v:[0,T] \times \mathbb{R}^q \rightarrow \mathbb{R}$ as the following function
        \begin{equation}\label{fv}
            f_v(t,{\bf x}) := f(t,{\bf x},v(t,{\bf x}),(\nabla v(t,{\bf x}))^{\top}\sigma(t,{\bf x})),
        \end{equation}
        where $f$ is the driver of the FBSDE \eqref{eq:FBSDE_Diff} and $\sigma$ is the diffusion coefficient of the SDE \eqref{eq:FBSDE_SDE}.
    We define the following two random variables
\begin{equation}\label{Psi definition}
        \Psi(s,y,g_1,g_2) := \int_s^T g_1(r,X^{s,\bf{y}}_r)dr + g_2(X^{s,\bf{y}}_T),\quad
        \Psi^N(s,y,g_1,g_2) := \int_s^T g_1(r,X^{N,s,\bf{y}}_r)dr + g_2(X^{N,s,\bf{y}}_T),
\end{equation}
        where $X^{s,{\bf y}}$ (respectively, $X^{N,s,{\bf y}}$) denotes the diffusion process starting from ${\bf y}$ at time $s$ (respectively, its approximation with $N$-time steps) which solves the forward SDE \eqref{eq:FBSDE_SDE}.
        Finally, we define $c_{i,j}(\phi):=\sum_{k,l=0}^{i,j}|\nabla^k (\partial_t)^l \phi|_{\infty}$, with $c_0(\phi) := c_{0,0}(\phi)$ for $\phi\in \mathcal{C}_b^{i,j}$.



We first look at the case where the FBSDE has a driver which is independent of $Y$ and $Z$, that is, $f(t, {\bf x}, y, {\bf z}) = f(t,{\bf x})$, which means the FBSDE is linear.
We demonstrate the adaptive control variate approach (one of the two main ingredients) in this simple setting.
It is then extended to the general FBSDE via Picard iterations (the second of the two main ingredients).
From Theorem \ref{Th:Equivalence_PDE}, we see that in order to solve a linear FBSDE, one can equivalently solve the linear PDE \eqref{eq:PDE_semilinear}, for which an adaptive control variate scheme is developed in \cite{OptControl_Emm_Syl}.
Here, one aims to numerically solve the linear PDE:
\begin{equation}\label{eq:PDE_linear}
		        ((\partial/\partial t)+ \mathcal{L}_t)v + f = 0, \quad v(T,\cdot) = \Phi(\cdot).
		    \end{equation}
Using the Feynman-Kac formula, the probabilistic solution of this PDE can be expressed as the following conditional expectation in accordance with \eqref{Psi definition}:
		    \begin{equation}\label{eq:Feynman}
		        v(t,{\bf x}) = \mathbb{E}_{t,{\bf x}}\left[\Phi(X_T) + \int_t^T f(s,X_s)ds\right] = \mathbb{E}\left[\Psi(t,{\bf x} , f, \Phi )\right],
		    \end{equation}
		    where $X^{t,{\bf x}}$ is the solution to the standard SDE \eqref{eq:FBSDE_SDE} starting from $\bf x$ at time $t$.
		    One wishes to compute a sequence of solutions $(v_r)_{r\in\{0,1,\cdots\}}$ by writing 
		    \begin{equation*}
		        v_{r+1} = v_r + \left(\text{Monte Carlo evaluations of the error} \, (v-v_r)\right),
		    \end{equation*}
		    as this sequence is proved to converge to the solution $v$ of \eqref{eq:PDE_linear}.
This approach is backed up by the probabilistic representation of the error term:
\begin{equation*}
v(t,{\bf x}) -  v_r(t,{\bf x}) = \mathbb{E}\left[\Psi\left( t, {\bf x}, f + ((\partial/\partial t) + \mathcal{L}_t)v_r, \Phi - v_r \right)\right]=:c_r(t,{\bf x}).
\end{equation*}

We now proceed to the second ingredient, which is the most general form of Picard iterations \cite{BSDE_AdapCont} and will be used to approximate the solution of the nonlinear FBSDE by the solutions of a sequence of (simple) linear FBSDEs (ones with driver independent of $Y$ and $Z$), which converge geometrically to $(Y,Z)$.
Define recursively the Picard iterative sequence $(Y^r,Z^r)_{r\in\{0,1,\cdots\}}$ with $(Y^0,Z^0)= (0,0)$, as follows:
	    \[ 
	        -dY_t^{r+1} = f(t,X_t,Y_t^r,Z_t^r)dt - Z_t^{r+1}dW_t, \quad Y_t^{r+1} = \Phi(X_T),\quad t\in [0,T].
	    \] 
	    for $r\in \{0,1,\cdots\}$ forward.
	     This sequence of linear FBSDEs converges to the unique solution $(Y,Z)$ of the original non-linear FBSDE $dt \otimes d\mathbb{P}\mbox{-}a.e.$
		    By writing $Y_t^r = v_r(t,X_t)$ and $Z_t^r = (\nabla v_r(t,X_t))^{\top}\sigma(t,X_t)$ and by Theorem \eqref{Th:Equivalence_PDE}, each linear FBSDE can be equivalently written as a PDE:
		    \begin{equation*}
		        ((\partial/\partial t)+ \mathcal{L}_t)v_{r+1} + f(\cdot,\cdot,v_r, (\nabla v_r)^{\top}\sigma) = 0, \quad v_{r+1}(T,\cdot) = \Phi(\cdot), \quad r \in \{ 0,1,\cdots\}.
		    \end{equation*}
		   The sequence of solutions of linear PDEs $(v_r,\nabla v_r)_{r\in \{0,1,\cdots\}}$ converges in a suitable $L_2$ norm to the solution $(v,\nabla v)$ of the semilinear PDE:
		    \begin{equation*}
		        ((\partial/\partial t)+ \mathcal{L}_t)v + f(\cdot,\cdot,v, (\nabla v)^{\top}\sigma) = 0, \quad v(T,\cdot) = \Phi(\cdot),
		    \end{equation*}
		    which then yields the solution of the non-linear FBSDE by setting $(Y,Z) = (v, (\nabla v)^{\top} \sigma)$.
The adaptive control variate method also converges geometrically, and so when combined, the two ingredients provide an algorithm with geometric convergence.
To describe the algorithm, we employ the notation $\mathcal{L}_n u(s,{\bf x}):= \frac{1}{2}{\rm tr}[\sigma^{\otimes 2}(\phi_n(s),{\bf x}){\rm Hess}(u(s,{\bf x}))] + \langle \mu(\phi_n(s),{ \bf x}),\nabla u(s,{\bf x})\rangle$, where $\phi_n(s):= \max\{t_k\in \pi_n : t_k \leq s\}$ with respect to the partition $\pi_n$ on the interval $[0,T]$.

\begin{algorithm}[H]
{\small
\SetAlgoLined
 \textbf{Initialization}: Set $v_0 \equiv 0$ and assume that an approximate solution $v_r$ of class $\mathcal{C}_b^{1,2}$ has been built at step $r$.
 Take $n$ points $(t_k^r,{\bf x}_k^r)_{k \in \{1, \cdots, n\}} \in [0,T] \times \mathbb{R}^q$ at each step of the iteration.
 \begin{itemize}
    \setlength{\parskip}{0cm}
    \setlength{\itemsep}{0cm}
     \item Evaluate $c_r(t_k^r,{\bf x}_k^r)$ using $M$ iid realizations by
\begin{equation*}
c_r^M(t_k^r, {\bf x}_k^r) = \frac{1}{M} \sum_{m=1}^M \Psi^N \left( t_k^r,{\bf x}_k^r, f_{v_r} + ((\partial/\partial t) + \mathcal{L}_n)v_r, \Phi - v_r\right)_m.
\end{equation*}
	    \item Build the global solution $c_r^M(\cdot)$ based on the values $(c_r^M(t_k^r,{\bf x}_k^r))_{k \in \{1, \cdots, n\}}$ using a linear approximation operator $\mathcal{P}^r$:
		    \begin{equation*}
		        \mathcal{P}^r c(\cdot) = \sum_{k=1}^n c(t_k^r, {\bf x}_k^r)w_k^r(\cdot),
		    \end{equation*}
		    where $(w^r_k)_{k\in\{1,\cdots,n\}}$ are suitable weight functions.
		    The approximation of $v$ at step $r+1$ is then computed as
		    \begin{equation*} \label{next_it}
		        v_{r+1}(t,{\bf x})= \mathcal{P}^r(v_r + c_r^M)(t, {\bf x}).
		    \end{equation*}
 \end{itemize}
 \caption{}\label{algorithm 412}}
\end{algorithm}

We here address a few key questions regarding the algorithm.
First, it is critical how to choose the grid points $(t_k^r,{\bf x}_k^r)_{k \in \{1, \cdots, n\}}$ at each iteration $r$.
There, at each iteration, we take a new grid of $n$ points $(t_k^r, {\bf x}_k^r)_{k \in \{1, \cdots, n\}}$ that are independent and uniformly distributed on $[0,T] \times [-a,a]^q$ for a suitable positive constant $a$.
As we want to solve the PDE on $[0,T] \times \mathbb{R}^q$, we must choose $a$ large enough.
Next, an appropriate norm needs to be prepared for measuring the error.
      By combining results on BSDEs stated in a norm (leading to the integration with respect to the measure $e^{\beta s}ds$ \cite{BSDEFin_Kar_Pen}) and results on the bounds for solutions of linear PDEs in weighted Sobolev spaces (leading to the integration with respect to the measure $e^{-\lambda \|{\bf x}\|}d{\bf x}$ \cite{mu}), the convergence of the algorithm is derived in the norm:
      \begin{equation*}
            \| V \|^2_{\lambda, \beta} := \mathbb{E} \left[ \int_0^T \int_{\mathbb{R}^q} e^{\beta s} \|V_s({\bf x})\|^2 e^{-\lambda \|{\bf x}\|} d{\bf x} ds \right] < +\infty,
      \end{equation*}
for $\lambda>0$, $\beta>0$ and the set of processes $V: \Omega \times [0,T] \times \mathbb{R}^q \rightarrow \mathbb{R}^d$ that are $\mathcal{P}r \otimes \mathcal{B}(\mathbb{R}^q)$-measurable (where $\mathcal{P}r$ is the $\sigma$-field of predictable subsets of $\Omega \times [0,T]$).
        This norm is then employed to measure the error $(Y - Y^r, Z - Z^r)$, corresponding to the approximation error made at step $r$ of the algorithm.
        The expectations appearing in $ \|Y - Y^r \|^2_{\lambda, \beta}$ and $ \|Z - Z^r \|^2_{\lambda, \beta}$ are computed with respect to the law of $X$, $X^n$ and all the possible random variables used to compute $v_r$.

The sequence of approximation operators $(\mathcal{P}^r)_{r\in\{0,1,\cdots\}}$ is allowed to change at each step of the iteration, under a number of technical conditions: measurability, linearity, regularity, boundedness, able to approximate functions and spatial derivatives well and stability and centering property for random function.
We refer the reader to \cite{BSDE_AdapCont} for explicit definitions of each property and an example of a sequence of operators $\mathcal{P}^r$ which satisfy the desired properties.
The operators used are kernel based estimators and are based on the non-parametric technique local averaging.
This scheme is further accelerated via parallel computing \cite{LabartLelong+2013+11+39} by replacing the kernel operator with an extrapolating operator for approximating functions and their derivatives.

The geometric convergence results for the algorithm can be described as follows \cite{BSDE_AdapCont}.
If Assumption \ref{Assumption_3} holds, $\Phi \in \mathcal{C}_{b}^{2+\alpha}$ and $f$ is bounded and Lipschitz, then there exists a constant $K(T)$ such that
        \begin{equation*}
            \|Y - Y^r\|^2_{\lambda, \beta} + \|Z-Z^r\|_{\lambda,\beta}^2 \leq S_k + K(T) \frac{c_{0,2}^2(v)}{n},
        \end{equation*}
        where $S_r \leq \eta S_{r-1} + \epsilon$ for suitable constants $\eta$ and $\epsilon$.
        Moreover, for $\beta$ and $\mathcal{P}$-parameters large enough so that $\eta < 1$, it holds that
        \begin{equation*}
            \limsup_{r \rightarrow +\infty} \|Y - Y^r\|^2_{\lambda, \beta} + \|Z-Z^r\|_{\lambda,\beta}^2 \leq \frac{\epsilon}{1 - \eta} + K(T) \frac{c_{0,2}^2(v)}{n}.
        \end{equation*}
        This result is derived by splitting the error $ \|Y-Y^r\|^2_{\lambda, \beta} + \|Z-Z^r\|_{\lambda,\beta}^2$ into its difference sources, such as an Euler scheme, Picard iteration, approximation operator $\mathcal{P}$ and Monte Carlo simulations. 

\textcolor{black}{In \cite{Mall_rep_method}, the authors propose an analytical approximation scheme which also centers on the representation of FBSDEs given in \cite{Mall_rep}, where a Malliavin calculus method is applied to the forward SDE in conbination with the Picard iteration scheme for the backward component.
Various numerical examples are presented to demonstrate the convergence of the proposed algorithm.
}
        
Before moving on, we mention that there is a line of research on deep learning based algorithms for solving FBSDEs and parabolic PDEs in high dimensions.
The main idea is to make an analogy between the FBSDE and reinforcement learning, where PDEs play an important role from an implementation point of view. 
As the class of deep learning based algorithms cannot be categorized simply as forward methods and has become a very active field of research on its own, we do not discuss it here but set up an individual section (Section \ref{section deep learning}) exclusively for this class of algorithms.

\subsection{Branching diffusion system based methods}
\label{subsection branching diffusion methods}

We next look at the branching diffusion system based numerical methods for nonlinear PDEs and corresponding BSDEs \cite{riskhenry2012, BSDE_branching_1}. 
Consider the following semilinear PDE of KPP (Kolmogorov-Petrovskii-Piskunov) type: 
\begin{equation} 
((\partial/\partial t) +\Delta) u(t,{\bf x})+\beta \left(\sum\nolimits_{k\in \mathbb{N}_0}p_k (u(t,{\bf x}))^k - u(t,{\bf x})\right)=0,\label{kpp_pde}
\end{equation} 
with $u(T,{\bf x})=\Phi({\bf x})$, which can be written as \eqref{eq:PDE_semilinear} with $\mathcal{L}_t=\Delta$ and $f(t,{\bf x},y,{\bf z})=\beta(\sum_{k\in \mathbb{N}_0} p_k y^k-y)$, where $\Delta$ denotes the Laplacian and $\{p_k\}_{k\in \mathbb{N}_0}$ is a probability mass sequence satisfying $p_k\in [0,1]$ and $\sum_{k\in \mathbb{N}_0} p_k=1$.
It is well known \cite{McKean, Skorokhod, watanabe} that the solution to the PDE \eqref{kpp_pde} admits a probabilistic representation based on the branching diffusion system in which every particle dies in an exponential time of parameter $\beta$ and creates $k$ iid descendants with probability $p_k$.
Then, every descendant dies and reproduces iid descendants independently after independent exponential times in accordance with the same mechanism.
In the other words, the solution to the PDE \eqref{kpp_pde} can be represented as
\begin{equation} 
u(t,{\bf x})=\mathbb{E}\left[\prod\nolimits_{k=1}^{N_T} \Phi(Z_T^k)\Big|\,(N_t,Z^1_t)=(1,{\bf x})\right], \label{branching_representation}
\end{equation} 
where $N_T$ is the number of particles alive at time $T$, $Z_T^k$ denotes the position of the $k$th particle at time $T$ and the condition indicates that the system is initialized at time $t$ with one particle at position ${\bf x}$.
For the semilinear PDE with the Laplacian $\Delta$ in \eqref{kpp_pde} replaced by the Ito generator $\mathcal{L}_t$ and a general function $f(y)$: 
\[
((\partial/\partial t)+\mathcal{L}_t) u(t,{\bf x})+f(u(t,{\bf x}))=0, \quad u(T,{\bf x})=\Phi({\bf x}),
\] 
or the corresponding FBSDE \eqref{eq:FBSDE_Diff} and \eqref{eq:FBSDE_SDE} with driver $f(t,{\bf x},y,{\bf z})=f(y)$, the so-called ``marked" branching diffusion method is proposed in \cite{riskhenry2012} to evaluate $u(0,{\bf x})$ based on the equivalence between \eqref{kpp_pde} and \eqref{branching_representation} after a polynomial approximation of the driver $f(y) \approx \beta (\sum_{k=0}^m (a_k/p_k) p_k y^k-y)$ of finite degree $m$, where $\{p_k\}_{k\in \{0,\cdots,m\}}$ here is such that $p_k\in [0,1]$ and $\sum_{k\in \{0,\cdots,m\}} p_k=1$ and the fraction $(a_k/p_k)$ represents the weight to count at vertices of type $k$.
The corresponding FBSDE can thus be evaluated via a fully forward-looking simulation of particles owing to the representation \eqref{branching_representation}.

The computation required for branching diffusion methods comes almost entirely from the construction of branching particles.
Hence, in terms of the problem dimension, those methods are not as prohibitive as Picard iteration methods (Section \ref{subsection picard iteration methods}), let alone the computation of conditional expectations (Section \ref{section computation of conditional expectations}), for which a few dimensions would be the best in practice.
It is reported in \cite{BSDE_branching_1} that eight-dimensional PDEs are successfully solved without an issue.
We note that the marked branching diffusion method is extended further in \cite{agarwal_claisse_2020} to elliptic semilinear PDEs by introducing absorption of particles and in \cite{10.1214/17-AIHP880} to address the case where the nonlinear driver depends on the $Z$ component.
In addition, the branching diffusion system is rich enough to provide probabilistic representations of semilinear PDEs beyond the standard form \eqref{eq:PDE_semilinear}.
For instance, those up to any higher order derivative in the driver (not only up to the first order $\nabla u$ as of \eqref{eq:PDE_semilinear}) have been addressed in \cite{nguwi_penent_privault_2022, penent2022fully} with numerical illustrations.
Moreover, probabilistic representations and associated numerical methods have been developed in \cite{penent2021elliptic, penent2021parabolic} for parabolic and elliptic semilinear PDEs \eqref{eq:PDE_semilinear}, but with the Ito generator $\mathcal{L}_t$ replaced by a fractional Laplacian that corresponds to subordination of the underlying Brownian motion by a stable subordinator.

Interestingly, the concept of the branching diffusion system can be employed for constructing backward numerical methods \cite{BSDE_branching, Indep_Z_branching}, where the driver depends on $Z$ \cite{BSDE_branching} or is independent of $Z$ \cite{Indep_Z_branching}.
Despite the backward nature lies outside the primary scope of Section \ref{section forward numerical methods}, we illustrate such a method in brief, in accordance with the latter for the sake of simplicity. 
The other main assumptions made are that $\Phi:\mathbb{R}^d \rightarrow \mathbb{R}$ is measurable and bounded and that $f:\mathbb{R}^d \times \mathbb{R} \rightarrow \mathbb{R}$ is measurable (as usual), uniformly Lipschitz continuous in its first argument and satisfies linear growth and Lipschitz continuity in its second argument.
	As a consequence, there exists a constant $M \geq 1$ such that $|\Phi(X_T)| \leq M$ and $\|X\| + |Y| \leq M$ on $[0,T]$ almost surely.
	The driver $f$ of the FBSDE is approximated by $f_{l_0}$ which has a local polynomial structure and is given as $f_{l_0}({\bf x},y_1,y_2):=\sum_{j=1}^{j_0} \sum_{l=0}^{l_0} a_{j,l}({\bf x})y_1^l \phi_j(y_2)$, where $(a_{j,l},\phi_j)_{(j,l)\in \{1,\cdots,j_0\}\times \{0,\cdots,l_0\}}$ is a family of continuous and bounded maps. For all $y_1,y_2 \in \mathbb{R}$, $j \in \{1,\cdots,j_0\}$, $l \in \{0,\cdots,l_0\}$ and a positive constant $C$, these maps satisfy $|a_{j,l}| \leq C$, $|\phi_j(y_1) - \phi_j(y_2)| \leq C|y_1 - y_2|$ and $|\phi_j| \leq 1$.
	If the driver was simply approximated by a polynomial, then typically, the approximating FBSDEs would explode in finite time and thus no convergence could be expected. The solution to the FBSDE $(\overline{Y},\overline{Z})$ with driver $f$ replaced by $f_{l_0}$ can approximate the true solution $(Y,Z)$ well whenever $f_{l_0}$ is a good approximation of $f$:
    \begin{equation*}
        \mathbb{E}\left[\sup_{t \in [0,T]}|Y_t-\overline{Y}_t|^2\right] + \mathbb{E}\left[\int_0^T \|Z_s-\overline{Z}_s\|^2ds\right] \leq C \mathbb{E}\left[\int_0^T |f - f_{l_0}|^2(X_s,Y_s,Y_s)ds\right],
    \end{equation*}
where the positive constant $C$ is independent of $f_{l_0}$.
Hence, for accurate approximation, one requires a driver that can be approximated well by polynomials.
	
	The solution to the new FBSDE $(\overline{Y},\overline{Z})$ is then approximated by means of a Picard type iteration scheme. To that end, define the process $\overline{Y}^m$ in a recursive way.
	That is, let $\overline{Y}_T^r := \Phi(X_T)$, and define, on each interval $[t_k,t_{k+1}]$, $(Y_.^r,Z_.^r)$ as the solution on $[t_k,t_{k+1}]$ of
     \begin{equation}\label{Picard_iteration_poly}
         Y_.^r = \overline{Y}_{t_{k+1}}^r + \int_.^{t_{k+1}} f_{l_0}(X_s,Y_s^r,\overline{Y}_s^{r-1})ds - \int_.^{t_{k+1}} Z_s^r dW_s.
     \end{equation}
     Then, $\overline{Y}^r = Y^r$ on $(t_k,t_{k+1}]$, and $\overline{Y}_{t_k}^r = (-M)\vee Y_{t_k}^m \wedge M$.
     The aim is thus to solve the above Picard iteration backwards on each interval $[t_k,t_{k+1}]$, and then using the now completely solved iteration on $[0,T]$, begin solving the next iteration backwards across the intervals. Importantly, this algorithm requires the truncation of the approximation of $Y$ at some given time steps in order to reduce the approximating driver to a globally Lipschitz driver.
     The error due to the Picard iteration scheme is given in the form $|\overline{Y}_t^r - \overline{Y}_t| \leq C(T-t)^r$ for all $t\in [0,T]$ and $r\in\mathbb{N}$, where $|\overline{Y}^r_t|$ is uniformly bounded in $t$ and $r$.
     
	
	Each step of the Picard iteration is conducted backwards on each interval $[t_k,t_{k+1})$ by using a representation of $\overline{Y}^r$ in terms of branching diffusion systems. A set of particles $(X^{(l)})_{l \in K}$ is constructed, where each particle is the solution to the forward SDE (but each with their own Brownian motion) with a killing time $T_l$. At this time, the particle splits into a random number of new particles which follow the same dynamics, but with their own killing times and Brownian motion. Define the set $ \overline{\mathcal{K}}_t^l$ as the collection of particles in the $l$-th generation that were born before or at time $t$, and the set $ \mathcal{K}_t^l$ as the collection of particles in $ \overline{\mathcal{K}}_t^l$ that are still alive at time $t$.
	Also, given $v^{r-1}$ and $v^r(t_{k+1}, \cdot)$, define
    \[
    V_{t,{\bf x}}^r := \left[\prod_{l \in \mathcal{K}_{t_{k+1}-t}} \frac{v^r(t_{k+1},X_{t_{k+1}-t}^{{\bf x},(l)})}{\overline{F}(t_{k+1}-t-T_{l-})} \right] \left[\prod_{l \in \overline{\mathcal{K}}_{t_{k+1}-t} \setminus \mathcal{K}_{t_{k+1}-t}} \frac{\sum_{j=1}^{j_0}a_{j,\xi_l}(X_{T_l}^{{\bf x},(l)})\phi_j(v^{r-1}(t+T_l,X_{T_l}^{{\bf x},(l)}))}{p_{\xi_l}\rho(\delta_l)} \right],
    \]
    for all $(t,{\bf x}) \in [t_k,t_{k+1}) \times \Theta$, where $\Theta$ is a compact subset of $\mathbb{R}^d$.
    Finally, set $v^r(t,{\bf x}) := \mathbb{E}[V^r_{t,{\bf x}}]$ for $(t,{\bf x}) \in (t_k,t_{k+1}) \times \Theta$ and $r \in \mathbb{N}$, and $v^r(t_k,{\bf x}):= (-M) \vee \mathbb{E}[V^r_{t_k,{\bf x}}] \wedge M$, for ${\bf x} \in \Theta$ and $r \in \mathbb{N}$.
    Then, it can be proven that $\overline{Y}_{\cdot}^r = v^r(\cdot, X)$ on $[0,T]$. By then approximating $\mathbb{E}[V]$ by $\widetilde{\mathbb{E}}[\widetilde{V}]$, a Monte Carlo approximation on a finite functional space, and putting everything together, one obtains a numerical method for approximating $(Y,Z)$.
    When implementing the algorithm in practice, it is best to modify the algorithm to avoid lots of expensive Picard iterations.
    Two modifications are suggested which use only one Picard iteration, but provide an accurate estimate nevertheless.

\subsection{Asymptotic expansion}
\label{subsection asymptotic expansion}

The method of asymptotic expansion provides a tailor-made approximation for BSDEs by expanding a nonlinear BSDE into by a sequence of linear BSDEs in a flexible way, yet with a rigorous error analysis.
By taking advantage of simple computation involved, the method offers a fast computing scheme for BSDEs. 
Asymptotic expansion is categorized here as a forward methods since the resulting linear BSDEs here can be computed by forward Monte Carlo schemes alone, without the need for backward or regression schemes. 
We discuss the method of asymptotic expansion here separately from Picard based forward methods,
as PDEs are certainly useful here (as we describe shortly) but not fully essential.

We here describe the method in accordance with \cite{
Malliavin_expansions}.
On a filtered probability space $(\Omega,{\cal F},\mathbb{F},\mathbb{P})$, consider the following FBSDE with a small perturbation $\varepsilon \in [0,1]$: 
\begin{align}
dX_t&=\mu(t,X_t)dt+\sigma(t,X_t) dW_t, \label{AEmethod_forward} \\
-dY_t^\varepsilon&=\varepsilon f(t,X_t,Y_t^\varepsilon,Z_t^\varepsilon)dt-Z_t^\varepsilon dW_t, \quad Y_T^\varepsilon=\Phi(X_T), \label{AEmethod_backward}
\end{align}
where 
the driver $f$ is sufficiently smooth, with corresponding PDE given by
\[
((\partial/\partial t)+{\cal L}_t)v^\varepsilon(t,{\bf x})+\varepsilon f(t,{\bf x},v^\varepsilon(t,{\bf x}),( \nabla v^\varepsilon (t,{\bf x}))^\top \sigma(t,{\bf x}) )=0, \ \ v^\varepsilon(T,{\bf x})=\Phi({\bf x}),
\]
due to Theorem \ref{Th:Equivalence_PDE}.
Recall that $Y_t^\varepsilon=v^\varepsilon (t,X_t)$ and $Z_t^\varepsilon=( \nabla v^\varepsilon (t,X_t))^\top \sigma(t,X_t)=:\nabla v^\varepsilon \sigma(t,X_t)$ hold. 
Now, we expand the process $(Y^{\varepsilon},Z^{\varepsilon})$ around the solution of null driver (linear) BSDE $(Y^{0},Z^{0})$: 
\[
-dY_t^0=-Z_t^0 dW_t, \quad Y_T^0=\Phi(X_T). 
\]
Note that $(Y^0,Z^0)$ is explicitly solvable as $Y_t^0=\mathbb{E}[\Phi(X_T) | {\cal F}_t ]=v^0(t,X_t)$ and $Z_t^0=( \nabla v^0 (t,X_t))^\top \sigma(t,X_t)=\nabla v^0 \sigma (t,X_t)$, where $v^0$ is the solution of the corresponding linear PDE $((\partial/\partial t)+{\cal L}_t )v^0(t,{\bf x})=0$ with $v^0(T,{\bf x})=\Phi({\bf x})$.
The following expansion is then proposed \cite{doi:10.1142/S0219024912500343, Malliavin_expansions}: for every $m \in \mathbb{N}$,  
\begin{equation}
(Y^\varepsilon,Z^\varepsilon) \approx (Y^0,Z^0)+\sum_{k=1}^{m} \frac{\varepsilon^k}{k!} (Y^{(k)},Z^{(k)}), \label{expansion_bsde}
\end{equation}
where each process $(Y^{(k)},Z^{(k)})$ evidently satisfies $(Y^{(k)},Z^{(k)})=((\partial^k/\partial \varepsilon^k)Y^{\varepsilon},(\partial^k/\partial \varepsilon^k)Z^{\varepsilon})|_{\varepsilon=0}$ and, moreover, can be characterized by the linear BSDE:
\[
-dY_t^{(k)}=k\frac{d^{k-1}}{d \varepsilon^{k-1}} f(t,X_t,Y_t^\varepsilon,Z_t^\varepsilon)|_{\varepsilon=0}dt-Z_t^{(k)} dW_t, \quad Y_T^{(k)}=0.
\]
It is essential that every $(Y^{(k)},Z^{(k)})$ is explicitly solvable since each solves a linear BSDE.
The expansion \eqref{expansion_bsde} is justified in \cite{Malliavin_expansions} as follows: for each $\beta>0$, $\lambda>0$ and $m\in\mathbb{N}$, there exists $C>0$ such that 
\[
\left\| v^\varepsilon-\left(v^0+\sum_{k=1}^m \frac{\varepsilon^k}{k!} v^{(k)} \right) \right\|_{\beta,\lambda}^2 + \left\| \nabla v^\varepsilon \sigma - \left(\nabla v^0 \sigma+\sum_{k=1}^m \frac{\varepsilon^k}{k!} w^{(k)}\right) \right\|^2_{\beta,\lambda} \leq C \varepsilon^{2(m+1)},
\]
for all $\varepsilon \in (0,1]$, where 
\[
 v^{(k)}(t,{\bf x}):=\frac{\partial^k}{\partial \varepsilon^k}v^{\varepsilon}(t,{\bf x})|_{\varepsilon=0},\quad w^{(k)}(t,{\bf x}):=\frac{\partial^k}{\partial \varepsilon^k}\nabla v^\varepsilon \sigma(t,{\bf x})|_{\varepsilon=0},
\] for $k\in \mathbb{N}$, with the norm $\| h \|_{\beta,\lambda}:=\int_0^T \int_{\mathbb{R}^q} e^{\beta s}\| h(s,{\bf x}) \|e^{-\lambda |{\bf x}|}d{\bf x}ds$ for $\beta>0$ and $\lambda>0$.
This error estimate is consistent with the expansion (\ref{expansion_bsde}) in the sense of $Y_t^{(k)}=v^{(k)}(t,X_t)$ and $Z_t^{(k)}=w^{(k)}(t,X_t)$ for $k\in\mathbb{N}$.

The development in \cite{Malliavin_expansions} is generalized in \cite{doi:10.1137/14100021X} 
to nonsmooth drivers. 
The development in \cite{doi:10.1137/14100021X} can be summarized as follows: for 
a driver $f$ which is Lipschitz in $(t,x)$ and is Gateaux-differentiable in $(y,z)$ in the sense that 
for any square-integrable predictable processes $\phi$ and $\psi
$, 
there exist $\kappa \in (0,1]$ and a square-integrable predictable process $Df_\cdot (\phi,\psi) $ such that $\mathbb{E}[\int_0^T |(f(t,X_t,Y_t^0+\varepsilon \phi_t,Z_t^0+\varepsilon \psi_t)-f(t,X_t,Y_t^0,Z_t^0)-\varepsilon Df_t(\phi,\psi))/\varepsilon |^2  dt ]=o(\varepsilon^{2\kappa})$,
it holds that
\[
\sup_{t\in [0,T]}\mathbb{E}\left[\left| Y_t^\varepsilon-\left( Y_t^0+ \varepsilon Y_t^{(1)}+\frac{\varepsilon^2}{2}Y_t^{(2)}\right) \right|^2 \right] + \mathbb{E} \left[ \int_0^T\left| Z_t^\varepsilon- \left( Z_t^0+\varepsilon Z_t^{(1)}+\frac{\varepsilon^2}{2}Z_t^{(2)} \right)\right|^2dt \right] = o(\varepsilon^{4+2\kappa}).  
\]
The effectiveness and practical accuracy of the schemes has been well demonstrated via numerical experiments in \cite{doi:10.1142/S0219024912500343, doi:10.1137/14100021X}. 
In practice, the expansion up to the first or second order often provides a sufficient accuracy even for non-smooth drivers.
An efficient implementation with interacting particles is developed in \cite{Fujii_Takahashi_particle_2015}.
We also refer to \cite{10.1007/978-3-319-33446-2_3} for nonlinear Monte Carlo scheme using asymptotic expansion method with interacting particles.  
 
The method of asymptotic expansion has found its application to a variety of models in mathematical finance. For instance, it is employed in \cite{doi:10.1142/S2010139212500152} for a complete market on BSDEs.
Numerical analyses are conducted for American options in \cite{Fujii:2015aa} based on asymptotic expansion method with interacting particles. 
Asymptotic expansion is applied in \cite{crepey_song} for numerical experiments in a counterparty risk model. 
In \cite{Armenti_Crepey_2017}, a dynamic framework is introduced for analyzing XVAs with asymptotic expansion employed in numerical experiments.
A polynomial expansion method is developed in \cite{fujii2016} by approximating target BSDEs via a recursive system of linear ODEs. 
Using asymptotic expansion, numerical approximations are derived in \cite{Agarwal_et_al_2019} for McKean-type anticipative BSDEs arising in initial margin requirement models. 
We close this section with even different lines of research, such as \cite{sekine} for quadratic BSDEs, \cite{doi:10.1080/17442508.2018.1521808} for BSDEs with jumps, \cite{fujii_takahashi_takahashi_2019, doi:10.1142/S2424786320500127, takahashi2021new} for deep learning schemes inspired by asymptotic expansion, and \cite{Fahrenwaldt_2016, FUJII20191492, Mall_rep_method} for similar but different types of expansions in nonlinear PDE and BSDEs. 

Based on asymptotic expansion, though not quite in the same spirit as above, statistical inference has been investigated for BSDEs, as from a practical point of view, it may well happen that some problem coefficients cannot be fully determined upon implementation.
For instance, in the case where the drift coefficient $\mu(t,{\bf x};\theta)$ of the forward process \eqref{eq:FBSDE_SDE} remains partially unspecified in its parameter $\theta$ and the diffusion coefficient $\epsilon \sigma(t,{\bf x})$ is kept small with $\epsilon \approx 0$, the approximation problem for the solution to FBSDEs is investigated in \cite{kutoyants_2016, KUTOYANTS2014111} with maximum likelihood estimation of the unknown parameter $\theta$ concurrently conducted.
The case where the diffusion coefficient is parameterized instead is investigated in \cite{Gasparyan_Kutoyants_2015}.

\subsection{Multilevel Picard approximation}\label{subsection multilevel picard}


\textcolor{black}{Multilevel Picard approximation is an emerging class of numerical methods, mainly consisting of the following three steps \cite{E_2021_nonlinearity}: 
\begin{enumerate}
\setlength{\parskip}{0cm}
\setlength{\itemsep}{0cm}
\item reformulation of the PDE as a stochastic fixed point problem,
\item approximation of the unique fixed point by Picard iterations,
\item and then approximation of the iterations by multilevel Monte Carlo methods.
\end{enumerate}
The method is said to be full history recursive in the sense that the realizations in the $m$-th iteration require those in the $1$st, $2$nd, $\cdots$, $(m-1)$-th iteration.
Literally, its name originates from the two major components, that is, Picard iterations and multilevel Monte Carlo methods. 
As only a few Picard iterations are required in practice and expectations are nested along the iterations, the conditional expectations may be approximated by the standard Monte Carlo method, leading to nested simulations in the number of iterations. 
The cost can be kept tractable by using only
a small number of samples in each nested layer of simulations, for which a very efficient variance reduction is needed, typically handled by the multilevel Monte Carlo approach.}

A major motivation for introducing multilevel Picard approximation is to solve high dimensional nonlinear PDEs in a similar spirit to deep learning based methods (Section \ref{section deep learning}). 
The category is initiated in the work \cite{German_Picard_pde_code}, in which the applicability of approximation methods based on Picard iterations and multilevel Monte Carlo methods is investigated 
for high dimensional nonlinear PDEs arising in physics and finance. 
It is subsequently shown \cite{doi:10.1098/rspa.2019.0630} that the complexity 
grows polynomially both in the dimension and in the reciprocal of the required accuracy in the case of semilinear heat equations with gradient-independent and globally Lipschitz continuous nonlinearities. 
For instance, its complexity is shown to be $\mathcal{O}(d \epsilon^{-4})$ for such semilinear heat equations, where $d$ is the dimension of the problem and $\varepsilon$ is the required precision.

Multilevel Picard approximation methods have been developed further in many papers. 
A class of multilevel Picard approximation is proposed in \cite{beck2020nonlinear} for computing iterated nested expectations.
\textcolor{black}{For very high-dimensional problems, there is another nested Monte Carlo algorithm for the PDE formulation, where the nesting is along a random time grid \cite{Warin+2018+225+247}}
Multilevel Picard approximation algorithms are shown to overcome the curse of dimensionality for high dimensional nonlinear heat equations with general time horizons and gradient-dependent nonlinearities \cite{overcoming1} and also overcome the curse of dimensionality in the $L^p$ sense for high-dimensional semilinear PDEs  \cite{hutzenthaler2021strong}.
It is shown in \cite{hutzenthaler2021overcoming} that a Monte Carlo type numerical method approximates the solution path of BSDEs with complexity which is at most polynomial in the model dimension and at most quadratic in the reciprocal of the prescribed approximation accuracy.
It is proved in \cite{doi:10.1137/17M1157015} that semilinear heat equations with gradient-dependent nonlinearities can be approximated under suitable technical conditions with polynomial complexity both in the dimension and the reciprocal of the accuracy.
In \cite{hutzenthaler2021speed}, Picard iterations for backward stochastic differential equations with globally Lipschitz continuous nonlinearity are shown to converge, at least, at the rate of square-root of factorial.
An extension of multilevel Picard approximation is applied to general forward diffusion in \cite{https://doi.org/10.48550/arxiv.2204.08511}.
Other numerical algorithms can be found in \cite{E:2021uo} for approximating solutions of general high-dimensional semilinear parabolic partial differential equations at single space-time points, and in \cite{becker2022learning} for parametric approximation problems by combining Monte Carlo algorithms with machine learning techniques to learn the random variables in Monte Carlo simulations including the multilevel Picard approximation.
Yet another type can be found in \cite{chassagneux2021learning} based on a Picard iteration scheme in which a sequence of linear-quadratic optimization problems are solved by means of the stochastic gradient descent algorithm.


\subsection{Further on forward numerical methods}\label{subsection further on forward numerical methods}

We also find a fully forward-looking numerical method in \cite{Fujii:2015aa}, where the non-linear driver of the FBSDE is treated as a perturbation and consequently the FBSDE is converted into a series of decoupled linear FBSDE.
The method proposed for approximating these linear FBSDEs uses an interacting particle method in order to perform the involved integration steps, as opposed to using direct Monte Carlo simulation.
The homotopy analysis method has been reported effective \cite{zhong_liao_2017} for FBSDEs.
The coupled FBSDEs are transformed into a control problem \cite{10.1214/EJP.v10-295} by discretizing the backward component in the forward direction and are simulated by combining the standard backward discretization with Picard iterations \cite{10.1214/07-AAP448}. 
Finally, a forward method can be found in \cite{briand2014} based on a Wiener chaos expansion, Picard iterations and Malliavin calculus techniques.

\section{Deep learning}
\label{section deep learning}

The method of deep learning is one of the most active emerging categories in the context of numerical methods for BSDEs. 
\textcolor{black}{These methods have attracted a great deal of attention for their unique feature of solving high dimensional nonlinear BSDEs and corresponding nonlinear PDEs, 
owning to neural network approximation by which the estimation of nested high-dimensional expectations and gradients is significantly eased.}
So far, a wide variety of numerical results have been developed for approximating the solution of nonlinear PDEs by a neural network from physics and finance and FBSDEs related to a nonlinear pricing model for financial derivatives.
Those numerical methods demonstrate the efficiency and accuracy of the algorithm for several 100-dimensional nonlinear PDEs and FBSDEs.
This effectiveness in high-dimensional problems, in addition to the intrinsic structure involving both aspects of backward and forward methods, makes the method of deep learning quite distinctive in sharp contrast to backward and forward methods (Sections \ref{section backward numerical methods} and \ref{section forward numerical methods}).
\textcolor{black}{Despite the complexity of deep learning methods not being able to be investigated in a general form at present, there exist a few theoretical results in the literature that support deep neural networks in overcoming the curse of dimensionality in numerical approximations of some linear and nonlinear PDEs, where the complexity is shown to grow at most polynomially in both the PDE dimension and the reciprocal of the approximation accuracy (for instance, \cite{BHJK2020,doi:10.1137/19M125649X,10.1093/imanum/drab027,https://doi.org/10.48550/arxiv.1809.02362,hutzenthaler_jentzen_kruse_nguyen_2020,jentzen_salimova_welti_2021,doi:10.1142/S0219530520500116,takahashi_yamada_2022}).}


In what follows, we describe and summarize such deep learning based methods.
To this end, we start with the basics on multilayer neural networks.
We reserve $L$ for the number of layers.
For $k\in\{0,1,\cdots,L-1\}$, let $\ell_k \in \mathbb{N}$ denote the size of input of the $k$-th layer and let $\ell_L$ be the size of the output layer.
For $k\in \{1,\cdots,L\}$, define $W_k:=(w^k_{i,j})_{i,j}$ in $\mathbb{R}^{\ell_{k} \times \ell_{k-1}}$ with $w^k_{i,j}:=\theta^k_{j+\ell_k (i-1)}$ called weight and $b_k:=(\theta^k_{\ell_{k-1} \times \ell_k+1},\cdots,\theta^k_{\ell_{k-1} \times \ell_k+k})$ in $\mathbb{R}^{\ell_k}$ called bias.
Let $\mathbb{A}_{k}: {\bf x}(\in \mathbb{R}^{\ell_{k-1}}) \mapsto W_k {\bf x} + b_k (\in  \mathbb{R}^{\ell_{k}})$ and define $A_e({\bf x}):=(a(x_1),\cdots,a(x_e))$ for ${\bf x} \in \mathbb{R}^e$, where $a$ denotes a real-valued nonlinear function on $\mathbb{R}$, called the activation function. 
With all those, we define the neural network $\mathscr{N}^{\Theta}_L: \mathbb{R}^{\ell_0} \to \mathbb{R}^{\ell_L}$ by $\mathscr{N}^{\Theta}_L({\bf x}):=\mathbb{A}_L \circ \cdots \circ A_{\ell_2} \circ \mathbb{A}_{2} \circ A_{\ell_1} \circ \mathbb{A}_{1}({\bf x})$ for ${\bf x} \in \mathbb{R}^{\ell_0}$, with $\Theta:=(\theta^1_{1},\cdots,\theta^L_{\ell_{L-1} \times \ell_L+L})$.
Note that every continuous function on a compact set can be approximated by a neural network to any desired precision.


\subsection{Deep BSDE}
\label{subsection deep bsde}

We first summarize a numerical scheme called the deep BSDE method \cite{BSDE_DeepLearning, Han8505}.
\textcolor{black}{The essence of deep BSDE methods is, broadly speaking, to make use of the gradient of the solution (with respect to the control process $Z$) as the policy function, 
and approximate it through a neural network as is done in the standard deep reinforcement learning.}
For the reader's convenience, we recall the FBSDE \eqref{eq:FBSDE_Diff} and \eqref{eq:FBSDE_SDE}:  
\begin{align*}
dX_t&=\mu(t,X_t)dt+\sigma(t,X_t)dW_t,\\
-dY_t&=f(t,X_t,Y_t,Z_t)dt-Z_tdW_t,
\end{align*}
with $X_0={\bf x}$ and $Y_T=\Phi(X_T)$, and Theorem \ref{Th:Equivalence_PDE}, which asserts $(Y_t,Z_t)=(v (t,X_t), (\nabla v (t,X_t))^\top \sigma(t,X_t))$, with $v$ satisfying the semilinear parabolic PDE under Assumption \ref{Assumption_3}:
\[
 ((\partial/\partial t)+{\cal L}_t)v(t,{\bf x})+f(t,{\bf x},v(t,{\bf x}), (\nabla v(t,{\bf x}))^\top \sigma(t,{\bf x}) )=0,
\]
with $v(T,{\bf x})=\Phi({\bf x})$.

In general, the deep learning-based method aims to approximate $Y_0=v(0,{\bf x})$. 
We start with the following stochastic optimization problem: 
\[
\inf_{y,Z} \mathbb{E}\left[| \Phi(X_T) - Y_T^{y,Z}|^2\right] \quad \text{subject to} \ \ \
Y_t^{y,Z}:=y-\int_0^t f(s,X_s,Y^{y,Z}_s,Z_s)ds+\int_0^t Z_s dW_s, \quad t\in [0,T].
\]
Evidently, to solve this minimization problem, the expectation $\mathbb{E}[| \Phi(X_T) - Y_T^{y,Z} |^2]$ needs to be computable or approximated.
For this purpose, define the sequence $\{\overline{X}_{t_k}^n\}_{k\in \{0,1,\cdots,n\}}$ as the Euler-Maruyama discretization of the forward process:
\begin{equation}\label{euler forward deep}
\overline{X}_{t_{k+1}}^n=\overline{X}_{t_k}^n+\mu(t_k,\overline{X}_{t_k}^n)(t_{k+1}-t_k)+\sigma(t_k,\overline{X}_{t_k}^n)\left(W_{t_{k+1}}-W_{t_k}\right), \quad k\in \{0,1,\cdots,n-1\},
\end{equation}
with initial state $\overline{X}_{t_0}^n={\bf x}$, and define the {\it forward} discretization $\{\overline{Y}_{t_k}^n\}_{k\in\{0,\cdots,n-1\}}$ of the backward process by
\begin{equation}\label{forward approximation of the backward process}
\overline{Y}_{t_{k+1}}^n=\overline{Y}_{t_k}^n-f(t_k,\overline{X}_{t_k}^n,\overline{Y}_{t_k}^n,z_k(\overline{X}_{t_k}^n))(t_{k+1}-t_k)+z_k(\overline{X}_{t_k}^n)\left(W_{t_{k+1}}-W_{t_k}\right),\quad k\in \{0,1,\cdots,n-1\},
\end{equation}  
with initial state $Y_{t_0}^n=y$, in order to approximate the minimization problem, as follows:
\begin{equation}\label{approximate minimization problem}
\inf_{y,Z} \mathbb{E}\left[| \Phi(X_T) - Y_T^{y,Z} |^2\right] \approx \inf_{y,\{z_k\}_k} \mathbb{E}\left[| \Phi(\overline{X}_T^n)- \overline{Y}^n_T |^2\right]  
\approx \inf_{y,\{ z_k \}_k} \frac{1}{M} \sum_{l=1}^M | \Phi(\overline{X}_{T,l}^n)- \overline{Y}_{T,l}^n |^2, 
\end{equation}
for sufficiently large $M$ for decent estimation quality, where $\{\overline{X}_{T,l}^n\}_{l\in\mathbb{N}}$ and $\{\overline{Y}_{T,l}^n\}_{l\in\mathbb{N}}$ are independent copies of $\overline{X}_T^n$ and $\overline{Y}_T^n$, respectively.
At the end, the obtained solutions, say, $y^*$ and $\{z^*_k\}_{k\in \{0,1,\cdots,n-1\}}$, provide the approximation $Y_0 \approx y^*$ and $Z_{t_k} \approx z^\ast_k(\overline{X}_{t_k}^n)$ for $k\in \{0,1,\cdots,n-1\}$.

Deep learning comes into play when solving the approximate minimization problem \eqref{approximate minimization problem}. 
For $k\in \{1,\cdots,n-1\}$, approximate the target continuous function $z_k$ by a neural network as ${\bf x} \mapsto z_k^{\theta_k}({\bf x}) \approx \mathscr{N}^{\theta_k}_{L}({\bf x})$ with parameter $\theta_k \in \mathbb{R}^{q(L)}$ for a suitable dimension $q(L) \in \mathbb{N}$, depending on the number of layers $L$.
For $\Theta=(\theta^1,\cdots,\theta^{d+1},\theta_1,\cdots,\theta_{n-1})$ with $\theta^k\in \mathbb{R}$ for $k\in \{1,\cdots,d+1\}$, the forward discretization of the backward process \eqref{forward approximation of the backward process} can be parameterized with $\Theta$ as  
\[
\overline{Y}_{t_{k+1}}^n(\Theta)=\overline{Y}_{t_k}^n(\Theta)-f(t_k,\overline{X}_{t_k}^n,\overline{Y}_{t_k}^n(\Theta),z^{\theta_k}_k(\overline{X}_{t_k}^n))(t_{k+1}-t_k)+z^{\theta_k}_k(\overline{X}_{t_k^n})\left(W_{t_{k+1}}-W_{t_k}\right),\quad k\in \{0,1,\cdots,n-1\},
\]  
with initial states $\overline{Y}_0^n(\Theta)=\theta^1$
and $z^{\theta_0}_0({\bf x})=(\theta^{2},\cdots,\theta^{1+d})$. 
We seek the parameter vector $\Theta$ via the following minimization problem:
\[
\inf_{\Theta } \frac{1}{M} \sum_{l=1}^M \left| \Phi(\overline{X}_{T,l}^n)- \overline{Y}_{T,l}^n(\Theta) \right|^2=:\inf_{ \Theta } \phi(\Theta), 
\]
by stochastic gradient descent, that is, $\Theta_{m}=\Theta_{m-1}- \gamma \nabla \phi(\Theta_{m-1})$ recursively for $m \in \mathbb{N}$, where $\gamma$ is a suitable learning rate. 
After a sufficiently large number of recursions, say $K$, we obtain the parameter vector $\Theta_K$ with which the approximation of the target is given by $Y_0=u(0,{\bf x}) \approx \overline{Y}_0^n(\Theta_K)$.

\textcolor{black}{As described in brief in the beginning of the present section, the significance of deep BSDE methods is to adopt deep neural networks for computing gradients of the solution and approximating the backward component forward in time so that high-dimensional PDEs and BSDEs can be solved in a realistic runtime.
An error analysis is conducted for deep BSDE methods in \cite{Han_2020}, where a posteriori estimates are derived for coupled FBSDEs which relate the quadratic terminal loss to the approximation error for the numerical solution of the FBSDE (a posteriori estimates for BSDEs can also be found in \cite{doi:10.1137/120878689}). 
We note that the approximation of the backward component $Y$ parametrically forward in time as above and then the minimization of an error criterion has been found valid even outside the realm of deep learning \cite{ankirchner:hal-01915772}.}


The deep BSDE method is modified in \cite{wang2022deep} by measuring the loss at the forward initial time (rather than the terminal time, as above), called Deep BSDE-ML method, for approximating linear decoupled FBSDEs, as well as is extended in \cite{https://doi.org/10.48550/arxiv.2211.04349} for BSDEs with jumps. 
\textcolor{black}{The backward deep solvers are developed in \cite{https://doi.org/10.48550/arxiv.1807.06622, https://doi.org/10.48550/arxiv.2006.07635} in order to apply the deep BSDE solver-based method to financial problems, such as pricing Bermudan swaption and nonlinear pricing in high-dimensional settings. 
High-dimensional coupled FBSDEs with non-Lipschitz diffusion coefficients are numerically solved in \cite{PUQR2021-025} using the deep BSDE method.} 
The deep BSDE method is also developed in \cite{andersson2022convergence} for strongly coupled FBSDEs stemming from stochastic control.




\subsection{Deep backward dynamic programming}
\label{subsection deep backward dynamic programming}


Deep backward dynamic programming is proposed in \cite{hure2020deep} and then further improved in \cite{germain2021approximation}, literally on the basis of backward dynamic programming arising from discretization methods of BSDEs for high dimensional nonlinear PDEs via the minimization of loss functions at each step, defined recursively by backward induction. 
\textcolor{black}{In contrast to deep BSDE methods (Section \ref{subsection deep bsde}), deep backward dynamic programming is built in reference to the backward resolution technique upon an implicit backward Euler scheme.
It employs machine learning techniques for estimating the solution and its gradient by minimizing a loss function on each time step, 
where such local problems are then solved recursively with a stochastic gradient algorithm backward in time.
Two schemes are proposed in \cite{hure2020deep} for dealing with those local problems: {\bf (DB1)} approximating both the solution and its gradient by a neural network, and {\bf (DB2)} approximating the solution alone by a neural network, with its gradient estimated directly with automatic differentiation.
It is worth mentioning that from a viewpoint of the computation of conditional expectations, deep backward dynamic programming may also be regarded as a backward Euler scheme where nonlinear regression is employed with a deep learning technique.}

To describe deep backward dynamic programming in brief, the solution $(X,Y,Z)$ of the FBSDE \eqref{eq:FBSDE_Diff} and \eqref{eq:FBSDE_SDE}  
is discretized, first through the Euler-Maruyama discretization of the forward component $\{\overline{X}^n_{t_k}\}_{k\in \{0,1,\cdots,n\}}$ with $\overline{X}^n_{t_0}={\bf x}$ on the time grid $(0=)t_0<t_1<\cdots<t_n(=T)$ in a similar manner to \eqref{euler forward deep} and the backward component, by starting with $\overline{Y}^n_{t_n}=\Phi(\overline{X}^n_{t_{n}})$ and then defining
\[
\overline{Y}^n_{t_k}=\overline{Y}^n_{t_{k+1}}+f(t_k,\overline{X}^n_{t_k},\overline{Y}^n_{t_k},\overline{Z}^n_{t_k})(t_{k+1}-t_k)-\overline{Z}^n_{t_k}(W_{t_{k+1}}-W_{t_k}),
\]
for $k\in \{n-1,\cdots,1,0\}$ backwards, where the integrand component $\{\overline{Z}^n_{t_k}\}_{k\in \{0,\cdots,n\}}$ is not explicit as of yet and to be found during the following procedure.
Along the discretized pair $\{(\overline{X}_{t_k}^n,\overline{Y}_{t_k}^n)\}_{k\in \{0,\cdots,n\}}$, we describe the aforementioned two schemes of deep backward dynamic programming: 

\noindent {\bf (DB1)} With the terminal condition ${\cal U}_{n}=\Phi$, iterate 
\[
({\cal U}_k,{\cal Z}_k)=\argmin_{(u_k,z_k)} \mathbb{E}\left[\left| {\cal U}_{k+1}(\overline{X}^n_{t_{k+1}})- \left(u_k(\overline{X}^n_{t_k})+f(t_k,\overline{X}^n_{t_k},u_k(\overline{X}^n_{t_k}),z_k(\overline{X}^n_{t_k}))(t_{k+1}-t_k)-z_k(\overline{X}^n_{t_k})(W_{t_{k+1}}-W_{t_k})\right)\right|^2\right],
\]
for $k\in \{n-1,\cdots,1,0\}$ backwards.

\noindent {\bf (DB2)} With the terminal condition ${\cal V}_n=\Phi$, iterate
\begin{multline*}
{\cal V}_k= \argmin_{u_k} \mathbb{E}\Bigg[\bigg| {\cal V}_{k+1}(\overline{X}^n_{t_{k+1}})\\
- \left(u_k(\overline{X}^n_{t_k})+f(t_k,\overline{X}^n_{t_k},u_k(\overline{X}^n_{t_k}),(\nabla u_k(\overline{X}^n_{t_k}))^{\top}\sigma(t_k,\overline{X}^n_{t_k}))(t_{k+1}-t_k)-(\nabla u_k(\overline{X}^n_{t_k}))^{\top} \sigma(t_k,\overline{X}^n_{t_k}) (W_{t_{k+1}}-W_{t_k}) \right)\bigg|^2\Bigg],
\end{multline*}
for $k\in \{n-1,\cdots,1,0\}$ backwards, where $\nabla u_k$ is the automatic differentiation of the network function $u_k$. 
At each step $k\in \{n-1,\cdots,1,0\}$, the network functions $({\cal U}_k,{\cal Z}_k)$ in {\bf (DB1)} and ${\cal V}_k$ in {\bf (DB2)} can be found literally by employing deep neural networks.
The effectiveness of deep backward dynamic programming has been supported by numerical results on nonlinear PDEs up to 50-dimension \cite{hure2020deep} as well as error analysis \cite{germain2021approximation}.

\subsection{Deep splitting}
\label{subsection deep splitting}

In a similar yet different line from the previous two methods, the so-called deep splitting method is proposed in \cite{doi:10.1137/19M1297919} on the basis of the splitting principle 
with deep neural networks. 
Consider the nonlinear parabolic PDE: 
\[
((\partial/ \partial t) +{\cal L}_t)u(t,{\bf x})+f(t,{\bf x},u(t,{\bf x}),(\nabla u^{\top} \sigma)(t,{\bf x}))=0, \quad (t,{\bf x})\in [0,T)\times \mathbb{R}^q,
\]
with $u(T,{\bf x})=\Phi({\bf x})$, which can be represented, 
by the nonlinear Feynman-Kac formula, as 
\[
u(t,X_t)=\mathbb{E}\left[ \Psi(X_T) +\int_t^T f(s,X_s,u(s,X_s),(\nabla u^{\top} \sigma)(s,X_s)) ds  \bigg|\, {\cal F}_t \right],
\]
where the forward component $\{X_t:\, t\in [0,T]\}$ is as given in Section \ref{subsection deep bsde}.
In light of this nonlinear Feynman-Kac representation, one wishes to construct the recursive approximation $V_k(\cdot)\approx u(t_k,\cdot)$ for $k\in \{n-1,\cdots,1,0\}$ backwards at discrete time points $(0=)t_0<t_1<\cdots<t_n(=T)$, given by
\[
V_k(\overline{X}^n_{t_k})=\mathbb{E}\left[V_{k+1}(\overline{X}^n_{t_{k+1}})+f(t_{k+1},\overline{X}^n_{t_{k+1}},V_{k+1}(\overline{X}^n_{t_{k+1}}),(\nabla V_{k+1}(\overline{X}^n_{t_{k+1}}))^{\top}\sigma(t_{k+1},\overline{X}^n_{t_{k+1}})) (t_{k+1}-t_k) \Big|\, \overline{X}^n_{t_k} \right], 
\]
where $\{\overline{X}^n_{t_k}\}_{k\in \{0,1,\cdots,n\}}$ is a sequence of discrete observations of the forward component by the Euler-Maruyama scheme with $\overline{X}^n_{t_0}={\bf x}$ in a similar manner to \eqref{euler forward deep}.
In the deep splitting method, this recursion is approximated by a continuous function on the support of the marginal via the minimization: 
\[
V_k\approx \argmin_{v \in \mathcal{C}({\rm supp}(\overline{X}^n_{t_k})) } \mathbb{E}\left[\left|v(\overline{X}^n_{t_k})-\left( V_{k+1}(\overline{X}^n_{t_{k+1}})+f(t_{k+1},\overline{X}^n_{t_{k+1}},V_{k+1}(\overline{X}^n_{t_{k+1}}),(\nabla V_{k+1}(\overline{X}^n_{t_{k+1}}))^{\top}\sigma(t_{k+1},\overline{X}^n_{t_{k+1}})) (t_{k+1}-t_k) \right) \right|^2 \right],
\]
so as to approximate the map ${\bf x} \mapsto V_k({\bf x})$ by a neural network $V_k(\cdot) \approx {\cal N}_L^\Theta(\cdot)$ with respect to parameter set $\Theta$. 
It is reported in \cite{doi:10.1137/19M1297919} that the deep splitting method succeeds to deal with as high-dimensional nonlinear PDEs as $10000$ dimensions.

\subsection{Deep Galerkin method and physics-informed neural networks}\label{subsection DGM PINN}

Deep Galerkin Method (DGM) \cite{DGM} and Physics-Informed Neural Networks (PINN) \cite{PINN} are proposed for solving high-dimensional nonlinear problems.
Despite both methods being based largely on deterministic PDEs, we summarize their essence here for the reason that those may provide solvers for nonlinear BSDEs as well.
Hereafter, we refer to DGM and PINN collectively as the deep learning based PDE solver and give a brief summary in accordance with \cite{DGM}. 

Let $D \subset \mathbb{R}^d$ and let $u: [0,T] \times D \to \mathbb{R}$ solve the following (nonlinear) PDE:   
\[
\begin{cases}
(\partial/\partial t) u(t,{\bf x}) +f(u(t,{\bf x}),\nabla u(t,{\bf x}),\mbox{Hess}(u(t,{\bf x})))=0,& (t,{\bf x}) \in [0,T] \times D, \\
u(0,{\bf x})=\Psi({\bf x}),& {\bf x} \in D,\\
u(t,{\bf x})=g(t,{\bf x}),& (t,{\bf x}) \in [0,T] \times \partial D, 
\end{cases}
\]
where $f$ is a nonlinear differential operator, and $\Psi: \mathbb{R}^d \to \mathbb{R}$ and $g: [0,T] \times \mathbb{R}^d \to \mathbb{R}$ represent the initial and boundary conditions, respectively.
In the deep learning based PDE solver, the solution $u$ is approximated by a function $u^{\Theta}={\cal N}_L^{\Theta}$ through a deep neural network with respect to parameter set $\Theta$, by minimizing the error between both sides of the nonlinear PDE above, evaluated at sampled points in the space-time domain and in its boundary.
An optimal parameter set $\Theta^*$ is here searched for in such a way to ideally minimize the loss function: 
\begin{equation}\label{loss function Ltheta}
\ell(\Theta):=\left\|(\partial/\partial t) u^{\Theta} +f(u^{\Theta},\nabla u^{\Theta},\mbox{Hess}(u^\Theta))\right\|_{[0,T] \times D,\nu_1}
+\left\| u^{\Theta}(0,\cdot)-\Psi(\cdot) \right\|_{D,\nu_2}+\left\| u^{\Theta} -g\right\|_{[0,T] \times \partial D,\nu_3},
\end{equation}
where $\| \varphi \|_{S,\nu}:=\int_S | \varphi({\bf x}) |^2\nu(d{\bf x})$ with the probability measure $\nu$ and its support $S$.
For minimization, one employs stochastic gradient decent on an approximate loss function:
\[
\Theta_{n+1}=\Theta_n-\gamma_n \nabla G_n(\Theta_n),
\]
where $\{\gamma_n\}_{n\in\mathbb{N}}$ is a sequence of learning rates.
Here, $G_n$ is an approximation of the loss function \eqref{loss function Ltheta}, defined by
\[
G_n(\Theta):=\left|(\partial/\partial t) u^{\Theta}(t_n,x_n) +f(u^{\Theta}(t_n,x_n),\nabla u^{\Theta}(t_n,x_n),\mbox{Hess}(u^{\Theta}(t_n,x_n)))\right|^2\\
+\left| u^{\Theta}(0,w_n)-\Psi(w_n)\right|^2+\left| u^{\Theta}(\tau_n,z_n) -g(\tau_n,z_n)\right|^2,
\]
where $(t_n, x_n)$, $w_n$ and $(\tau_n,z_n)$ are random elements, respectively, taking values in $[0, T] \times D$, in $D$ and in $[0,T] \times \partial D$, according to the probability measures $\nu_1$, $\nu_2$ and $\nu_3$.
The iterative procedure above is to be repeated until a suitable convergence criterion is met.

In addition, more generalized frameworks using Monte Carlo methods and an efficient implementation of the neural network are discussed in \cite{DGM}, along with various numerical results on high-dimensional American options, high-dimensional HJB equations and Burger's equations.
We do not go into further details and applications on the deep learning based PDE solver, since those would lie way outside the scope of the present survey on BSDEs.
Instead, we refer the reader to \cite{Nature, SIAMReview} for more applications and 
to \cite{nusken2021interpolating, raissi2018forwardbackward} for its close connection with BSDEs. 

\subsection{Further on deep learning based methods}

A wide variety of deep learning based schemes have been developed for solving high-dimensional BSDEs, some of which take advantage of the aforementioned backward and forward methods (Sections \ref{section backward numerical methods} and \ref{section forward numerical methods}).
For instance, in \cite{fujii_takahashi_takahashi_2019}, a deep learning scheme is introduced in combination with asymptotic expansion (Section \ref{subsection asymptotic expansion}), which is further extended in  \cite{doi:10.1142/S2424786320500127, takahashi2021new, Tsuchida:2022vo}.
The deep learning technique of \cite{BSDE_DeepLearning} is extended in \cite{ganesan2020pricing, yu2019deeplearning} to address both terminal and boundary conditions of PDEs.
As a computational framework for portfolio risk management problems, the so-called Deep xVA solver is proposed \cite{gnoatto2020deep} by recursively using the deep BSDE method for a coupled system of BSDEs.
A discretization scheme is employed in \cite{negyesi2021step} for solving BSDEs based on deep learning regressions (Section \ref{subsection LS regression}).
A deep signature/log-signature FBSDE algorithm is developed in \cite{feng2021deep} for approximating FBSDEs with state and path dependent features. 
A deep Runge-Kutta method is proposed in  \cite{https://doi.org/10.48550/arxiv.2212.14372}.
In the presence of constraints on the gains process, an approximation of BSDEs is obtained in \cite{KharroubiLimWarin} by neural network approximation.
Iterative diffusion optimization techniques are studied using deep learning techniques \cite{nusken_richter} for applications such as importance sampling and rare event simulation.
The Long Short Term Memory networks is applied in \cite{kapllani2021deep} to improve the Deep BSDE method \cite{BSDE_DeepLearning}.
A new algorithm is proposed in \cite{teng2021gradient} based on a $\theta$-discretization of the time-integrands with eXtreme Gradient Boosting (XGBoost) regression for efficiently computing conditional expectations. 
A deep learning-based stochastic branching algorithm is developed in \cite{nguwi_penent_privault_2022} for numerically solving fully nonlinear PDEs.
High-dimensional fully-coupled FBSDEs are solved in \cite{ji2020algorithms} with three algorithms based on deep learning.
A deep learning based method is proposed in \cite{teng_shi_zhu_2021} for solving forward-backward doubly stochastic differential equations.
Finally, we refer the reader to  \cite{BHJK2020, https://doi.org/10.1002/gamm.202100006, E_2021_nonlinearity, germain2021neural} for more recent developments and surveys on deep learning based methods, to \cite{chassagneux2021numerical} for singular BSDEs, to \cite{Beck:2019aa} for second-order BSDEs (Section \ref{section 2bsde}), to \cite{carmona2021convergence, 10.3389/fams.2020.00011, germain2021numerical, https://doi.org/10.48550/arxiv.2204.11924} for McKean-Vlasov BSDEs (Section \ref{subsection mckean-vlasov}), and to \cite{ji2022deep, pereira2022decentralized} for 
stochastic control problems.


\section{Discussion}\label{section discussion}

In the preceding sections (Sections \ref{section backward numerical methods}, \ref{section computation of conditional expectations}, \ref{section forward numerical methods} and \ref{section deep learning}), we have presented a systematic survey and categorization of various numerical methods for BSDEs.
However, these methods have thus far been presented in isolation and so not in an appropriate manner for drawing comparisons.
In the present section, we thus present those categories all on the table, with a brief description of each, so that a relevant contrast can be made against a few key factors regarding their implementation.
Such factors include the convergence property, the dimensionality, and the complexity, which directly influence, individually and/or in combination, what method would be chosen for the BSDE in question.
For instance, if a very accurate approximation is needed, on the one hand, then it is desirable to select a method with a strong theoretical convergence guarantee even if the method is more difficult to code and has a longer running time.
On the other hand, if a less rigorous approximation is enough for the time being (such as, for preliminary testing purposes), then it would be more reasonable to choose a method that is easier to code and runs faster, even if it may be less theoretically accurate.

In order to make comprehensive, yet, to-the-point contrasts, we must, out of sheer necessity, restrict ourselves to a representative method or two in displaying each category and do not claim that the resulting conclusions entirely hold for every numerical method in a certain category.
Moreover, throughout, we focus on the approximation of the backward component $Y$ of each BSDE model $(X,Y,Z)$ and the corresponding PDE solution (either the point $u(t,{\bf x})$ or the function $u(t,\cdot)$).
We reserve $d$ for the dimension of the forward component $X$ and assume $Y$ is one-dimensional, that is, ${\bf x}(\in \mathbb{R}^d ) \mapsto u(t,{\bf x}) (\in \mathbb{R})$, and write $C$ for constant multiples whose values change depending on the context.

\subsection{Backward numerical methods along with computation of conditional expectations}
\label{subsection discussion backward methods}

As backward discretization (Section \ref{section backward numerical methods}) is intrinsically built upon the computation of conditional expectations (Section \ref{section computation of conditional expectations}), we discuss representative methods consisting of those two methodological components in combination.
Here, we let $n$ represent the number of discretization points in the given time interval, that is, $(t_0,t_1,\cdots,t_n)$ and let $M$ indicate the number of iid replications for Monte Carlo methods involved in each subinterval. 

\begin{itemize}
\item LSMC (Section \ref{subsection LS regression}) + Backward Euler (Section \ref{subsection backward euler methods}): 
Let $Y_{t_k}^{n,M,K}$ denote an approximation at time $t_k$ in accordance with Section \ref{subsection LS regression}, and let $n$ be the number of time discretization, 
instead of the mesh $\pi$.
It is shown under suitable conditions and choice of basis functions \cite{Lin_regress1} that 
\begin{equation}\label{rate_lsm}
\max_{k\in \{0,\cdots,n\}} \left\|Y_{t_k}-Y_{t_k}^{n,M,K} \right\|_2 = \mathcal{O} \left( n^{-1/2} \right),
\end{equation}
provided that $K \approx n^d$ and $M \approx n^{d+3}$.
Despite the availability of theoretical analyses \cite{https://doi.org/10.48550/arxiv.0806.4447, backward_discrete_5, Lin_regress1}, it is still an open question as to the choice of basis functions for a given BSDE.
It is worth mentioning that LSMC methods are effective in approximating the function  $u(t,\cdot)$ (rather than a single point $u(t,{\bf x})$ alone). 

\item Malliavin (Section \ref{subsection malliavin calculus based methods}) + Backward Euler (Section \ref{subsection backward euler methods}):  
Let $Y_{t_k}^{n,M}$ denote an approximation at time $t_k$ in accordance with 
\cite{backward_discrete_4}.
It is shown under suitable conditions that 
\begin{equation}
\max_{k\in \{0,\cdots,n\}} \left\|Y_{t_k}-Y_{t_k}^{n,M} \right\|_2 \leq C n^{-1/2} + C \frac{n^{d/4}}{M^{1/2}},\label{rate_malliavin}
\end{equation} 
which suggests $M \approx n^{d/2+1}$ to bound the overall rate by $\mathcal{O}(n^{-1/2})$.
In general, the required computation can be heavy for dealing with terms involving Skorohod integrals (such as \eqref{malliavin fraction}),
for which an improved algorithm is proposed in \cite{Mall_improve}, while as many as $2^d$ iid standard normal random variables need to be generated.

\item Quantization (Section \ref{subsection quantization methods}) + Backward Euler (Section \ref{subsection backward euler methods}):
Let $Y_{t_k}^{n,N}$ denote an approximation at time $t_k$ in accordance with \cite{Quantiz2}, where $N=1+N_1+\cdots+N_n$ with each $N_k$ the number of points in $\mathbb{R}^d$ used to make up the space grid at the $k$-th discretization step.
It is shown under suitable conditions that 
\begin{equation}
\max_{k\in\{0,\cdots,n\}} \left\|Y_{t_k}-Y_{t_k}^{n,N}\right\|_2 \le C n^{-1/2}+ C\frac{n^{1+1/d}}{N^{1/d}},  \label{rate_quantization}
\end{equation}
which suggests $N \approx n^{(3/2)d+1}$ to bound the overall rate by $\mathcal{O}(n^{-1/2})$.
As described in the above and Section \ref{subsection quantization methods}, Monte Carlo simulation is required for completing the quantization method with care on the number of iid replications $M$.
Quantization methods work in relatively low dimensions, while generally being effective for the approximation of the function $u(t,\cdot)$.

\item Tree (Section \ref{subsection tree based methods}) + Backward Euler (Section \ref{subsection backward euler methods}): 
Let $Y_{t_k}^{n}$ be an approximation at time $t_k$ by the standard (deterministic) binomial tree  method \cite{ECP1030, Indep_Z_Discrete}. 
In order to satisfy the order $1/2$, that is,
\begin{equation}
\max_{k\in \{0,\cdots,n\}} \left\|Y_{t_k}-Y_{t_k}^{n}\right\|_2 = \mathcal{O}\left(n^{-1/2} \right), \label{rate_tree}
\end{equation}
the required complexity is approximately $2^{dn}$, which is very large, whereas the complexity drops down to $(n+1)$ with a recombination tree in one-dimensional problems.

\item Cubature (Section \ref{subsection cubature methods}) + Second-order discretization (Section \ref{subsection higher-order methods}):  
Let $Y_0^{n,N}$ be an approximation of the initial state $Y_0$ by a cubature method and the tree-based branching algorithm (TBBA) measure with the cubature degree of $7$ and $N$ particles at every time point $t_k$, while $n$ here denotes the number of time discretization by the second order method  \cite{Crisan2014}. 
It is known under suitable conditions that 
\[
\mathbb{E}\left[\left|Y_0-Y_0^{n,N}\right|\right] \le  C n^{-2}+ C \frac{n}{N^{1/2}}.
\]
We refer the reader to \cite{chassagneux_trillos_2020} for relevant complexity analysis.

\end{itemize}

With numerical methods for conditional expectations (Section \ref{section computation of conditional expectations}) now collectively aligned in conjunction with backward discretization methods (Section \ref{section backward numerical methods}), we make a few relevant remarks before moving on to forward numerical methods (Section \ref{subsection discussion forward methods}).
First, recall that the backward Euler method $Y^n$ with forward Euler-Maruyama scheme $X^n$ is shown to be of order $1/2$ for the backward process $Y$ of a Markovian FBSDE $(X,Y,Z)$ under minimal Lipschitz conditions \cite{backward_discrete_3}, that is, $\max_k \| Y_{t_k}-Y_{t_k}^n \|_2 \leq C n^{-1/2}$, whereas the better rate $\| Y_{0}-Y_{0}^n \| \leq C n^{-1}$ at the initial time has also been reported elsewhere \cite{backward_discrete_1, backward_discrete_6}.
This is not an essential disagreement but a simple difference in the evaluation point, due to the error on the forward component of order $1/2$ for all time points but the initial time (that is, $\| X_{t_k}-X_{t_k}^n \|_2=\mathcal{O}(n^{-1/2})$ for all $k\in \{1,\cdots,n\}$) and  the deterministic equality at the initial time (that is, $X_0=X_0^n={\bf x}$). 
Thus, the rates described in \eqref{rate_lsm}, \eqref{rate_malliavin}, \eqref{rate_quantization} and \eqref{rate_tree} can be conservative if the approximation is focused on the initial state $Y_0$ alone.

Next, we make a note on the complexity and dimensionality, which are undoubtedly important  upon implementation individually, as well as inextricably bound up together.
LSMC can work in up to 10 dimensions and tends to be more efficient than Malliavin calculus based methods.
For example, it is known \cite{Lin_regress1} that the complexity of LSMC is $\mathcal{O}(\epsilon^{-d-4})$ to achieve the squared error $\epsilon^2$ in some instances, while the Malliavin calculus based method bears the complexity $\mathcal{O}(\epsilon^{-d-13})$, provided that the forward component is a geometric Brownian motion \cite{backward_discrete_4}.
In general, the application of quantization or tree based methods is limited to considerably low-dimensional problems due to the use of the safety grid. 
Cubature methods with multi-linear interpolation result in the complexities $\mathcal{O}(\epsilon^{-3d/2+1/2})$ if the backward Euler discretization is employed and $\mathcal{O}(\epsilon^{-d+1/2})$ with the second-order discretization with the Richardson-Romberg extrapolation \cite{chassagneux_trillos_2020}.
As such, in most cases, the complexity of backward methods (in combination with computation for conditional expectations) tends to explode exponentially as the dimension of the problem increases.

We also remark on the standing assumption of each category. 
Least-squares regression based methods generally work under minimal assumptions on BSDEs (and SDEs), while a few extra regularity assumptions on the coefficients are required to allow the derivation of robust estimates for the involved error using various regression tools and also to ensure the stability of the algorithm (for instance, \cite{chassagneux2014AAP, LSMDP_noMal, pelsser2019, TENG2020117}).
Malliavin calculus based methods tend to require lots of extra assumptions, including various regularity assumptions and in particular, assumptions regarding the Malliavin weights (as in \cite{Mall_LSregres, Indep_Z_LSMC, naito_yamada_2019}).
Quantization and tree based methods work under the standard Lipschiz conditions on forward and backward SDEs, while cubature methods require some extra smoothness on coefficients of the forward SDE and a smooth driver with Lipschitz terminal condition in most instances.

\subsection{Forward numerical methods} 
\label{subsection discussion forward methods}

Next, we make a contrast of the four categories of forward numerical methods (Section \ref{section forward numerical methods}), again, by providing a brief summary of the convergence property, 
and key features.

\begin{itemize}
\item Picard iteration (Section \ref{subsection picard iteration methods}): 
For an approximate solution $(Y^{r,n}_t)_{t\in [0,T]}$ in accordance with the Picard iteration method \cite{Picard_it1}, it is shown that  
\[
   \sup_{t\in [0,T]} \mathbb{E}\left[ \left|Y_{t} - Y_t^{r,n}\right|^2 \right] \leq C n^{-1}+C \left(\frac{1}{2} + Cn^{-1} \right)^r,
\]
for sufficiently large $n$.
This scheme needs to be combined with a numerical method for computing conditional expectations (Section \ref{section computation of conditional expectations}) and, by and large, similar limitations apply as described in Section \ref{subsection discussion backward methods}.
As such, in its implementation, variance reduction techniques play an importance role, for instance, importance sampling \cite{importance_samp1, FBSDE_Importance} and control variates \cite{BSDE_AdapCont}, in solving high-dimensional problems.

\item Branching diffusion system (Section \ref{subsection branching diffusion methods}): 
For a BSDE $(Y,Z)$ with driver $f$ (independent of ${\bf z}$), it is shown under suitable conditions \cite{BSDE_branching_1} that 
\[
 \left|Y_{0} - Y_0^{n,\ell}\right| \le C \left\| f-f^{n,\ell} \right\| + Cn^{-1/2}, 
\]
where $Y_{\cdot}^{n,\ell}$ denotes a branching diffusion approximation to the BSDE with a discretized $\ell$-th order polynomial driver $f^{n,\ell}$ \cite{BSDE_branching, Indep_Z_branching}.
After the branching diffusion approximation has been applied, it is the usual practice to employ Monte Carlo simulation for implementation, such as a nonlinear Monte Carlo simulation \cite{10.1007/978-3-319-33446-2_3}, which is known to be effective in high-dimensional nonlinear pricing problems.

\item Asymptotic expansion (Section \ref{subsection asymptotic expansion}):
For a BSDE $(Y^\varepsilon,Z^\varepsilon)$ with a perturbed driver $\varepsilon f$, it holds under suitable conditions \cite{doi:10.1137/14100021X,Malliavin_expansions} that 
\[
\sup_{t\in [0,T]}\mathbb{E}\left[\left| Y_t^\epsilon-\left( Y_t^{(0)}+ \epsilon Y_t^{(1)}+\cdots+({\epsilon^m}/{m!}) Y_t^{(m)}\right) \right|^2 \right]^{1/2}  = \mathcal{O}(\epsilon^{m+1}),
\]
where $Y_t^{(k)}=((d^k/d \epsilon^k)Y_t^{\epsilon})|_{\epsilon=0}$ for $k\in \{0,1,\cdots\}$.
After the expansion has been performed, similarly to the aforementioned branching diffusion system, often Monte Carlo simulation is used for implementation \cite{10.1007/978-3-319-33446-2_3}.

\item Multilevel Picard approximation (Section \ref{subsection multilevel picard}): 
Consider a solution $u$ of a $d$-dimensional semilinear PDE 
corresponding to a BSDE model, 
and let $u^{\theta}_{n,M}:[0,T]\times \mathbb{R}^d \times \Omega \to \mathbb{R}$ denote an approximate solution to the PDE $u(\cdot,\cdot)$ by the multilevel Picard approximation in accordance with \cite{E_2021_nonlinearity,doi:10.1098/rspa.2019.0630}.
For each ${\bf x}\in \mathbb{R}^d$, there exist $K:(0,1] \to \mathbb{N}$ and $c>0$ such that for all $d \in \mathbb{N}$ and $\epsilon \in (0,1]$, 
\[
\mathbb{E}\left[\left|u(T,{\bf x})-u^0_{K(\varepsilon),K(\varepsilon)}(T,{\bf x})\right|^2\right]^{1/2} \leq \epsilon,
\]
with the complexity $\mathcal{O}(\varepsilon^{-c}d^{c})$, which suggests that for each initial state of the forward component, the initial state $Y_0$ of the backward component can be approximated without suffering from the curse of dimensionality.

\end{itemize}

To conclude the categories of forward numerical methods, we make a short remark on the standing assumption, as we have done towards the end of Section \ref{subsection discussion backward methods}.
As mentioned throughout the survey, forward numerical methods do not inherently work backwards to avoid the computation of conditional expectations, often with the aid of the FBSDE and PDE equivalence (Theorem \ref{Th:Equivalence_PDE}).
In return, it is rather evident that numerical methods based on solving the equivalent PDE (Section \ref{section forward numerical methods}) must impose Assumption \ref{Assumption_3} so that the FBSDE and PDE equivalence holds. 
In addition, a few extra regularity assumptions are often made.
For instance, the branching diffusion system based method (as of \cite{BSDE_branching, Indep_Z_branching, BSDE_branching_1}) often requires the key assumption that the driver can be represented as the sum of a power series to be approximated by polynomials, which obviously does not hold for every BSDE. 
The asymptotic expansion method (as of \cite{crepey_song, doi:10.1142/S0219024912500343, Fujii_Takahashi_particle_2015, Malliavin_expansions})  
requires some smoothness on the driver to expand a nonlinear BSDE by linear BSDEs, while the condition can be relaxed in some instances \cite{doi:10.1137/14100021X}.
Finally, to employ the multilevel Picard iteration method, one needs to carefully verify the relevant conditions on the driver as well as check the structure of the forward component (see, for example, \cite{overcoming1, doi:10.1098/rspa.2019.0630, doi:10.1137/17M1157015, hutzenthaler2021speed, https://doi.org/10.48550/arxiv.2204.08511}).

\subsection{Summary}

To conclude the comparisons and discussions made above, and by considering the enormous numerical examples on BSDEs in the literature, we summarize in Table \ref{table 0} key aspects of each category upon implementation, that is, what is approximated and the problem dimension.
In particular, the dimension of the problem at hand is one of the more restrictive conditions for which method can be chosen, in particular, most methods are not efficient in high-dimensional problems.
Note that we did not make a comparative discussion on deep learning based methods (Section \ref{section deep learning}) in the form of a separate subsection (like Sections \ref{subsection discussion backward methods} and \ref{subsection discussion forward methods}), as their error and complexity analyses are still evolving, some in infancy, and awaiting major advances, as opposed to the effectiveness in very high-dimensional spatial approximation problems.
Still, for the sake of completeness, we align the category of deep learning based methods at the bottom of Table \ref{table 0} with a very broad perspective.

\begin{table}[H]
\centering{
\begin{tabular}{|c||c|c|} 
\hline
methods & target & dimension\\ \hline\hline
LSMC + Backward Euler & function $u(t,\cdot)$ & $d\approx 10$\\ \hline
Malliavin + Backward Euler & function $u(t,\cdot)$ & a few \\ \hline
Quantization + Backward Euler & function $u(t,\cdot)$ & a few \\ \hline
Tree + Backward Euler & function $u(t,\cdot)$ & one or two \\ \hline
Cubature + Second-order discretization & function $u(t,\cdot)$ & a few \\ \hline \hline
Picard iteration & point $u(t,{\bf x})$ & higher than a few \\ \hline
Branching diffusion system & point $u(t,{\bf x})$ & higher than a few \\ \hline
Asymptotic expansion & point $u(t,{\bf x})$ & higher than a few \\ \hline
Multilevel Picard approximation & point $u(t,{\bf x})$ & $d>100$  \\ \hline \hline
Deep learning & function $u(t,\cdot)$ &  $d>100$ \\ \hline
\end{tabular}
\caption{Two key aspects upon implementation}
\label{table 0}
}\end{table}

Without a doubt, computing time is a vital factor to consider when selecting a method, as there may be a need for very quick calculations.
However, if no time pressure exists, then the choice of method is less restrictive. 
We remark that various attempts have been made so far in the literature by wisely splitting the required computation for massive parallelization on highly multicore GPUs.
Parallelization here has naturally proven effective because the primary issue does not lie in the computational time but memory consumption requirements, for instance, by algorithms that require many sample paths at once on memory, such as binomial lattice based methods \cite{5662503, 6128469, Parallel_FD}, the multistep method \cite{kapllani2019multistep}, LSMC \cite{gobet_lopezsalas_vasquez_2020, backward_discrete_7, proceedings2019021044}, the four-step scheme \cite{parallel_fourstep}, and the forward Picard iteration \cite{LabartLelong+2013+11+39}.

\color{black}
\section{Numerical methods for BSDEs with nonstandard features}\label{section BSDEs with nonstandard features}

To address problems in stochastic control, finance, and partial differential equations, BSDEs often need to be equipped with extra features for better capturing relevant properties under consideration.
There has been an increasing interest in those classes of BSDEs and in developing numerical methods exclusively for each or some combinations of those features.
Despite most of what follows having already appeared in their respective subsections, we again categorize and summarize those numerical methods here in terms of the class of BSDEs in brief, without going into much detail in order to avoid overloading the paper with lengthy technical intricacies.

\subsection{Coupled FBSDEs}
\label{subsection coupled FBSDEs}

Consider the following FBSDE, where the coefficients of the forward SDE can depend on the backward components $(Y,Z)$:
\[ 
 \begin{cases} 
 X_t = {\bf x} + \int_0^t \mu(s,X_s,Y_s,Z_s)ds + \int_0^t\sigma(s,X_s,Y_s)dW_s,\\
 Y_t = \Phi(X_T) + \int_t^T f(s, X_s, Y_s,Z_s)ds- \int_t^T Z_sdW_s,
 \end{cases}
\]
for $t\in [0,T]$.
The FBSDE in this form is called a coupled FBSDE in the sense that the backward components are allowed to couple in the drift and the diffusion of the forward component, in sharp contrast to the decoupled formulation \eqref{eq:FBSDE_Diff}-\eqref{eq:FBSDE_SDE}.
Hence, the standard Ito generator $\mathcal{L}_t$ and the diffusion coefficient $\sigma$ in the PDE \eqref{eq:PDE_semilinear} for decoupled FBSDEs needs to replaced accordingly in order to carry the backward components, resulting in the following PDE under suitable technical conditions:
\[ 
    \begin{cases} 
      (\partial/\partial t)v(t,{ \bf x})+ \langle \mu(t,{ \bf x},v(t,{\bf x}), (\nabla v(t,{\bf x}))^{\top}\sigma(t,{\bf x},v(t,{\bf x}))), \nabla v(t,{\bf x})\rangle \\
      \quad +\frac{1}{2}{\rm tr}[\sigma^{\otimes 2}(t,{\bf x},v(t,{\bf x})){\rm Hess}(v(t,{\bf x}))]
         + f(t, {\bf x}, v(t,{\bf x}), (\nabla v(t,{\bf x}))^{\top}\sigma(t,{\bf x},v(t,{\bf x}))) = 0, & (t,{\bf x})\in [0,T)\times \mathbb{R}^q,\\
        v(T,{\bf x}) = \Phi({\bf x}),& {\bf x}\in \mathbb{R}^q.
    \end{cases}
\]
As such, numerical approximation of a coupled FBSDE or the equivalent PDE is evidently far more intricate than that of a decouple counterpart, of which the forward component can be treated separately from the backward components.  
In particular, if the forward component depends on $Z$ as well, then the gradient of the solution $\nabla u(t,{\bf x})$ needs to be addressed inside the forward component.
Such coupled FBSDEs naturally appear in the utility maximization problem with general utility function in economics \cite{BSDE_utility1}.


As for numerical methods, the four-step scheme \cite{Four_step1} opened the door for coupled FBSDEs as early as in the 1990s, followed by various studies 
through the corresponding semi- or quasi-linear parabolic PDE. 
Based on this four-step scheme, the authors in \cite{Finite_diff2} develop an implementable numerical method.
Specifically, for the PDE part, they use the combined characteristics and finite difference method to approximate its solution.
Almost sure uniform convergence and weak convergence of the method are proven, with the rate of convergence being comparable to that of the approximation of the forward SDE (done using an Euler type scheme).
In \cite{delarue2006} and moreover \cite{delarue2008}, we find improved alternatives of the method described in \cite{Finite_diff2}, which weakens regularity assumptions required.
The authors in \cite{Finite_diff3, Milsteinbook2004, 10.1093/imanum/drl019} also give a numerical method based on the four-step scheme, but propose a different method for solving the involved PDE. They approximate the solution of the PDE by using layer methods, which are constructed by means of a probabilistic approach. The derivatives of the solution to the PDE are found by using finite differences. 
The four-step scheme once again appears and is reconstructed in \cite{parallel_fourstep} with new conditions. It is then associated with the idea of domain decomposition methods (associated with the Schawrz waveform relaxation method). This approach is used in order to parallelize the related equations.
Finally, in \cite{GONG2015220}, we find a numerical method based on the four-step scheme which focuses on approximating the solution to coupled FBSDEs.

We add that various numerical methods have also been applied in tailored ways to coupled FBSDEs, such as Picard iterations \cite{10.1214/07-AAP448}, a transform into a control problem via a fully forward discretization \cite{10.1214/EJP.v10-295}, multistep schemes \cite{liu2020fully, doi:10.1137/130941274}, a defferred correction method for ODEs \cite{tang_zhao_zhou_2017}, Fourier methods for computing conditional expectations \cite{HUIJSKENS2016593}, and more recently, deep learning based methods \cite{andersson2022convergence, Han_2020, ji2020algorithms, ji2022deep}.
 








\subsection{Reflected BSDEs}
\label{subsection reflected BSDEs numerical methods}

We first review reflected backward stochastic differential equations (RBSDEs) \cite{BSDE_American}, which form an important class of BSDEs in the sense that they provide a deep insight into many practical problems, such as optimal stopping problems, American option pricing, stochastic optimal controls and differential games.
Their solutions are reflected at a given stochastic process, called the obstacle, and provide a probabilistic representation for the unique viscosity solution of an obstacle problem for a nonlinear parabolic partial differential equation.

Now, consider the following BSDE, where the process $(Y_t)_{t\in [0,T]}$ is forced to stay above the process $(h(t,X_t))_{t\in [0,T]}$: 
\[ 
 \begin{cases} 
         Y_t = h(T,X_T) + \int_t^T f(s, X_s, Y_s,Z_s)ds +K_T-K_t- \int_t^T Z_sdW_s, \\
         Y_t\ge h(t,X_t),\\
         \int_0^T(Y_s-h(s,X_s))dK_s=0,\\
          X_t = {\bf x}+ \int_0^t \mu(X_s)ds + \int_0^t\sigma(X_s)dW_s,
    \end{cases}
\]
for $t\in [0,T]$, where the stochastic process $(K_t)_{t\in [0,T]}$ is continuous and non-decreasing with $K_0=0$, and grows only when $Y_t=h(t,X_t)$.
The process $(Y_t)_{t\in [0,T]}$ can then be written in the form of the solution to an optimal stopping problem, as
\[
 Y_t=\esssup_{\tau\in \mathcal{T}_t}\mathbb{E}\left[\int_t^{\tau} f(s,X_s,Y_s)ds +h(\tau,X_{\tau})\Big|\,\mathcal{F}_t\right], \quad t\in [0,T],
\]
where $\mathcal{T}_t$ denotes the set of $(\mathcal{F}_s)_{s\in [t,T]}$-stopping times taking values in $[t,T]$, and provides a probabilistic interpretation $Y_t=v(t,X_t)$ where $v$ is a solution to the obstacle problem
\[
 \begin{cases}
  \max\left\{((\partial/\partial t)+\mathcal{L}_t)v(t,{\bf x})+f(t,{\bf x},v(t,{\bf x})),\,h(t,{\bf x})-v(t,{\bf x})\right\}=0,&(t,{\bf x})\in [0,T)\times \mathbb{R}^q,\\
 v(T,{\bf x})=h(T,{\bf x}),& {\bf x}\in \mathbb{R}^q, 
 \end{cases}
\]
in a similar yet different manner from Theorem \ref{Th:Equivalence_PDE}.
For theoretical details on RBSDEs, we refer the reader to, for instance, the monograph \cite[Chapter 14]{delong}.    
    
The most popular application of RBSDEs is, perhaps, the pricing of American options \cite{BSDE_American, BSDE_American2}, due to its structure with lower obstacles.
In fact, the RBSDE framework is not strictly necessary for the pricing of American options, as the classical binomial pricing model for example works well for its low-dimensional problems, while least-square Monte-Carlo methods \cite{LSMC_forward} remain efficient for high-dimensional problems, both outside the RBSDE framework.
Nonetheless, the RBSDE framework is useful for the pricing of American options, for instance, through detailed error analysis by approximating RBSDEs in quantization methods \cite{Quantiz2, Quantiz1} and its stability results for approximating its payoff function \cite{10.1007/978-3-642-25746-9_7}. 


As for discretization methods, it is proved in \cite{Quantiz2} that the backward Euler type discretization attains the order $1/2$ for BSDEs with Lipchitz coefficients and the order $1$ for those with semi-convex obstacles. 
Other backward Euler methods are developed for RBSDEs in \cite{Katarzyna, doi:10.1080/07362994.2011.610162} based on the random walk approximation, in \cite{BSDE_Discrete} based on the lattice tree and in \cite{backward_discrete_4} in combination with LSMC regression and the Malliavin approach (Sections \ref{subsection LS regression} and \ref{subsection malliavin calculus based methods}) for approximating conditional expectations. 
The stability analysis for backward Euler methods for RBSDEs is conducted in \cite{bouchard2008discrete, ma2005representations}.  
Most recently, the stability and convergence analysis is refined in \cite{sun2020quantitative} with a focus on quadratic RBSDEs.
 
After discretizing RBSDEs, more complex equations need to be solved concurrently with approximation of conditional expectations in the presence of the obstacle term.
This difficulty has been addressed by discretizing the state space, such as the so-called quantization, before employing the dynamic programming principle. 
Quantization algorithms for RBSDEs, initially developed in \cite{BALLYPAGESPRINTEMS+2001+21+34}, have been investigated intensively in \cite{Quantiz2, Quantiz1} to deal with multi-dimensional RBSDEs 
with an application to the pricing of an American option, and 
further extended, for example, in combination with the higher-order discretization \cite{mcwalter2018recursive} and for those whose driver depends on the control process \cite{nmeir2021quantizationbased}. 
In addition, the primal-dual method is generalized in \cite{German_primal_dual} to a backward dynamic programming equation associated with time discretization schemes for RBSDEs.

In the literature, other advanced types of RBSDEs have also attracted attention, such as discrete RBSDEs, RBSDEs with multiple obstacles and RBSDEs with jumps. 
To deal with those classes, the concept of penalization, which was originally employed for proving the existence and uniqueness of RBSDEs \cite{BSDE_American}, proves useful for numerical discretization \cite{Memin:2008aa} as well, and can be combined with Fourier transform techniques \cite{Hyndman:2017aa}.
It has been applied to optimal switching problems in real option pricing \cite{hamadene2007starting}, followed by RBSDEs with oblique reflections \cite{chassagneux2012, hu2010multi}, doubly RBSDEs \cite{XU20111137}, and doubly RBSDEs with jumps \cite{DUMITRESCU2016827}. 
A numerical method is developed in \cite{risks8030072} for doubly RBSDEs whose generator depends on the future values of the solution with application to the default risk modeling.




\subsection{BSDEs with jumps}
\label{subsection BSDEs with jumps}

BSDEs may contain jumps in their backward component as $-dY_t=f(t,X_t,Y_t,Z_t,U_t)-Z_tdW_t-U_tdN_t$ with a suitable jump process $(N_t)_{t\in [0,T]}$ and an control process $(U_t)_{t\in [0,T]}$, and often in their forward component as well.
On top of non-trivial issues on existence, uniqueness and stability \cite{papapantoleon2021stability}, the presence of jumps increases complexity of numerical implementation, for instance, by paying a great deal of attention to the jump size and timing.  
For instance, discretization schemes are developed, based on Malliavin calculus on the Gaussian space \cite{Mall_Calc1}, and the so-called shot noise representation of the Poisson random measure \cite{10.1214/20-PS359} in \cite{massing2021approximation}.
A forward scheme for BSDEs \cite{briand2014} is extended \cite{GEISS20162123} to BSDEs with jumps by Wiener chaos expansion and Picard iterations.
A numerical algorithm is developed in \cite{Fancy_BSDEs} to approximate the solution to a decoupled FBSDE driven by a pure L\'evy process with its small jumps approximated by a suitable Brownian motion.
An explicit prediction-correction scheme is developed for solving decoupled FBSDEs with jumps \cite{EAJAM-6-253}.
In \cite{lejay:inria-00357992}, Brownian and Poissonian components are independently approximated by random walks.
More complex structures are also of interest.
A Fourier-based method \cite{doi:10.1137/16M1099005} and asymptotic expansion \cite{doi:10.1080/17442508.2018.1521808} are proposed, by which non-Poissonian local jump measures can be taken care of.
Numerical methods for doubly reflected BSDEs with jumps are developed, with \cite{DUMITRESCU2016206} and without \cite{DUMITRESCU2016827} penalization. 
A single jump (not systematic jumps) is considered in the context of BSDEs in \cite{KharroubiLim+2015+81+109}, where a single jump can represent a default time in credit risk or counterparty risk.
A deep learning based method is developed in \cite{https://doi.org/10.48550/arxiv.2211.04349} and applied to BSDEs with jumps.
We refer the reader to the monograph \cite{delong} for details as well as applications in insurance and finance.

\subsection{BSDEs with non-global Lipschitz conditions}
\label{subsection Lipschitz BSDEs}

The majority of key properties of numerical methods for BSDEs are derived upon the global Lipschitz conditions on the drivers and terminal condition, whereas such restrictive conditions do not hold all the time in important applications, such as the Fisher-KPP and FitzHugh-Nagumo equations.
In the literature, considerable attempts have been made to develop viable numerical methods where those Lipschitz conditions are relaxed. 
The $L^2$-time regularity of the $Z$-component is studied in \cite{GOBET20101105} for an irregular terminal condition (such as an indicator function), as an extension of the $L^2$-time regularity result \cite{backward_discrete_3} for a Lipschitz terminal, where the order of convergence is explicitly connected to the rate of decrease of the expected conditional variance of the terminal.
Two discretization methods are studied in \cite{turkedjiev2015} for a class of locally Lipschitz Markovian BSDEs:
One is a backward Euler scheme (Section \ref{subsection backward euler methods}) by approximating a projection of the $Z$ process, while the other is the Malliavin weight scheme (Section \ref{subsection malliavin calculus based methods}, especially \eqref{Z_formula} and \eqref{Mall_weights}) by directly approximating the marginals of the $Z$ process, where advanced a priori estimates and stability results for this class of BSDEs are derived by extending the representation theorem \cite{Mall_rep} and employed for obtaining competitive convergence rates.
It is proved \cite{cheridito_stadje_2013} that a backward stochastic difference equation (BS$\Delta$E) admits a solution and its sequence is convergent to a BSDE as the time-grid gets finer even when their drivers are not Lipschitz in the $Z$ component.
An error analysis of a time discretization method is performed in \cite{lionnet2015} for systems of BSDEs with drivers of polynomial growth and monotone in the state variable, where a tamed version of the explicit Euler scheme is shown to converge despite that the standard version may diverge unlike with the canonical Lipschitz driver.
A general framework is developed in \cite{lionnet2018} for explicit numerical schemes for BSDEs with drivers of polynomial growth, where the convergence of some modified explicit scheme is of the same rate as implicit schemes and has comparable computing cost to the standard explicit scheme.
Finally, for drivers not Lipschitz in the backward component, the Picard iteration is applied to discretized FBSDEs \cite{borkowski_janczakborkowska_2022} and is employed for deriving the existence and uniqueness for $G$-BSDEs \cite{zhang_jiang_2021}.

\subsection{Quadratic BSDEs}
\label{subsection qBSDEs}

BSDEs with generators of quadratic growth (with respect to the variable ${\bf z}$ as in $|f(t,{\bf x},{\bf y},{\bf z})|\le C(1+\|{\bf y}\|+\|{\bf z}\|^2)$) have been studied actively under the name of quadratic BSDEs.
Quadratic BSDEs play important roles in mathematical finance, such as utility optimization with exponential utility functions and risk minimization for the entropic risk measure.
Their numerical aspects have thus been investigated ever since the theoretical formalization \cite{10.1214/aop/1019160253}.
In contrast to standard BSDEs, the treatment of the quadratic driver is not trivial, particularly in placing the upper bound estimate of the control process $Z$.
In the development of numerical methods, this issue is typically addressed by truncating the quadratic driver.

The well-known path regularity theorem \cite{backward_discrete_3} is extended in \cite{IMKELLER2010348} to the quadratic-growth setting along with convergence rate of the distance between the solution of a quadratic BSDE and its approximation by truncation.
In \cite{Imkeller_DosReis_Zhang_2010}, the Cole-Hopf exponential transformation is applied to the approximate solution. 
Time-discretization schemes are developed in \cite{richou2011} for quadratic BSDEs based on truncation and non-uniform time partition, and in \cite{chassagneux2016} with the error of order less than $1/2$.
For quadratic BSDEs with reflection, a truncated discrete-time numerical scheme is proposed in \cite{sun2020quantitative} along with some practical examples.
In \cite{doi:10.1142/S2010139212500152}, an asymptotic expansion technique (Section \ref{subsection asymptotic expansion}) is developed for a quadratic-growth FBSDE appearing in an incomplete market with stochastic volatility.
An approximation method is constructed in \cite{FUJII20191492} for quadratic BSDEs using semi-analytic asymptotic expansion. 
With a view towards utility maximization problems, a special class of backward stochastic partial differential equations (Section \ref{subsection BSPDEs}) is investigated in \cite{doi:10.1142/S0219024911006437}, along with numerical simulation for relevant portfolio optimization problems.  

\subsection{Second-order BSDEs}\label{section 2bsde}

Second-order BSDEs (2BSDEs) are a type of BSDE whose nonlinear drift contains the second order derivative of the corresponding PDE, widely applied in financial modeling such as the uncertain volatility model, transaction cost model, illiquid market model and the pricing for passport options. 
In general, a process $(Y_t,Z_t,\Gamma_t,\alpha_t)_{t\geq 0}$ 
is called a 2BSDE if it solves the system
\begin{align*}
-dY_t=\left[f(X_t,Y_t,Z_t,\Gamma_t)-\frac{1}{2}\mbox{tr}\left[\sigma^{\otimes 2}(t,X_t)\Gamma_t\right]\right]dt-Z_t \sigma(t,X_t)dW_t,\quad
dZ_t=\alpha_t dt+\Gamma_t\sigma(t,X_t)dW_t,\quad t\in [0,T],
\end{align*}
with $Y_T=g(X_T)$ and a suitable forward process $X$ in the form of \eqref{eq:FBSDE_SDE}.
This system can be expressed as $Y_t=u(t,X_t)$, $Z_t=\nabla u(t,X_t)$, $\Gamma_t=\mbox{Hess}(u(t,X_t))$ and $\alpha_t=((\partial/\partial t)+{\cal L})\nabla u(t,X_t)$ for $t\in [0,T]$ on the basis of the parabolic PDE:
\[
\frac{\partial}{\partial t} u(t,x)+f(t,x,u(t,x),\nabla u(t,x),\mbox{Hess}(u(t,x)))=0,\quad u(T,x)=g(x),
\]
whose full nonlinearity makes the class of 2BSDE distinct from the standard BSDEs. 
Its existence and the uniqueness are proved in \cite{cheridito2007second} under Lipschitz continuity of the driver and suitable conditions, along with the numerical scheme: starting with $\overline{Y}_T=g(\overline{X}_T)$ and $\overline{Z}_T=\nabla g(\overline{X}_T)$, proceed backwards for $k\in \{1,\cdots,n\}$, 
\begin{gather*}
\overline{Y}_{t_{k-1}}=\mathbb{E}[ \overline{Y}_{t_k} | \mathscr{F}_{t_{k-1}} ]+\left[f(t_{k-1},\overline{X}_{t_{k-1}},\overline{Y}_{t_{k-1}},\overline{Z}_{t_{k-1}},\overline{\Gamma}_{t_{k-1}})
-\frac{1}{2}\mbox{tr}\left[\sigma^{\otimes 2}(t_{k-1},\overline{X}_{t_{k-1}})\overline{\Gamma}_{t_{k-1}}\right]  \right](t_k-t_{k-1}),\\
\overline{Z}_{t_{k-1}}=\frac{1}{t_k-t_{k-1}} (\sigma^{-1}(t_{k-1},\overline{X}_{t_{k-1}}))^{\top} \mathbb{E}\left[ \overline{Y}_{t_k} (W_{t_k}-W_{t_{k-1}}) | \mathscr{F}_{t_{k-1}} \right],\quad \overline{\Gamma}_{t_{k-1}}=\frac{1}{t_k-t_{k-1}} \mathbb{E}\left[ \overline{Z}_{t_k} (W_{t_k}-W_{t_{k-1}})^{\top} | \mathscr{F}_{t_{k-1}} \right]\sigma^{-1}(t_{k-1},\overline{X}_{t_{k-1}}),
\end{gather*}
based on the aforementioned nonlinear Feynman-Kac representation.
In \cite{10.1214/10-AAP723}, 
an implementable explicit method is developed to approximate the conditional expectations appearing in the above scheme.
A variance reduction method is proposed in \cite{Alanko_Avellaneda} for the computation of involving conditional expectations.
\textcolor{black}{2BSDEs can also be effectively solved numerically with a monotone scheme \cite{10.1214/14-AAP1030}.}
We mention that deep learning based approximation algorithms (Section \ref{section deep learning}) are developed in \cite{Beck:2019aa} for 2BSDEs.

\subsection{McKean-Vlasov FBSDEs}
\label{subsection mckean-vlasov}

When dealing with a stochastic differential game with mean field interactions, also known as a mean field game, one often encounters the following optimal control problem:
\begin{equation}\label{mean field game}
\inf_{\alpha \in {\cal A}} \mathbb{E}\left[\int_t^T f(s,X_s,\mu_s,\alpha_s)ds + g(X_T,\mu_T)\Big|\,X_t={\bf x}\right]=:u(t,{\bf x}),
\end{equation}
subject to the so-called McKean-Vlasov stochastic differential equation
\begin{equation}\label{MKV SDE}
 dX_s=b(s,X_s,\mu_s,\alpha_s)ds+\sigma(s,X_s,\mu_s,\alpha_s)dW_s,
\end{equation}
where $(\mu_t)_{t\in [0,T]}$ is a deterministic flow of probability measures and the infimum in \eqref{mean field game} is taken over the set $\mathcal{A}$ of all progressively measurable processes.
The McKean-Vlasov SDE \eqref{MKV SDE} can be interpreted as a limit of the empirical measure of increasing individual particles, known as the propagation of chaos \cite{chaintron2021propagation} where the particles tend to be independent of each other in this limit as the impact of each decreases.
Under suitable technical conditions, the solution $u$ to the optimization problem \eqref{mean field game} corresponds to that of a suitable HJB equation, as well as providing an insightful probabilistic interpretation through the representation $u(t,X_t)=Y_t$ in accordance with the McKean-Vlasov FBSDE, in general: for $s\in [t,T]$,
\begin{align*}
dX_s&=b(s,X_s,{\cal L}(X_s),{\cal L}(Y_s),{\cal L}(Z_s))ds+\sigma(s,X_s,{\cal L}(X_s),{\cal L}(Y_s),{\cal L}(Z_s)) dW_s, \\
-dY_s&=f(s,X_s,{\cal L}(X_s),{\cal L}(Y_s),{\cal L}(Z_s))ds-Z_sdW_s, 
\end{align*} 
with $X_t=\xi$ and $Y_T=g(X_T,{\cal L}(X_T))$, where ${\cal L}(F)$ denotes the law of the random element $F$.
We note that two approaches (of Pontryagin and weak types) result in distinct McKean-Vlasov FBSDEs and thus numerical methods \cite{esaimprocs}.

As for relevant numerical methods, a variety of deterministic ones can be found in the literature for mean field games \cite{Achdou2020, lauriere2021numerical}, for instance, based on finite differences or variational approaches, whereas McKean-Vlasov FBSDEs above pave the way towards probabilistic numerical methods and analysis \cite{esaimprocs}, such as a tree method \cite{10.1214/18-AAP1429}, first-order and Crank-Nicolson schemes \cite{Zhang:2022vq}, a higher order discretization method for decoupled cases \cite{de2015cubature}, McKean-Vlasov anticipative FBSDEs arising in initial margin requirements \cite{Agarwal_et_al_2019} and a posteriori error estimates for fully coupled McKean-Vlasov FBSDEs \cite{reisinger2020posteriori}, and recently very actively, deep learning based methods \cite{campbell2021deep, carmona2021convergence, 10.3389/fams.2020.00011, germain2021numerical, https://doi.org/10.48550/arxiv.2204.11924}.

\subsection{BSPDEs}
\label{subsection BSPDEs}

Backward stochastic partial differential equations are of emerging interest from an implementation point view in its intimate relation to stochastic optimal and utility maximization problems.
A few existing developments of numerical methods are as follows.
For forward-backward stochastic heat equations from stochastic optimal control, discretization methods are developed via spatial discretization \cite{doi:10.1137/15M1022951}, space-time discretization \cite{prohl_wang_2021}, and the local discontinuous Galerkin method \cite{doi:10.1137/20M1383252}.
A special class of BSPDEs is investigated in \cite{doi:10.1142/S0219024911006437} via their reduction to BSDEs from both theoretical and numerical perspectives in the context of utility maximization and portfolio optimization problems.
A deep learning based method is developed in \cite{teng_shi_zhu_2021} for forward-backward doubly stochastic differential equations, where the additional diffusion terms in its BSDE component yields an equivalence to semilinear parabolic SPDEs as well as optimal control problems.





\section{Concluding remarks}\label{section concluding remarks}

In this survey, we have focused on examining a wide array of numerical methods for BSDEs along with the main assumptions made, key convergence properties, as well as summaries of the key advantages and disadvantages.
In particular, we have characterized broadly into three categories; backward (Section \ref{section backward numerical methods}), forward (Section \ref{section forward numerical methods}) and deep learning based methods (Section \ref{section deep learning}), along with the computation of conditional expectations (Section \ref{section computation of conditional expectations}) to complete the categories of backward numerical methods (Section \ref{section backward numerical methods}).
Within those categories, we have further categorized the involved methods.
For backward methods, we have sub-categorized the methods as backward Euler methods (Section \ref{subsection backward euler methods}) and higher-order methods (Section \ref{subsection higher-order methods}).
For the computation of conditional expectations to complement backward methods (Section \ref{section backward numerical methods}), we have sub-categorized the methods as least-squares regression based methods (Section \ref{subsection LS regression}), Malliavin calculus based methods (Section \ref{subsection malliavin calculus based methods}), Malliavin weights dynamic programming with regression methods (Section \ref{subsection dynamic programming}), quantization methods (Section \ref{subsection quantization methods}), tree based methods (Section \ref{subsection tree based methods}), and cubature methods (Section \ref{subsection cubature methods}). 
For forward methods, we have sub-categorized the methods as Picard iteration methods (Section \ref{subsection picard iteration methods}), branching diffusion system based methods (Section \ref{subsection branching diffusion methods}), asymptotic expansion (Section \ref{subsection asymptotic expansion}), and multilevel Picard iterations (Section \ref{subsection multilevel picard}).
For deep learning based methods, we have sub-categorized the methods as deep BSDE (Section \ref{subsection deep bsde}), deep backward dynamic programming (Section \ref{subsection deep backward dynamic programming}), deep splitting (Section \ref{subsection deep splitting}), and deep Galerkin methods and physics-informed neural networks (Section \ref{subsection DGM PINN}).

As is clear from the present survey, numerical methods for BSDEs are intrinsically involved, and further, have been presented in a scattered manner in the literature.
Hence, for the sake of collective comparison, we have devoted Section \ref{section discussion} to displaying those surveyed categories against a few important practical aspects, such as convergence properties and implementation aspects, each of which is crucial for understanding the power and limitations of numerical methods from a practical viewpoint.

In addition, BSDEs of complex structure, and thus the numerical methods employed for them, have become of growing interest in the years to come.
To enhance this line of research, we have also surveyed existing numerical methods in terms of BSDEs with nonstandard features in Section \ref{section BSDEs with nonstandard features}, such as coupled (Section \ref{subsection coupled FBSDEs}), reflection (Section \ref{subsection reflected BSDEs numerical methods}), jumps (Section \ref{subsection BSDEs with jumps}), non-Lipschtiz (Section \ref{subsection Lipschitz BSDEs}), quadratic (Section \ref{subsection qBSDEs}), second-order (Section \ref{section 2bsde}) and McKean-Vlasov BSDEs (Section \ref{subsection mckean-vlasov}), and BSPDEs (Section \ref{subsection BSPDEs}).
To the best of our knowledge, numerical aspects have not been explored for other types as of yet, such as ergodic BSDEs, delayed BSDEs, and BSDEs with regime-switching.
As such, we hope this survey can serve as an initial yet insightful guideline when selecting and developing appropriate numerical methods for the interested BSDE, now and in the future.

With the advent of deep learning based methods (Section \ref{section deep learning}), advances in numerical methods for BSDEs appear to have entered a new era. 
In sharp contrast to the conventional approaches (Sections \ref{section backward numerical methods}, \ref{section computation of conditional expectations}, and \ref{section forward numerical methods}), deep learning based methods have opened the door wide for very high-dimensional BSDEs and nonlinear PDEs by employing the policy functions for readily approximating gradient dependent nonlinearities, to which the existing methods have devoted substantial effort.
Further lines of numerical methods are still expected to come, which will benefit from being built upon the distinctive developments and advantageous features of the existing numerical methods for BSDEs and PDEs, to say nothing of expanding fields of application in relation to high-dimensional and nonlinear problems.
However, it also remains to address a variety of unexplored theoretical aspects of the existing numerical methods, for instance, ones for BSDEs with nonstandard features (Section \ref{section BSDEs with nonstandard features}), not only from the traditional perspective of applied mathematics and numerical analysis, but also from the viewpoint of pure mathematics.
Synergy and complementarity of those distinctive expertise would be beneficial and indispensable for making further advances on numerical methods for BSDEs.

\small
\bibliographystyle{abbrv}
\bibliography{BSDEsurvey_ver3.bib}

\end{document}